\documentclass[10t,a4paper]{article}                     % onecolumn (standard format)
\usepackage{color}   %May be necessary if you want to color links
\usepackage[bookmarks=true,pdfencoding=unicode]{hyperref}
\hypersetup{
    backref=true,
    citecolor=blue,
    colorlinks=true, %set true if you want colored links
    linktoc=all,     %set to all if you want both sections and subsections linked
    linkcolor=blue,  %choose some color if you want links to stand out
}

\usepackage{bookmark}
\usepackage[latin1]{inputenc}
\usepackage[T1]{fontenc}
\usepackage{dsfont}
\usepackage{mathptmx}      % use Times fonts if available on your TeX system
\usepackage{amsmath}
\usepackage{amssymb}
\usepackage{subfig}
\usepackage{mdframed}
\usepackage{authblk}

\newcommand{\qed}{\hfill \ensuremath{\Box}}
\newcommand{\eod}{\hfill \ensuremath{\rule[1.0pt]{4.4pt}{4.4pt}}}
\newcommand{\ewl}{\hfill \ensuremath{\blacktriangle}}

\DeclareMathAlphabet{\mathcal}{OMS}{cmsy}{m}{n}

\newtoggle{LatexWithProofs}
\togglefalse{LatexWithProofs}

\newtoggle{LatexFull}
\toggletrue{LatexFull}
\newtheorem{pLemma}{Lemma}
\newtheorem{pProp}{Proposition}
\newtheorem{pArgument}{Argument}
\newtheorem{pExample}{Example}
\newtheorem{pCorol}{Corollary}
\newtheorem{pLongDef}{Definition}
\newtheorem{pTheorem}{Theorem}

\newcommand{\pRef}[1]{(\ref{#1})}

\newcommand{\peLongDef}[1]{$\eod$ \end{pLongDef}}
\newcommand{\peLongDefX}[1]{\end{pLongDef}}

\newcommand{\peArgument}[1]{ \hyperlink{#1}{\ewl} \end{pArgument}}

\newcommand{\peExample}[1]{ \hyperlink{#1}{\ewl} \end{pExample}}
\newcommand{\peExampleX}[1]{\end{pExample}}

\newcommand{\pePlainExample}[1]{$\eod$ \end{pExample}}
\newcommand{\pePlainExampleX}[1]{\end{pExample}}

\newcommand{\peCorol}[1]{\hyperlink{#1}{\ewl} \end{pCorol}}
\newcommand{\peCorolX}[1]{\end{pCorol}}

\newcommand{\peLemma}[1]{ \hyperlink{#1}{\ewl} \end{pLemma}}
\newcommand{\peLemmaX}[1]{\end{pLemma}}

\newcommand{\peProp}[1]{ \hyperlink{#1}{\ewl} \end{pProp}}
\newcommand{\pePropX}[1]{\end{pProp}}

\newcommand{\peTheorem}[1]{ \hyperlink{#1}{\ewl} \end{pTheorem}}

\newcommand{\pLink}[1]{\eqno \hyperlink{#1}{\ewl}}

\iftoggle{LatexFull}
{
\newcommand{\pFullLink}[1]{\eqno \hyperlink{#1}{\ewl}}
\newcommand{\peFullProp}[1]{ \hyperlink{#1}{\ewl} \end{pProp}}
\newcommand{\peFullCorol}[1]{ \hyperlink{#1}{$\ewl$} \end{pCorol}}
\newcommand{\peFullLemma}[1]{ \hyperlink{#1}{$\ewl$} \end{pLemma}}
\newcommand{\peFullExample}[1]{ \hyperlink{#1}{$\ewl$} \end{pExample}}
}
{
\newcommand{\pFullLink}[1]{\eqno \eod}
\newcommand{\peFullProp}[1]{$\eod$ \end{pProp}}
\newcommand{\peFullCorol}[1]{$\eod$ \end{pCorol}}
\newcommand{\peFullLemma}[1]{$\eod$ \end{pLemma}}
\newcommand{\peFullExample}[1]{$\eod$ \end{pExample}}
}

\iftoggle{LatexWithProofs}
{
\definecolor{darkpink}{rgb}{0.91, 0.33, 0.5}
\definecolor{darksalmon}{rgb}{0.91, 0.59, 0.48}
\definecolor{desertsand}{rgb}{0.93, 0.79, 0.69}
\definecolor{celadon}{rgb}{0.67, 0.88, 0.69}
\definecolor{darkcyan}{rgb}{0.0, 0.55, 0.55}

\newcommand{\pbClaim}[1]{{\color{red} ?}}
\newcommand{\peClaim}[1]{{\color{red} ?`}}
\newcommand{\pbArgument}[1]{\begin{pArgument} \label{#1}  {\color{darksalmon} #1}}

\newcommand{\pbExample}[1]{\begin{pExample} \label{#1}  {\color{darksalmon} #1}}

\newcommand{\pbExampleB}[1]{
\begin{pExample} \label{#1}
{\color{darksalmon} #1}
\hypertarget{{T#1}}{}
\bookmark[
rellevel=1,
keeplevel,
dest=T#1
]{Example \ref{#1}}
}

\newcommand{\pbExampleBT}[2]{
\begin{pExample}[#2] \label{#1}
{\color{darksalmon} #1}
\hypertarget{{T#1}}{}
\bookmark[
rellevel=1,
keeplevel,
dest=T#1
]{Example \ref{#1}: {#2}}
}

\newcommand{\pbCorol}[1]{\begin{pCorol} \label{#1}  {\color{darksalmon} #1}}

\newcommand{\pbCorolB}[1]{
\begin{pCorol} \label{#1}  {\color{darksalmon} #1}
\hypertarget{{T#1}}{}
\bookmark[
rellevel=1,
keeplevel,
dest=T#1
]{Corollary \ref{#1}}
}

\newcommand{\pbCorolBT}[2]{
\begin{pCorol}[#2]\label{#1}  {\color{darksalmon} #1}
\hypertarget{{T#1}}{}
\bookmark[
rellevel=1,
keeplevel,
dest=T#1
]{Cor. \ref{#1}: {#2}}
}

\newcommand{\pbLemma}[1]{\begin{pLemma} \label{#1}  {\color{darksalmon} #1}}

\newcommand{\pbLemmaB}[1]{
\begin{pLemma} \label{#1}  {\color{darksalmon} #1}
\hypertarget{{T#1}}{}
\bookmark[
rellevel=1,
keeplevel,
dest=T#1
]{Lemma \ref{#1}}
}

\newcommand{\pbLemmaBT}[2]{
\begin{pLemma}[#2] \label{#1}  {\color{darksalmon} #1}
\hypertarget{{T#1}}{}
\bookmark[
rellevel=1,
keeplevel,
dest=T#1
]{Lem. \ref{#1}: {#2}}
}

\newcommand{\pbProp}[1]{\begin{pProp} \label{#1}  {\color{darksalmon} #1}}

\newcommand{\pbPropB}[1]{
\begin{pProp} \label{#1}  {\color{darksalmon} #1}
\hypertarget{{T#1}}{}
\bookmark[
rellevel=1,
keeplevel,
dest=T#1
]{Prop. \ref{#1}}
}

\newcommand{\pbPropBT}[2]{
\begin{pProp}[#2] \label{#1}  {\color{darksalmon} #1}
\hypertarget{{T#1}}{}
\bookmark[
rellevel=1,
keeplevel,
dest=T#1
]{Prop. \ref{#1}: {#2}}
}

\newcommand{\pbTheorem}[1]{\begin{pTheorem} \label{#1}  {\color{darksalmon} #1}}

\newcommand{\pbTheoremB}[1]{
\begin{pTheorem} \label{#1}  {\color{darksalmon} #1}
\hypertarget{{T#1}}{}
\bookmark[
rellevel=1,
keeplevel,
dest=T#1
]{Theorem \ref{#1}}
}

\newcommand{\pbTheoremBT}[2]{
\begin{pTheorem}[#2] \label{#1}  {\color{darksalmon} #1}
\hypertarget{{T#1}}{}
\bookmark[
rellevel=1,
keeplevel,
dest=T#1
]{Theor. \ref{#1}: {#2}}
}

\newcommand{\pbLongDef}[1]{\begin{pLongDef}
\label{#1} \hypertarget{{#1}}{} {\color{darksalmon} #1}
\bookmark[
rellevel=1,
keeplevel,
dest=#1
]{Definition \ref{#1}}
}

\newcommand{\pbLongDefB}[2]{\begin{pLongDef}
\label{#1} \hypertarget{{#1}}{} {\color{darksalmon} #1}
\bookmark[
rellevel=1,
keeplevel,
dest=#1
]{\ref{#1}: #2}
}

\newcommand{\pbLongDefBT}[2]{\begin{pLongDef}[#2]
\label{#1} \hypertarget{{#1}}{} {\color{darksalmon} #1}
\bookmark[
rellevel=1,
keeplevel,
dest=#1
]{\ref{#1}: #2}
}

\newcommand{\pbChain}[1]{\[ }
\newcommand{\peChain}[1]{\ {\color{red} \checkmark  \wrm{#1}} \] }

\newcommand{\pePClaim}[1]{ \checkmark}

\newcommand{\pbTClaim}[1]{\begin{equation} \label{#1}  }
\newcommand{\peTClaim}[1]{ \ {\color{red} \checkmark  \wrm{#1}} \end{equation} }
\newcommand{\peEClaim}[1]{ \ {\color{red} \checkmark  \wrm{#1 \bullet} }\end{equation} }

\newcommand{\pbDef}[1] {\hypertarget{{#1}}{} \begin{equation} \label{#1} }
\newcommand{\peDef}[1]{ {\color{blue} \bigstar \wrm{#1}} \end{equation}}

\newcommand{\pbHypot}[1]{ \begin{equation} \label{#1}  }
\newcommand{\peHypot}[1]{  \ \ {{\color{blue} \bigstar \wrm{#1}}}  \end{equation} }

\newcommand{\pbProofB}[2]{\newpage {\color{darksalmon} Proof of {#1} {#2}} \ref{#2} \hypertarget{{#2}}{}
\bookmark[
rellevel=1,
keeplevel,
dest=#2
]
{{#1} \ref{#2}}
}

\newcommand{\peProof}[2]{\qed{} {\color{darksalmon} End of proof of #1 #2} \ref{#2}}

\newcommand{\pbVerifyB}[1]{\newpage {\color{darksalmon} Verification of Example {#1}} \ref{#1} \hypertarget{{#1}}{}
\bookmark[
rellevel=1,
keeplevel,
dest=#1
]
{Example \ref{#1}}
}

\newcommand{\peVerify}[1]{\qed{} {\color{darksalmon} End of verification of Example #1} \ref{#1}}

}
{
\newcommand{\pbChain}[1]{\[}
\newcommand{\peChain}[1]{\]}

\newcommand{\pbClaim}[1]{}
\newcommand{\peClaim}[1]{}
\newcommand{\peTClaim}[1]{ \end{equation}}
\newcommand{\peEClaim}[1]{ \end{equation}}
\newcommand{\pbTClaim}[1]{\begin{equation} \label{#1}}

\newcommand{\peDef}[1]{  \end{equation}}

\newcommand{\lpPlainClaim}[1]{}
\newcommand{\lpEndPlainClaim}[1]{}
\newcommand{\pbDef}[1]{ \hypertarget{{#1}}{} \begin{equation} \label{#1}}

\newcommand{\pbHypot}[1]{\begin{equation} \label{#1}  }
\newcommand{\peHypot}[1]{ \end{equation}}

\newcommand{\pbCorol}[1]{\begin{pCorol} \label{#1}}

\newcommand{\pbCorolB}[1]{\begin{pCorol} \label{#1}
\hypertarget{{T#1}}{}
\bookmark[
rellevel=1,
keeplevel,
dest=T#1
]
{Corollary \ref{#1}}
}

\newcommand{\pbCorolBT}[2]{
\begin{pCorol}[#2] \label{#1}
\hypertarget{{T#1}}{}
\bookmark[
rellevel=1,
keeplevel,
dest=T#1
]
{Cor. \ref{#1}: {#2}}
}

\newcommand{\pbArgument}[1]{\begin{pArgument} \label{#1}}
\newcommand{\pbExample}[1]{\begin{pExample} \label{#1}}

\newcommand{\pbExampleB}[1]{\begin{pExample} \label{#1}
\hypertarget{{T#1}}{}
\bookmark[
rellevel=1,
keeplevel,
dest=T#1
]{Example \ref{#1}}
}

\newcommand{\pbExampleBT}[2]{
\begin{pExample}[#2] \label{#1}
\hypertarget{{T#1}}{}
\bookmark[
rellevel=1,
keeplevel,
dest=T#1
]{Example \ref{#1}: {#2}}
}

\newcommand{\pbLemma}[1]{\begin{pLemma} \label{#1}}
\newcommand{\pbLemmaBT}[2]{\begin{pLemma}[#2] \label{#1}
\hypertarget{{T#1}}{}
\bookmark[
rellevel=1,
keeplevel,
dest=T#1
]
{Lem. \ref{#1}: {#2}}
}

\newcommand{\pbLemmaB}[1]{\begin{pLemma} \label{#1}
\hypertarget{{T#1}}{}
\bookmark[
rellevel=1,
keeplevel,
dest=T#1
]
{Lemma \ref{#1}}
}

\newcommand{\pbProp}[1]{\begin{pProp} \label{#1}}
\newcommand{\pbPropBT}[2]{\begin{pProp}[#2] \label{#1}
\hypertarget{{T#1}}{}
\bookmark[
rellevel=1,
keeplevel,
dest=T#1
]
{Prop. \ref{#1}: {#2}}
}

\newcommand{\pbPropB}[1]{\begin{pProp} \label{#1}
\hypertarget{{T#1}}{}
\bookmark[
rellevel=1,
keeplevel,
dest=T#1
]
{Prop. \ref{#1}}
}

\newcommand{\pbTheorem}[1]{\begin{pTheorem} \label{#1}}
\newcommand{\pbLongDef}[1]{\begin{pLongDef} \label{#1}}

\newcommand{\pbLongDefB}[2]{\begin{pLongDef}
\label{#1} \hypertarget{{#1}}{}
\bookmark[
rellevel=1,
keeplevel,
dest=#1
]{\ref{#1}: #2}
}

\newcommand{\pbLongDefBT}[2]{\begin{pLongDef}[#2]
\label{#1} \hypertarget{{#1}}{}
\bookmark[
rellevel=1,
keeplevel,
dest=#1
]{\ref{#1}: #2}
}

\newcommand{\pbProofB}[2]{ {\bf Proof of {#1} \ref{#2}} \hypertarget{{#2}}{}
\bookmark[
rellevel=1,
keeplevel,
dest=#2
]
{{#1} \ref{#2}}
}

\newcommand{\peProof}[2]{\qed{}}

\newcommand{\pbVerifyB}[1]{ {\bf Verification of Example \ref{#1}}  \hypertarget{{#1}}{}
\bookmark[
rellevel=1,
keeplevel,
dest=#1
]
{Example \ref{#1}}
}

\newcommand{\peVerify}[1]{\qed{}}

}

\newcommand{\wabs}[1]{\left|#1\right|}

\newcommand{\wcal}[1]{\mathcal{#1}}

\newcommand{\wceil}[1]{\lceil #1 \rceil}

\newcommand{\wdf}[2]{{#1}'\left(#2\right)}

\newcommand{\wdsf}[2]{{#1}''\!\!\left(#2\right)}

\newcommand{\wf}[2]{{#1}\left(#2\right)}
\newcommand{\wfc}[2]{{#1}\!\left(#2\right)}

\newcommand{\wflm}{\wrm{fl}}
\newcommand{\wflmt}{\wrm{Fl}}
\newcommand{\wflmk}[1]{\wrm{fl}_{#1}}

\newcommand{\wfl}[1]{\wfc{\wrm{fl}}{#1}}
\newcommand{\wflx}[1]{\wfc{\tilde{\wrm{fl}}}{#1}}
\newcommand{\wflmx}{\tilde{\wrm{fl}}}
\newcommand{\wflmtx}{\tilde{\wrm{Fl}}}
\newcommand{\wflmtxf}[1]{\wfc{\tilde{\wrm{Fl}}}{#1}}
\newcommand{\wflmxk}[1]{\tilde{\wrm{fl}}_{#1}}

\newcommand{\wfltx}[1]{\wfc{\tilde{\wrm{Fl}}}{#1}}
\newcommand{\wflxkf}[2]{\wfc{\tilde{\wrm{fl}}_{#1}}{#2}}

\newcommand{\wflt}[1]{\wfc{\wrm{Fl}}{#1}}

\newcommand{\wflk}[2]{\wfc{\wrm{fl}_{#1}}{#2}}

\newcommand{\wflr}[2]{\wfc{\wrm{#1}}{#2}}

\newcommand{\wflrk}[3]{\wfc{\wrm{#1}_{#2}}{#3}}

\newcommand{\wfloor}[1]{\lfloor {{#1}} \rfloor }

\newcommand{\wi}[1]{\wrm{i}}

\newdimen\CdotAxis
\newcommand*{\CdotAux}[3]{%
  {%
    \settoheight\CdotAxis{$#2\vcenter{}$}%
    \sbox0{%
      \raisebox\CdotAxis{%
        \scalebox{#1}{%
          \raisebox{-\CdotAxis}{%
            $\mathsurround=0pt #2#3$%
          }%
        }%
      }%
    }%
    \dp0=0pt %
    \sbox2{$#2\bullet$}%
    \ifdim\ht2<\ht0 %
      \ht0=\ht2 %
    \fi
    \sbox2{$\mathsurround=0pt #2#3$}%
    \hbox to \wd2{\hss\usebox{0}\hss}%
  }%
}
\newcommand{\wlr}[1]{\left( #1 \right)}

\newcommand{\wn}{\mathds N}

\newcommand{\wvnorm}[2]{{\left\| \wvec{#1}_{#2} \right\|}}

\newcommand{\wrone}{\mathds R}
\newcommand{\wrn}[1]{{\mathds R}^{#1}}
\newcommand{\wrm}[1]{\mathrm{#1}}

\newcommand{\wset}[1]{{\left\{ #1 \right\}}}

\newcommand{\wtr}{\mathrm{T}}
\newcommand{\wsign}[1]{\wfc{\wrm{sign}}{#1}}

\newcommand{\wvec}[1]{\mathbf{#1}}

\newcommand{\wones}{\mathds{1}}

\newcommand{\wz}{\mathds Z}

\newcommand{\wfpa}{\alpha}

\newcommand{\wfpes}[1]{\wcal{E}_{\wcal{A},e}}

\newcommand{\wfpmu}{\mu}
\newcommand{\wfpnu}{\nu}

\newcommand{\wfpsumk}[1]{S_{#1}}
\newcommand{\wfpsumkf}[2]{\wfc{S_{#1}}{#2}}

\newcommand{\wfpbin}{\wcal{E}}

\newcommand{\wfpf}{\wcal{F}}
\newcommand{\wfpsm}{\wcal{M}}
\newcommand{\wfpsi}{\wcal{I}}
\newcommand{\wfpsc}{\wcal{P}}

\newcommand{\wfps}{\wcal{S}}

\newcommand{\wfpemin}{e_{\alpha}}

\begin{document}
\title{Floating point numbers are real numbers}
\author[1]{Walter F. Mascarenhas \thanks{walter.mascarenhas@gmail.com}}
\affil[1]{
IME, Universidade de S\~{a}o Paulo}

\maketitle

\begin{abstract}
Floating point arithmetic allows us to use a finite machine,
the digital computer, to reach conclusions about models
based on continuous mathematics. In this article we
work in the other direction, that is,
we present examples in which continuous mathematics
leads to sharp, simple and new results about the
evaluation of sums, square roots and dot products
in floating point arithmetic.
\end{abstract}

\maketitle

\section{Introduction}
\label{secIntro}
According to Knuth \cite{Knuth}, floating point arithmetic has been
used since Babylonia (1800 B.C.). It
played an important role in the beginning of modern computing,
as in the work of Zuse in the late 1930s. Today we have several models for floating point
arithmetic. Some of these models are based on algebraic
structures, like Kulisch's Ringoids \cite{Kulisch}.
Others models validate numerical software
and lead to automated proofs of results
about floating point arithmetic \cite{Boldo}.

There are also models based on continuous mathematics, which are used
intuitively. For example, when analysing algorithms based
on the floating point operations $\wrm{op} \in \wset{+,-,*,/}$,
executed with machine precision $\epsilon$,
one usually argues that
\pbDef{defEpsArg}
\wfl{x \, \rm{op} \, y} = \wlr{\, x \,\wrm{op} \, y \, } \wlr{1 + \delta}
\hspace{1cm} \wrm{with} \hspace{1cm} \wabs{\delta} \leq  \epsilon,
\peDef{defEpsArg}
where $\wfl{z}$ is the rounded value of $z$.
Equation \pRef{defEpsArg} is called ``the $\wlr{1 + \epsilon}$ argument.''
It may not apply in the presence of underflow,
but  lead to many results in the
hands of Wilkinson \cite{WilkinsonA,WilkinsonB}.
The effectiveness of the $(1 + \epsilon)$ argument is illustrated
by Equation 3.4 in Higham \cite{Higham}, which
expresses the dot product $\hat{d}_n$ of the vectors
$\wvec{x},\wvec{y} \in \wrn{n}$ as
\pbDef{highamsDot}
\hat{d}_n = x_1 y_1 \wlr{1 + \theta_n} +
x_2 y_2 \wlr{1 + {\theta_n}'} + x_3 y_3 \wlr{1 + \theta_{n-1}}
+ \dots + x_n y_n \wlr{1 + \theta_2}.
\peDef{highamsDot}
The $\theta_k$ above are bounded in terms of the
unit roundoff $u$ as
\pbDef{highamsTheta}
\wabs{\theta_k} \leq \frac{k u}{1 - k u} =: \gamma_k,
\peDef{highamsTheta}
and Equations \pRef{highamsDot} and \pRef{highamsTheta} are a good
example of the use of continuous mathematics to analyze
floating point operations. They express well the effects of rounding
errors on dot products, and will suffice for most people
interested in their numerical evaluation.

The purpose of this article is to simplify and extend
the results about floating point arithmetic obtained
using the $(1 + \epsilon)$ argument. We argue that
by thinking of the set of floating point
numbers as a subset of the real numbers we can
use techniques from continuous mathematics to derive and prove
non trivial results about floating point arithmetic.
For instance, we show that in many circumstances we
can replace Higham's $\gamma_k$ by its linearized
counterpart $k u$ and still obtain rigorous bounds.
For us, the replacement of $\gamma_k$ by $k u$
is interesting because it leads to simpler
versions of our articles \cite{Masc,MascCam,arxiv}, and arguments
by other people could be simplified as well. We could,
for example, replace some of Wilkinson's $1.06$ factors by $1$.
In fact, we can even replace $\gamma_k$ by $k u / \wlr{1 + k u}$
when estimating the effects
of rounding errors in the evaluation of the sum
$\wfl{\sum_{k=0}^n x_k}$
of $n + 1$ numbers. Instead of
\pbTClaim{naiveSum}
\wabs{\wfl{\sum_{k=0}^n x_k} -
 \sum_{k=0}^n x_k} \leq \frac{n u}{1 - n u} \sum_{k=0}^n \wabs{x_k},
\peTClaim{naiveSum}
we prove the sharper bound
\pbTClaim{sharpSumB}
\wabs{\wfl{\sum_{k = 0}^n x_k} -  \sum_{k=0}^n x_k}
\leq \frac{n u}{1 + n u} \sum_{k = 0}^n \wabs{x_k},
\peTClaim{sharpSumB}
for arithmetics with subnormal numbers, when we
round to nearest with unit roundoff $u$ and $20 n u \leq 1$.
This bound grows slightly less than linearly
with $n u$, that is, the right hand side is a strictly concave function of
$n u$. Due to this concavity, we can rigorously
conclude from Equation \pRef{sharpSumB} that
\pbTClaim{sharpSumS}
\wabs{ \wfl{\sum_{k = 0}^n x_k} - \sum_{k = 0}^n x_k} \leq n u \sum_{k = 0}^n \wabs{x_k},
\peTClaim{sharpSumS}
and Equation \pRef{sharpSumS} is simpler than Equation \pRef{naiveSum}.
When $x_k \geq 0$, Equation \pRef{sharpSumS} can be improved to
\pbTClaim{sharpSumSB}
x_k \geq 0 \Rightarrow
\wabs{ \wfl{\sum_{k = 0}^n x_k} - \sum_{i = 0}^n x_k} \leq u \sum_{k = 1}^n \sum_{i = 0}^k x_i,
\peTClaim{sharpSumSB}
which is also simple and does not have higher order terms in $u$.
We also analyze dot products, and derive simple and rigorous bounds like
\pbTClaim{sharpDotIntro}
\wabs{ \wfl{\sum_{k = 0}^n x_k y_k} - \sum_{k = 0}^n x_k y_k} \leq
\wlr{n + 1} u \sum_{k = 0}^n \wabs{x_k y_k} \leq \wlr{n + 1} u \sqrt{\sum_{k = 0}^n x_k^2} \sqrt{\sum_{k = 0}^n y_k^2},
\peTClaim{sharpDotIntro}
provided that $\sum_{k = 0}^n \wabs{x_k y_k}$ is not too small,
for the formal concepts of small presented in the statement of our
results.

The bound in Equation \pRef{sharpSumB} is new, and was derived with the theory introduced
in the present article, but Equations \pRef{sharpSumS} and
\pRef{sharpDotIntro} are not new. Actually,
stronger versions of them were already proved by C.-P. Jeannerod and S. Rump in \cite{Rump1} and \cite{Rump2},
and we present them here only as an improvement of the older results in \cite{Higham}.
When comparing the content of these references with
the present work, please note that there is a slight difference in our bound \pRef{sharpDotIntro}
and similar bounds in them: our dot products
involve $n + 1$ pairs of numbers, whereas the sums and
dot products in \cite{Higham,Rump1,Rump3} are defined for $n$ numbers, or pairs of numbers.
Therefore, to leading order in $u$, Equation \pRef{sharpDotIntro} states exactly
the same as the analogous equations (3.7) in \cite{Higham}:
\pbDef{HighamFoo}
\wabs{ \wvec{x}^{\wtr{}} \wvec{y} - \wfl{\wvec{x}^{\wtr{}} \wvec{y}}} \leq
n u \wabs{\wvec{x}}^{\wtr{}} \wabs{\wvec{y}} + \wfc{O}{u^2},
\peDef{HighamFoo}
because $n$ in \cite{Higham} is the same as $n + 1$ for us. Our contribution
is to show that there is no need for the $\wfc{O}{u^2}$ term in Equation \pRef{HighamFoo}
or the 1.06 factor in Wilkinson's Equation 6.11 \cite{WilkinsonB},
implicit in his exponent $t_2$.
We also extend it to situations in which we may have underflow because,
as one may expect after reading \cite{DemmelA,DemmelB,HighamSum,Higham,Neumaier}, the
bounds above must be corrected in order to handle underflow or
arithmetics without subnormal numbers. In the rest of the article we describe such
corrections.

Equations \pRef{sharpSumS} and \pRef{sharpSumSB} are simpler than Equation \pRef{naiveSum},
but the proofs we present for them are definitely not.
However, we hope that after our bounds are validated
via the usual peer review process or by automated tools,
people  will be able to use them without reading their
complicated proofs. For this reason we divided the article
in three parts (besides this introduction.)
In Section \ref{secDef} we define the terms which allow
us to treat floating point numbers as particular cases
of real numbers. In Section \ref{secSharp} we illustrate the use of the definitions
in Section \ref{secDef} to derive sharper and simpler bounds for the effects
on rounding errors in fundamental operations in floating point
arithmetic, like sums, products, square roots and dot products.

Readers should focus on Section \ref{secSharp}. It would
be nice if they could find better proofs for the results stated
in that section. In fact, we are glad that after we posted the first
version of this manuscript M. Lange and S. Rump \cite{Rump3}
derived stronger versions of some results presented here,
using more direct arguments.
This does not contradict the effectiveness
of the use of continuous mathematics to analyze
floating point arithmetic. Our point is that
we can deduce the results thinking in continuous terms,
and their formal proofs is just the last step
in the discovery process.

In the last part of the article we prove our results.
We try to handle all details in our proofs,
and this makes them long and tedious. For this reason, we wrote
two versions of the article. We plan to publish the long one,
and the very long one will be available at arxiv.org.
While reading any of these versions, we ask the reader
not to underestimate how easily ``short and intuitive''
arguments about floating point arithmetic can be wrong.
For example, in the appendix of our article \cite{Extended}
we argue that we can gain intuition about what would happen
if we were to round upward instead of to nearest by replacing
$u$ by $2u$. This argument is correct in that context, because
we verified each and every floating point operation in our
computations. However, this intuitive argument is
not rigorous in general. In fact, by replacing $u$ by $2u$ in the bounds
for rounding to nearest in the present article one will not
obtain rigorous bounds for arithmetics which round upward
or downward.

%%%%%%%%%%%%%%%%%%%%%%%%%%%%%%%
% full version
%%%%%%%%%%%%%%%%%%%%%%%%%%%%%%%
\iftoggle{LatexFull} {

Finally, this extended version of the article is meant to
be read using a software like the Adobe Acrobat Reader,
so that you can click on the hyperlinks (anything in blue)
and follow them. For example, the statement of our
lemmas end with a blue triangle. By clicking
on this triangle you will access the proof of the corresponding
result, and by clicking on the ``back button''
you will return to the statement.  Please, do use this feature
of your reader in order to select which arguments
to follow in more detail. Otherwise, you will find this article
to be unbearably long.

} %% latexFull
%%%%%%%%%%%%%%%%%%%%%%%%%%%%%%%
% full version
%%%%%%%%%%%%%%%%%%%%%%%%%%%%%%%

\section{Definitions}
\label{secDef}
This section presents models of floating point arithmetic
which extend the floating point operations to all real
numbers. In the same way that one can use complex
analysis to study integer arithmetic,
and Sobolev spaces and distributions to learn about regular
solutions of differential equations, by thinking of the
set of floating point numbers as a subset of the set
of real numbers we can use abstract
arguments from optimization theory, point set topology
and convex analysis to
reason about the floating point arithmetics implemented
in real computers.

Most of our floating point numbers have the form
$x = \pm \beta^{e} \wlr{\beta^{\mu} + r}$ where
$\beta \in \wset{2, 3, 4, \dots}$ is the base,
$e$ is an integer exponent, the exponent $\mu$ is a positive integer, and
the remainder $r$ is an integer in $[0,\wlr{\beta - 1} \beta^{\mu} )$.
We also define zero as a floating point
number and, finally, our models account for
subnormal numbers $s = \pm \beta^{e}r $, for
an integer $r \in [1, \beta^{\mu})$.
We now define floating point numbers more formally.

\pbLongDefBT{longDefBase}{Base}
A base is an integer greater than one.
\peLongDef{longDefBase}

\pbLongDefBT{longDefEps}{Unit roundoff}
The unit roundoff associated to the base $\beta$ and
the positive integer exponent $\mu$ is $u := u_{\beta,\mu} := 1/\wlr{2 \beta^\mu}$
(We omit the subscript from $u_{\beta,\mu}$ when
$\beta$ and $\mu$ are evident given the context.)
\peLongDef{longDefEps}

The unit roundoff is our measure of rounding errors.
It is equal to half the distance from $1$ to the next floating
point number. Some authors express their results in terms of ulps (units in the
last place) or the machine precision,
and our $u$ correspond to half of the ulp or the machine epsilon used by them.
However,  the reader must be aware that there are conflicting definitions
of these terms in the literature, and  there is
no universally accepted convention.
A choice must be made, and we prefer to follow Higham \cite{Higham}
and use the unit roundoff $u$ in Definition \ref{longDefEps}.

Our models are based on {\it floating point systems},
which are subsets of $\wrone{}$ to which we round real numbers.
The simplest floating point systems are the perfect ones,
which are defined below.
\pbLongDefBT{longDefMinusSet}{Minus set}
For $\wcal{A} \subset \wrone{}$, we define $-\wcal{A} := \wset{-x, \ \wrm{for} \  x \in \wcal{A}}$.
\peLongDef{longDefMinusSet}
\pbLongDefBT{longDefSign}{Sign function}
The function
$\wrm{sign}: \wrone{} \rightarrow \wrone{}$ is given by
$\wsign{0} := 1$ and $\wsign{x} = \wabs{x}/x$ for $x \neq 0$,
that is, we define the sign of $0$ as one.
\peLongDef{longDefMinusSet}

\pbLongDefBT{longDefBinade}{Equally spaced range (E)}
The equally spaced range associated to the integer exponent $e$,
the base $\beta$ and the positive exponent $\mu$ is
\[
\wfpbin_{e,\beta,\mu} := \wset{\beta^e \wlr{\beta^{\mu} + r} \ \ \wrm{for} \ \ r = 0,1,2,3,\dots, \wlr{\beta -1 } \beta^{\mu} - 1}
\]
(We write simply $\wfpbin_{e}$ when $\beta$ and $\mu$ are evident given the context.)
\peLongDef{longDefBinade}

\pbLongDefBT{longDefPerfect}{Perfect system}
The perfect floating point system associated to the base $\beta$ and the positive exponent $\mu$ is
\[
\wfpsc := \wfpsc_{\beta,\mu} := \wset{0} \,
\bigcup \, \wlr{\bigcup_{e = -\infty}^{\infty}   \wfpbin_{e,\beta,\mu}}
\bigcup \, \wlr{\bigcup_{e = -\infty}^{\infty} - \wfpbin_{e,\beta,\mu}}.
\eqno \eod
\]
\peLongDefX{longDefPerfect}

Perfect floating point systems are convenient for proofs, but ignore
underflow and overflow and are not practical.  It is our opinion that
the best compromise to handle overflow is to
assume that it does not happen, that is, to formulate
models which do not take overflow into account and shift the
burden to handle overflow to the users of the model.
This opinion is not due to laziness, but to the fact
that verifying the absence of overflow in particular cases is
simpler than dealing with floating point systems in which there
is a maximum element.

Underflow is more subtle than overflow, and it may
be difficult to avoid it even in simple cases. Therefore,
handling underflow in each particular case would be too
complicated, and it is a better compromise to have models
that take underflow into account.
Such models are formulated by limiting the range of the exponents
$e$ in Definition \ref{longDefPerfect}.
\pbLongDefBT{longDefMPFR}{MPFR system}
The MPFR system associated to the base $\beta$,
the positive integer $\mu$ and the integer exponent $\wfpemin < -\mu$ is
\[
 \wfpsm :=  \wfpsm_{\wfpemin,\beta,\mu} := \wset{0} \,
 \bigcup \,  \wlr{\bigcup_{e = \wfpemin}^{\infty} \wfpbin_{\beta,\mu}}
 \bigcup \,  \wlr{\bigcup_{e = \wfpemin}^{\infty} - \wfpbin_{\beta,\mu}}.
 \eqno \eod
\]
\peLongDefX{longDefMPFR}
The name MPFR is a tribute to the MPFR library \cite{MPFR},
which has been very helpful in our studies
of floating point arithmetic. This library does not use subnormal numbers,
but allows for very wide exponent ranges (the
minimal exponent is $1- 2^{30}$ in the default configuration.)
As a result, underflow is very unlikely and when it does happen
its consequences are minimal.
\pbLongDefBT{longDefSubNormal}{Subnormal numbers}
\index{subnormal number}
The set of positive subnormal numbers associated to
the base $\beta$, the positive integer exponent $\mu$
 and the integer exponent $\wfpemin$ is
$\wfps_{\wfpemin} := \wfps_{\wfpemin,\beta,\mu} :=
\wset{\beta^{\wfpemin} r, \ \ \wrm{with} \ \ r = 1,2,\dots, \beta^{\mu} - 1}$.
\peLongDef{longDefSubNormal}
\pbLongDefBT{longDefIEEE}{IEEE system}
The IEEE system associated to the base $\beta$, the
positive exponent $\mu$ and  the integer exponent $\wfpemin$, with $\wfpemin < -\wfpmu$, is
\[
\wfpsi := \wfpsi_{\wfpemin, \beta, \mu} := \wset{0} \,
\bigcup \, \wfps_{\wfpemin,\beta,\mu} \,
\bigcup \, -\wfps_{\wfpemin,\beta,\mu} \,
\bigcup \, \wlr{\bigcup_{e = \wfpemin}^{\infty} \wfpbin_{e,\beta,\mu}} \,
\bigcup \, \wlr{\bigcup_{e = \wfpemin}^{\infty} -\wfpbin_{e,\beta,\mu}}.
\]
The elements of $\wfps_{\wfpemin,\beta,\mu} \cup - \wfps_{\wfpemin,\beta,\mu}$
are the subnormal numbers for $\wfpsi$.
\peLongDef{longDefIEEE}

The name IEEE is due to the
IEEE 754 Standard for floating point arithmetic \cite{IEEE},
which contemplates subnormal numbers.

\pbLongDefBT{defFloatSys}{Floating point system}
\index{floating point system}
There are three kinds of floating point systems:
\begin{itemize}
\item The perfect ones in Definition \ref{longDefPerfect}, which
do not contain subnormal numbers.
\item The unperfect ones, which can be either
\begin{itemize}
\item The IEEE systems in Definition \ref{longDefIEEE},
which have subnormal numbers, or
\item The MPFR systems in Definition \ref{longDefMPFR},
which do not have subnormal numbers.
\end{itemize}
\end{itemize}
For brevity, we refer to ``the floating point system $\wfpf$'' as ``the system $\wfpf$,''
and throughout the article the letter $\wfpf$ will always refer to a floating
point system.
\peLongDef{longDefFloatSys}

Please pay attention to the technical detail that,
in order to avoid pathological cases, our
definitions require that $\beta \geq 2$ and $\mu > 0$,
so that the mantissas of our floating point numbers
have at least two bits and $u \leq 1/4$. Additionally, the minimum exponent
$\wfpemin$ for unperfect systems is smaller than $-\mu$, so that
$1$ and $1/\beta$ are floating point numbers. By limiting the exponent range,
we also limit the size of the smallest positive floating
point numbers, which are quantified by the numbers
$\wfpa$ and $\wfpnu$ below.

\pbLongDefBT{longDefAlpha}{Alpha}
For a perfect system we define $\wfpa := 0$;
the IEEE system $\wfpsi_{\wfpemin,\beta,\mu}$ has $\wfpa := \beta^{\wfpemin}$, and
$\wfpa := \beta^{\wfpemin + \wfpmu}$ for the MFPR system $\wfpsm_{\wfpemin,\beta,\mu}$
(Informally, the set of non negative elements of a system begins at $\alpha$.)
\peLongDef{longDefAlpha}

\pbLongDefBT{longDefNu}{Nu}
For a perfect system we define $\wfpnu := 0$ and
the unperfect system $\wfpf_{\wfpemin,\beta,\mu}$ has $\wfpnu := \beta^{\wfpemin + \mu}$.
(Informally, the  {\bf N}ormalized range for a system is formed by the numbers $z$ with
$\wabs{z} \geq \nu$, and $\nu$ is the Greek {\bf N}.)
\peLongDef{longDefNu}

\pbLongDefBT{longDefERange}{Exponent for F}
Any integer $e$ is an exponent for a perfect system,
and $e \in \wz{}$ is an exponent for the unperfect system
$\wfpf_{\wfpemin}$ if $e \geq \wfpemin$.
\peLongDef{longDeERange}

This article is about rounding to nearest, as we now formalize.

\pbLongDefBT{longDefRound}{Rounding to nearest}
A function $\wflm: \wrone{} \rightarrow \wrone{}$ rounds to nearest
in the floating point system $\wfpf$ if $\wfl{z} \in \wfpf$ and
$\wabs{\wfl{z} - z} \leq \wabs{x - z}$ for $x \in \wfpf$ and $z \in \wrone{}$.
\peLongDef{longDefRound}

\pbLongDefBT{longDefTies}{Breaking ties}
When $\wflm$ rounds to nearest in $\wfpf$,
we say that $\wflm$ breaks ties downward if,
for $x \in \wfpf$ and $z \in \wrone{}$,
$\wabs{x - z} = \wabs{\wfl{z} - z} \Rightarrow x \geq \wfl{z}$.
Similarly, $\wflm$ breaks ties upward if
$\wabs{x - z} = \wabs{\wfl{z} - z} \Rightarrow x \leq \wfl{z}$.
\peLongDef{longDefTies}

We now model the numerical sum  $\wfl{\sum_{k = 0}^n y_k}$ of $n+1$
real numbers. For technical reasons, it is important to allow for
the use of different rounding functions in the evaluation of the partial
sums $s_k = \wlr{\sum_{i = 0}^{k-1} y_i} + y_{k}$.
With this motivation, we state the last definitions in this section.

\pbLongDefBT{longDefTuple}{Rounding tuples}
A tuple of functions
$\wflmt = \wset{\wflmk{1},\dots,\wflmk{n}}$
rounds to nearest in $\wfpf$
if  its elements
round to nearest in $\wfpf$. In this case we say that $\wflmt$ is
a rounding $n$-tuple, $n$ is $\wflmt$'s dimension and $\wfpf$ is
$\wflmt$'s range.
\peLongDef{longDefTuple}

\pbLongDefBT{longDefProj}{Projection}
Let $\wcal{A}$ be a set and $\wcal{A}^n$ the Cartesian product
$\wcal{A} \times \dots \times \wcal{A}$ with $n$ factors.
For $k = 1,\dots,n$, we define $\wrm{P}_k: \wcal{A}^n \rightarrow \wcal{A}^k$ as
the projection on the first $k$ coordinates, that is
$\wfc{\wrm{P}_k}{x_1,\dots,x_n} := \wlr{x_1,\dots,x_k}$.
When $\wcal{A}$ is a vector space with zero element $\wvec{0}$, we define
$P_0: \wcal{A}^n \rightarrow \wset{\wvec{0}}$ as
$\wfc{\wrm{P}_0}{x_1,\dots,x_n} := \wvec{0}$.
\peLongDef{longDefProj}

\pbLongDefBT{longDefPSum}{Floating point sum}
Let $\wcal{R}$ be the set of all functions from $\wrone$ to $\wrone{}$,
and $f_0$ its zero element.
We define $S_0: \wset{0} \times \wset{f_0} \rightarrow \wrone{}$
as $\wfpsumkf{0}{0,f_0} := 0$.
For $n > 0$ we define $S_n: \wrn{n} \times \wcal{R}^n \rightarrow \wrone{}$
recursively as
$\wfpsumkf{n}{\wvec{z}, \, \wflmt}
 \,  := \,  \wflk{n}{\wfc{S_{n-1}}{\wrm{P}_{n-1} \wvec{z}, \wrm{P}_{n-1} \wflmt} + z_{n}}.
$
\peLongDef{longDefPSum}

As a convenient notation, given a rounding $n$-tuple $\wflmt$ we write
\[
\wflt{\sum_{k = 0}^n x_k} := \wfc{S_n}{\wlr{x_0 + x_1,x_2,x_3,x_4,\dots,x_n},\wflmt},
\]
and when $\wflmt = \wset{\wflm,\wflm,\dots,\wflm}$ has all its elements equal to $\wflm$ we write
\[
\wfl{\sum_{k = 0}^n x_k} := \wflt{\sum_{k = 0}^n x_k}.
\]
We ask the reader to forgive us for the inconsistency in these
expressions: neither $\wflt{\sum_{k = 0}^n x_k}$ nor
$\wfl{\sum_{k = 0}^n x_k}$ is the value of a
function $\wfc{\wflmt}{s}$ at $s = \sum_{k = 0}^n x_k$, but rather the
value obtained by rounding the partial sums using the elements of $\wflmt$.
Note also that $\wflt{\sum_{k = 0}^n x_k}$ is defined in terms of $x_0 + x_1$, that is, the first
term in the sum is treated differently from the others. The same
detail is present in Equation \pRef{highamsDot}, in which
$x_1 y_1$ and $x_2 y_2$ are treat differently from the other terms.

We emphasize that we define ``the floating point sum of $n + 1$ real numbers'',
and not the ``the sum of $n+1$ floating point numbers.'' As a result,
our rounded sums apply to all real numbers, not only to the ones
in the system $\wfpf$, in the spirit of the first paragraph of this section.
Dot products are similar to sums:
\pbLongDefBT{longDefFpDot}{Dot product}
The dot product of the vectors $\wvec{x},\wvec{y} \in \wrn{n + 1}$
evaluated with the rounding tuples $\wflmt = \wset{\wflm_1,\dots,\wflm_n}$
and $\wrm{R} = \wset{\wrm{r}_0,\dots,\wrm{r}_n}$ is
\[
\wfc{\wrm{dot}_{\wflmt,\wrm{R}}}{\sum_{k = 0}^n x_k y_k} :=
\wflt{\sum_{k=0}^n \wflrk{r}{k}{x_k y_k}}
\eqno \eod
\]
\peLongDefX{longDefFpDot}
We also analyze dot products evaluated with the
fused multiply add operations available in modern hardware and programming languages:
\pbLongDefBT{longDefFmaDot}{Fma dot product}
The fma dot product of the vectors $\wvec{x}, \wvec{y} \in \wrn{n+1}$
evaluated with the rounding tuple $\wflmt = \wset{\wflm_0,\dots,\wflm_n}$ is
\[
\wfc{\wrm{fma}_{\wflmt}}{\sum_{k=0}^n x_k y_k} :=
\wfpsumkf{n+1}{\wlr{x_0 y_0, x_1 y_1, x_2 y_2, \dots, x_n y_n},\wflmt}.
\eqno \eod
\]
\peLongDefX{longDefFmaDot}

\section{Sharp error bounds}
\label{secSharp}
This section presents sharper versions of the $(1 + \epsilon)$ argument. In summary,
we argue that when rounding to nearest
with  unit roundoff $u$,
in many situations we can use
\pbTClaim{bestEps}
\epsilon = \frac{u}{1 + u}
\peTClaim{bestEps}
in the $(1 + \epsilon)$ argument, and this
value is better than $u$ or $u/(1 - u)$.
The section has four parts. The first part
describes the advantages of the $\epsilon$ in
Equation \pRef{bestEps} when dealing with
a few floating point  operations. The next one
generalizes our results to sums of many numbers, by proving
the bound \pRef{sharpSumB}.
Section \ref{secCumSum} presents
bounds on the errors in sums which are expressed in
terms of $\sum_{k = 1}^n \wabs{\sum_{i = 0}^k x_i}$.
Section \ref{secDot} is about dot products.
It shows that by working with real numbers from the start
it is easy to adapt results derived for sums in order
to obtain bounds for the errors in dot products.

\subsection{Basics}
\label{secBasics}
This section is about the $(1 + \epsilon)$ argument
for a few floating point operations.
When rounding a floating point number,
our first lemma states that the $\epsilon$ in Equation \pRef{bestEps} can
be used when the real number $z$ is in the normal range, ie.,
the  absolute value of $z$ is greater than the number
$\wfpnu$ in Definition \ref{longDefNu}.
\pbLemmaBT{lemUNear}{A better epsilon}
If $\wrm{fl}$ rounds to nearest in $\wfpf$ and $\wabs{z} \geq \wfpnu_{\wfpf}$ then
\pbTClaim{rhoNear}
\wabs{\wfl{z} - z} \leq \frac{\wabs{z} u}{1 + u}.
\peTClaim{rhoNear}
In particular, if $\wfpf$ is perfect then
Equation \pRef{rhoNear} holds for all $z \in \wrone{}$.
\peLemma{lemUNear}
Lemma \ref{lemUNear} is sharp in the sense that for any $\epsilon$ smaller than
 $u/\wlr{1+u}$ there exists a real number $z$ near $1 + u$ for
which Equation \pRef{rhoNear} does not hold. It has been known for
a long time \cite{Knuth}, but it leads to bounds
slightly stronger than the ones in \cite{Higham} for instance, because
\[
\frac{u}{1 + u} < u < \frac{u}{1 - u}
\]
and when the result  of the operation
$x \, \wrm{op} \, y \neq 0$ is in the normal range we have the bound
\pbDef{nearRound}
\frac{1}{1 + u} \leq
\frac{\wfl{x \, \wrm{op} \, y}}{x \, \wrm{op} \, y} \leq \frac{1 + 2 u}{1 + u},
\peDef{nearRound}
instead of the usual bound
\pbDef{usualRound}
1 - u \leq \frac{ \wfl{x \, \wrm{op} \, y}}{x \, \wrm{op} \, y} \leq \frac{1}{1 - u}.
\peDef{usualRound}
As a result, we could use the same argument as Higham to
conclude that in a perfect floating point system
\pbDef{sharpSum}
\wflt{\sum_{k=0}^n x_k} = \wlr{x_0 + x_1} \xi_0^n + \sum_{i = 2}^n x_k \xi_k^{n - i + 1}
\hspace{0.5cm} \wrm{with} \hspace{0.5cm}
\frac{1}{1 + u} \leq \xi_k \leq \frac{1 + 2u}{1 + u}
\peDef{sharpSum}
and
\[
\wflt{\sum_{k=0}^n x_k y_k} = x_0 y_0 \xi_0^{n + 1} + \sum_{i = 1}^n x_k y_k \xi_k^{n - i + 2}
\hspace{0.5cm} \wrm{with} \hspace{0.5cm}
\frac{1}{1 + u} \leq \xi_k \leq \frac{1 + 2u}{1 + u}.
\]

The underlying reason as to why
\pbDef{concF}
\wfc{f}{u} := \frac{1 + 2 u}{1 + u}
\hspace{0.5cm} \wrm{is \ a  \ better \ upper\  bound \ than}  \hspace{0.5cm}
\wfc{h}{u} := \frac{1}{1 - u}
\peDef{concF}
is the difference between concavity and convexity.
The function $\wfc{f_\tau}{x} := \wfc{f}{u}^{\tau}$ has second derivative
\[
\wdsf{f_\tau}{u} = \frac{\tau \wfc{f_{\tau-2}}{u}}{\wlr{1 + u}^4} \wlr{\tau - 3  - 4 u}
\]
and is concave for $\tau \leq 3 + 4 u$. On the other hand, $\wfc{h_{\tau}}{x} := \wfc{h}{u}^{\tau}$ has second
derivative
\[
\wdsf{h_{\tau}}{u} = \tau \wlr{\tau + 1} \wfc{h_{\tau-2}}{u}
\]
and is convex for all $\tau > 0$ and $0 < u < 1$. As a result, we can linearize
rigorously an upper bound based on $f_{\tau}$, with $0 \leq \tau \leq 3 + 4 u$,
whereas linearizing an upper bound based on $h_{\tau}$ is
correct only to leading order. For instance,
using the bound \pRef{nearRound} we can prove the next corollary,
and similar results combining multiplications and divisions,
but we could not prove such results based only on the usual bound \pRef{usualRound}.
\pbCorolBT{corFourProd}{Three products} Let $x,y,z$ and $w$ be real numbers. If
$\hat{p}_1 := \wfl{x * y}$,
$\hat{p}_2 := \wfl{\hat{p}_1 * z}$ and
$\hat{p}_3 := \wfl{\hat{p}_2 * w}$,
$p_i \neq 0$, and $\wabs{p_i}$ satisfy Equation \pRef{rhoNear} for
$k = 1$, $2$ and $3$ then
\[
1 - k u \leq \frac{\hat{p}_k}{p_k} \leq 1 + k u,
\]
for $p_1 := x * y$, $p_2 := x * y * z$ and $p_3 := x * y * z * w$.
\peFullCorol{corFourProd}
Lemma \ref{lemUNear} also yields a simple proof of a well known
result about square roots when $\beta = 2$  \cite{Cody,Kahan}, and solves
an open problem for arbitrary bases $\beta$ \cite{BoldoB}:
\pbCorolBT{corSqrt}{Square roots}
For the base $\beta = 2$, if $x \in \wfpf$ is such that $x^2 \geq \wfpnu$ and
$\wfl{x^2}$ and  $\wfl{\sqrt{\wfl{x^2}}}$ are evaluated
rounding to nearest then $\wfl{\sqrt{\wfl{x^2}}} = \wabs{x}$.
Moreover,
\pbDef{thSqrt}
\wfl{ \frac{\wabs{x}}{ \wfl{\sqrt{\wfl{x^2}}} }} \leq 1
\peDef{thSqrt}
for a general base $\beta$, under the same hypothesis on $\wrm{fl}$ and $x$.
\peCorol{corSqrt}
The next two lemmas show that there are other conditions
besides $\wabs{z} \geq \wfpnu$ in which we can use
the bound in Equation \pRef{rhoNear}:

\pbLemmaBT{lemSmallSum}{Exact sums}
If $x, y \in \wfpf$ are such that $\wfpa \leq \wabs{x + y} \leq \beta \wfpnu$ then
$z := x + y \in \wfpf$, that is, the sum $x + y$ is exact.
In particular, $z$ satisfies Equation \pRef{rhoNear}.
\peLemma{lemSmallSum}

\pbLemmaBT{lemIEEESum}{IEEE sums}
Let $\wfpsi$ be an IEEE system
and $x,y \in \wfpsi$. If $0 < \wabs{x + y} \leq \beta \nu$
then $\wabs{x + y} \geq \wfpa$,
and $z := x + y \in \wfpsi$ and satisfies Equation \pRef{rhoNear}.
\peLemma{lemIEEESum}
The last two lemmas combined with Lemma \ref{lemUNear} imply that
we can use the bound \pRef{rhoNear} for every real number $z$ which
is the sum of two floating point numbers in an IEEE system.
This is yet another instance in which subnormal numbers lead to
simpler results, and corroborates Demmel's arguments \cite{DemmelA,DemmelB}
and the soundness of the decision to include subnormal numbers
in the IEEE standard for floating point arithmetic \cite{IEEE}.
Another instance is the fundamental Sterbenz's Lemma, which
must be modified by the inclusion of the term $\wfpa$ in
its hypothesis in order to hold for MPFR systems:

\pbLemmaBT{lemSterbenz}{Sterbenz's Lemma}
If $a,b \in \wfpf$ and $\wfpa \leq b - a \leq a$ then
$b - a \in \wfpf$.
\peFullLemma{lemSterbenz}

\subsection{Norm one bounds}
\label{secOneNorm}

This subsection extends Lemma \ref{lemUNear} to sums with many parcels.
Our results are described by the next lemma and its corollaries.
In particular, we show that underflow does not affect
sums of positive numbers. Therefore,
there is no need for terms involving the smallest positive
floating point number when bounding the errors in such sums.

\pbLemmaBT{lemNormOneBound}{Norm one bound}
If $\wflmt = \wset{\wflmk{1},\dots,\wflmk{n}}$
rounds to nearest in a perfect system,
$20 n u \leq 1$ and $y_0,\dots y_n \in \wrone{}$ then
\pbTClaim{normOneBound}
\wabs{
\, \wflt{\sum_{k = 0}^n y_k} - \sum_{k = 0}^n y_k
\,
} \leq
\frac{n u}{1 + n u} \sum_{k = 0}^n \wabs{y_k}.
\peTClaim{normOneBound}
\peLemma{lemNormOneBound}

\pbCorolBT{corNormOneIEEE}{IEEE norm one bound}
If $\wflmt = \wset{\wflmk{1},\dots,\wflmk{n}}$
 rounds to nearest in an IEEE system $\wfpsi$,
 $20 n u \leq 1$ and $y_0,\dots, y_n \in \wfpsi$
then Equation \pRef{normOneBound} is satisfied.
\peFullCorol{corNormOneIEEE}

\pbCorolBT{corNormOneMPFR}{MPFR norm one bound}
If $\wflmt$
rounds to nearest in a MPFR system $\wfpsm$,
$20 n u \leq 1$, $\wvec{y} \in \wfpsm^{n+1}$
and $y_k \geq 0$ for all $k$
then Equation \pRef{normOneBound} holds.
\peFullCorol{corNormOneMPFR}

\pbCorolBT{corNormOneUnperfect}{Unperfect norm one bound}
If $\wflmt = \wset{\wflmk{1},\dots,\wflmk{n}}$
rounds to nearest in an unperfect system,
$y_0,\dots,y_n \in \wrone{}$ and $20 n u \leq 1$ then
\pbTClaim{thSharpSumUnperfect}
\wabs{
\, \wflt{\sum_{k = 0}^n y_k} - \sum_{k = 0}^n y_k
\,
} \leq
\frac{n \alpha}{2} + \frac{n u}{1 + n u} \wlr{\frac{n \wfpa}{2} + \sum_{k = 0}^n \wabs{y_k}}.
\peTClaim{thSharpSumUnperfect}
If, additionally, $u \sum_{k = 0}^n \wabs{y_k} \geq \wfpa$
then Equation \pRef{normOneBound} is satisfied.
\peCorol{corNormOneUnperfect}

Note that Corollaries \ref{corNormOneIEEE} and \ref{corNormOneMPFR} have
different hypothesis regarding the floating point numbers $y_0,\dots,y_n$:
in the IEEE case, in which we have subnormal numbers, Equation \pRef{normOneBound}
holds for all such $y_k$. In the MPFR case, due to the absence of subnormal numbers,
we must assume that $y_k \geq 0$, for
Equation \pRef{normOneBound} does not hold for instance when $\beta = 2$, $x_0 = 3 \wfpa/2$,
$x_1 = -\wfpa$, $n = 1$ and we break ties upward.
Note also that the number $\wfpa$ in Equation \pRef{thSharpSumUnperfect}
for an IEEE system is much smaller than the
$\wfpa$ for the corresponding MPFR system.
The next example shows that the bound \pRef{normOneBound} is sharp:
\pbExampleBT{exSharpSum}{The norm one bound is sharp}
If $\wrm{fl}$ rounds to nearest in the perfect system $\wfpsc_{2,\mu}$,
breaking ties downward, and $x_0 := 1$ and $x_k := u$ for  $k = 1, \dots,n$ then
\[
\wfl{\sum_{k = 0}^n x_k} = 1 = \sum_{k=0}^n x_k - n u =
\sum_{k=0}^n x_k  - \frac{n u}{1 + n u} \sum_{k = 0}^n x_k.
\]
If $\wrm{fl}$ breaks ties upward for the same $x_k$ and $2 n u < 1$ then
\[
\wfl{\sum_{k = 0}^n x_k}  = 1 + 2 n u = \sum_{k=0}^n x_k + n u = \sum_{k=0}^n x_k
  + \frac{n u}{1 + n u} \sum_{k = 0}^n x_k.
\eqno \eod
\]
\pePlainExampleX{exSharpSum}

As in Lemma \ref{lemUNear}, the bound
\pRef{normOneBound} has concavity properties which
allow us to linearize rigorously bounds resulting from
a couple of its applications:
\pbLemmaBT{lemConvexity}{Convexity} For $k \in \wn{}$ and $i = 1,\dots,k$, let $n_i$ be a
positive number and define functions $f_k, g_k: (0,\infty) \rightarrow \wrone{}$ by
\[
\wfc{f_k}{u} =
\prod_{i = 1}^k \frac{1 + 2 n_i u}{1 + n_i u}
\hspace{1cm} \wrm{and} \hspace{1cm}
\wfc{g_k}{u} =
\prod_{i = 1}^k \frac{1}{1 + n_i u}.
\]
The functions $f_k$ are  strictly concave for $k = 1$, $2$ and $3$
and the functions $g_k$ are convex for all $k$.
In particular, for $k \leq 3$,
\[
1 - \wlr{\sum_{i = 1}^k n_i}  u \leq \wfc{g_k}{u} \leq \wfc{f_k}{u} \leq 1 + \wlr{\sum_{i = 1}^k n_i} u.
\pFullLink{lemConvexity}
\]
\peLemmaX{lemConvexity}

As a final point for this section, we note that Lemma \ref{lemNormOneBound} implies that
\pbTClaim{sharpSumNearMax}
\wabs{
\, \wflt{\sum_{k =0}^n y_k} - \sum_{k = 0}^n y_k
\,
} \leq
\frac{n \wlr{n + 1} u}{1 + n u} \max_{k = 0,\dots,n} \wabs{y_k},
\peTClaim{sharpSumNearMax}
and it is natural to ask whether the quadratic term in $n$
in the right hand side of Equation \pRef{sharpSumNearMax} is necessary. The next
example shows that bounds in terms of $\max \wabs{y_k}$ do need
a quadratic term in $n$ (or large constant factors):

\pbExampleBT{exQuadGrowth}{Quadratic growth}
If $\wflm$ rounds to nearest in the perfect system $\wfpsc_{2,\mu}$,
breaking ties downward,
$y_0 := 1 + u$ and $y_k := 1 + 2^{\wfloor{\wfc{\log_2}{k + 1}}} u$
for $k = 1,\dots, n := 2^m - 1$, where $m \in \wn{}$ is such that $2^m u < 1$,
then
\[
\wfl{\sum_{k = 0}^n y_k} = \sum_{k = 0}^n y_k - \frac{n^2 + 2 n + 3}{3} u \leq
\sum_{k = 0}^n y_k - \frac{n^2 + 2 n + 3}{6} u \max_{k = 0,\dots, n} \wabs{y_k}.
\pFullLink{exQuadGrowth}
\]
\peExampleX{exQuadGrowth}

\subsection{Cumulative bounds}
\label{secCumSum}

Although Lemma \ref{lemNormOneBound} leads to the simple bound
\pRef{sharpSumS} on the error in the evaluation of sums,
it is not as good from the qualitative
view as the result one would obtain from the version
of Higham's Equation  3.4 for sums, or from our Equation \pRef{highamsDot}.
We believe that Higham and Wilkinson would write this equation as
\pbDef{wilkSum}
\wfl{\sum_{k = 1}^n x_k} = \wlr{x_1 + x_2} \wlr{1 + \theta_{n-1}} +  x_3 \wlr{1 + \theta_{n-2}}
+ \dots + x_n \wlr{1 + \theta_1}.
\peDef{wilkSum}
for $\theta_n$ in Equation \pRef{highamsTheta}.
In fact, Wilkinson presents an expression similar to Equation \pRef{wilkSum}
for sums using a double precision accumulator in page 117 of \cite{WilkinsonB}.
Equation \pRef{wilkSum} gives a better intuition regarding the effects of rounding
errors in the corresponding sum than the bound
in Lemma \ref{lemNormOneBound}. Therefore, it is natural to look for bounds that
take into account the stronger relative influence of the first
parcels in Equation \pRef{wilkSum}.
The next examples are relevant in this context.

\pbExampleBT{exSharpSumNearB}{Minimum cumulative bound}
If $\wrm{fl}$ rounds to nearest in the perfect system $\wfpsc_{2,\mu}$,
breaking ties downward,
and $x_k := u^{-k}$ for $k = 1,\dots,n$ then
\[
\wfl{\sum_{k = 0}^n x_k} = u^{-n} = \sum_{k = 0}^n x_k  - \kappa_n u \sum_{k = 1}^n \sum_{i = 0}^k x_i
\]
for
\[
1 - u < \kappa_n
 := \frac{\wlr{1 - u}\wlr{{1 - u^{n}}}}{1 - u^{n} - n u^{n+1} \wlr{1 - u}}
%:= \frac{\wlr{1 - u}\wlr{{1 - u^{n}}}}{1 - u^{n+2} - n u^{n+1} \wlr{1 - u}}
< \wlr{1 - u} \wlr{1 + u^n} < 1.
\pFullLink{exSharpSumNearB}
\]
\peExampleX{exSharpSumNearB}

\pbExampleBT{exSharpSumNearC}{Maximum cumulative bound}
If $\wrm{fl}$ rounds to nearest in the perfect system $\wfpsc_{\beta,\mu}$,
breaking ties upward,
$1 = e_1 < e_2 \dots < e_n$ are integer exponents,
$x_0 := u$, $x_1 := 1$,
and $x_k := \beta^{e_k} \wlr{1 + u} - \beta^{e_{k-1}} \wlr{1 + 2 u}$
for $k = 2,\dots, n$ then
\pbDef{ssBnC}
\wfl{\sum_{k=0}^n x_k} \leq
\sum_{i = 0}^n x_i + \tau_n u  \sum_{k=1}^n \sum_{i = 0}^k x_k,
\peDef{ssBnC}
for
\[
\tau_n := \frac{1}{1 + u \wlr{\frac{\beta - 2}{\beta - 1} + \frac{n}{\beta^{n} - 1}}}.
\]
Additionally, if $e_k = k - 1$ for $k \geq 1$ then
we have equality in Equation \pRef{ssBnC}.
\peFullExample{exSharpSumNearC}

These examples indicate that there is an asymmetry between the upper and lower
bounds on the errors $\delta := \wfl{\sum_{k=0}^n x_k} - \sum_{i = 0}^n x_i$
in terms of $\sum_{k = 1}^n \sum_{i = 0}^k x_k$:
The constants $\kappa_n$ in Example \ref{exSharpSumNearB} and $\tau_n$ in Example \ref{exSharpSumNearC} are
equal to $1/\wlr{1 + u}$ for $n = 1$ but as $n$ increases,
$\kappa_n$ decreases toward $1 - u$ whereas $\tau_n$ increases toward $1$.
Therefore, the worst lower and upper values for
$\delta$ are reached in different
situations, and are due to distinct causes.
In fact, the lower bound for $\delta$ in the next Lemma
is a straightforward consequence of the convexity of the functions
$\wlr{1 + u}^{-k}$ and Equation
\pRef{sharpSum}, whereas the upper bound is a non
trivial consequence of the concavity of $\wlr{1 + 2u}/\wlr{1 + u}$.
\pbLemmaBT{lemPositiveBound}{Positive cumulative bound}
If $\wflmt = \wset{\wflmk{1},\dots,\wflmk{n}}$
rounds to nearest in a perfect system,
$y_0,y_1,\dots y_n \in \wrone{}$, with $y_k\geq 0$ for $k = 0,\dots,n$,
and $20 n u \leq 1$ then
\pbTClaim{thPositiveBound}
- \frac{u}{1 + u} \sum_{k = 1}^n \sum_{i = 0}^k y_k
 \leq
\wflt{\sum_{k=0}^n y_k}
 - \sum_{k = 0}^n y_k
\,
 \leq \tau_n u \sum_{k = 1}^n \sum_{i = 0}^k y_k,
\peTClaim{thPositiveBound}
for $\tau_n$ in Example \ref{exSharpSumNearC}.
\peLemma{lemPositiveBound}
\pbCorolBT{corPositiveUnperfect}{Unperfect cumulative bound}
If $\wflmt$ rounds to nearest in an unperfect system $\wfpf$,
$20 n u \leq 1$, $\wvec{y} \in \wfpf^{n+1}$ and
$y_k \geq 0$, then Equation \pRef{thPositiveBound} holds.
\peCorol{corPositiveUnperfect}
The next example shows that Lemma
\ref{lemPositiveBound} does not apply to sums of numbers with mixed signs,
and Lemma \ref{lemSignedSum} and its corollary show that the example
is nearly worst possible.
\pbExampleBT{exSignedSumB}{Mixed signs}
If $\wrm{fl}$ rounds to nearest in a
perfect system $\wfpsc_{2,\mu}$,
breaking ties upward,
$x_0 := u$, $x_1 := 1$,
$x_k := - 2^{1 - k} \wlr{1 + 3 u}$ for $k > 1$
and $2^n u \leq 1$ then
\pbDef{thXSSA}
\wfl{\sum_{k = 0}^n x_k} - \sum_{k = 0}^n x_k = 2 \wlr{1 - 2^{- n}} u =
 \frac{\kappa_n u}{1 - \wlr{n - 2} u}  \sum_{k = 1}^{n} \wabs{\sum_{i = 0}^k x_i},
\peDef{thXSSA}
for
\[
1 - u \leq \kappa_n :=
\frac{\wlr{1 - 2^{-n}}\wlr{1 - \wlr{n - 2} u}}{\wlr{1 -  2^{-n}} \wlr{1 + 3u} - n u} \leq 1.
\pFullLink{exSignedSumB}
\]
\peExampleX{exSignedSumB}
\pbLemmaBT{lemSignedSum}{Signed cumulative bound}
If $\wflmt = \wset{\wflmk{1},\dots,\wflmk{n}}$
rounds to nearest in a perfect system,
$y_0,y_1,\dots y_n \in \wrone{}$ and $20 n u < 1$ then
\pbTClaim{thSignedSum}
\wabs{ \,
\wflt{\sum_{k=0}^n y_k}
 - \sum_{k = 0}^n y_k
\,
} \leq
 \frac{u}{1 - \wlr{n-2}u} \sum_{k = 1}^n \wabs{\sum_{i = 0}^k y_i}.
\peTClaim{thSignedSum}
\peLemma{lemSignedSum}
\pbCorolBT{corSignedSum}{Unperfect signed cumulative bound}
If $\wflmt = \wset{\wflmk{1},\dots,\wflmk{n}}$
rounds to nearest in an unperfect system,
$y_0,y_1,\dots y_n \in \wrone{}$ and $20 n u \leq 1$ then
\pbTClaim{thCorSignedSum}
\wabs{
\, \wflt{\sum_{k = 0}^n y_k} - \sum_{k = 0}^n y_k
\,
} \leq
\wlr{1 + 2 n u} \frac{n \alpha}{2} +
\frac{u}{1 - \wlr{n - 2} u} \sum_{k = 1}^n \wabs{\sum_{i = 0}^k y_i},
\peTClaim{thCorSignedSum}
If, additionally, $u \sum_{k = 1}^n \wabs{\sum_{i = 0}^k y_i} \geq n \alpha$
then
\pbTClaim{thCorSignedSumL}
\wabs{
\, \wflt{\sum_{k = 0}^n y_k} - \sum_{k = 0}^n y_k
\,
} \leq
\frac{3}{2}\wlr{1 + \frac{n u}{2}} u \sum_{k = 1}^n \wabs{\sum_{i = 0}^k y_i}.
\peTClaim{thCorSignedSumL}
\peCorol{corSignedSum}

\subsection{Dot products}
\label{secDot}
This section presents bounds on the errors in the numerical evaluation of dot products.
These bounds are derived from the ones for sums
presented in Section \ref{secOneNorm}. This derivation is possible
because some of our previous bounds apply to general real numbers,
and a numerical dot product is simply a numerical sum of real numbers,
which may or may not be floating point numbers. If our analysis of
sums were restricted to floating point numbers then the extensions
presented here would be harder to derive. For example, the
next corollaries follow directly from Lemma \ref{lemNormOneBound}
and the Definition \ref{longDefFmaDot} of numerical dot products
using fused multiply adds (these corollaries are proved in the extended
version of the article):

\pbCorolBT{corFmaPerfect}{Dot prod. with fma}
If $\wflmt = \wset{\wflmk{0},\dots,\wflmk{n}}$
rounds to nearest in a perfect system,
$20  n u \leq 1$ and $\wvec{x},\wvec{y} \in \wrn{n+1}$ then
\pbTClaim{thSharpFmaDotNear}
\wabs{
\, \wf{\wrm{fma}_{\wflmt}}{\sum_{k = 0}^n x_k y_k} - \sum_{k = 0}^n x_k y_k
\,
} \leq
\frac{\wlr{n + 1} u}{1 + \wlr{n + 1} u}  \sum_{k = 0}^n \wabs{x_k y_k}.
\peTClaim{thSharpFmaDotNear}
\peFullCorol{corFmaPerfect}

\pbCorolBT{corFmaUnperfect}{Unperfect Dot prod. with fma}
If $\wflmt = \wset{\wflmk{0},\dots,\wflmk{n}}$ rounds to nearest in an unperfect system,
$20 n u \leq 1$ and $\wvec{x},\wvec{y} \in \wrn{n+1}$ then
\[
\wabs{
\, \wf{\wrm{fma}_{\wflmt}}{\sum_{k = 0}^n x_k y_k} - \sum_{k = 0}^n x_k y_k
\,
} \leq
\wlr{n + 1} \frac{\alpha}{2}
\]
\pbTClaim{thFmaUnperfect}
+   \frac{\wlr{n + 1} u}{1 + \wlr{n + 1} u} \wlr{\frac{\wlr{n + 1} \wfpa}{2} + \sum_{k = 0}^n \wabs{x_k y_k}}.
\peTClaim{thFmaUnperfect}
If, additionally, $ u \sum_{k = 0}^n \wabs{x_k y_k} \geq \wfpa$
then Equation \pRef{thSharpFmaDotNear} holds.
\peFullCorol{corFmaUnperfect}

When we evaluate dot products rounding each product $x_k y_k$, the
bounds are slightly worse, but can still be obtained
with the theory in Section \ref{secOneNorm}:

\pbCorolBT{corDotPerfect}{Dot prod.}
If $\wflmt = \wset{\wflmk{1},\dots,\wflmk{n}}$
and $\wrm{R} = \wset{\wrm{r}_0,\dots,\wrm{r}_n}$
round to nearest in a perfect system,
$20 n u \leq 1$ and $\wvec{x},\wvec{y} \in \wrn{n+1}$ then
\pbTClaim{thSharpDotNear}
\wabs{
\, \wfc{\wrm{dot}_{\wflmt,\wrm{R}}}{\sum_{k = 0}^n x_k y_k} - \sum_{k = 0}^n x_k y_k
\,
} \leq \beta_n u \sum_{k = 0}^n \wabs{x_k y_k} \leq \wlr{n + 1} u \sum_{k = 0}^n \wabs{x_k y_k},
\peTClaim{thSharpDotNear}
where
\[
\beta_n := \frac{n + 1 + 3 n u}{1 + \wlr{n + 1} u + n u^2} \leq \frac{n + 1}{1 + n u/2}
\hspace{1cm} \wrm{and} \hspace{1cm}
\beta_n \leq \frac{n + 1}{1 + \wlr{n - 3} u}.
\pFullLink{corDotPerfect}
\]
\peCorolX{corDotPerfect}

\pbCorolBT{corDotIEEE}{IEEE dot prod.}
If $\wflmt = \wset{\wflmk{1},\dots,\wflmk{n}}$
and $\wrm{R} = \wset{\wrm{r}_0,\dots,\wrm{r}_n}$
round to nearest in an IEEE system,
$20 n u \leq 1$ and $\wvec{x},\wvec{z} \in \wrn{n+1}$ then
\pbTClaim{thSharpDotNearIEEEB}
\wabs{
\, \wfc{\wrm{dot}_{\wflmt,\wrm{R}}}{\sum_{k = 0}^n x_k y_k} - \sum_{k = 0}^n x_k y_k
\,
} \leq
1.05 \wlr{n + 1} \frac{\wfpa}{2} + \beta_n u \sum_{k = 0}^n \wabs{x_k y_k},
\peTClaim{thSharpDotNearIEEEB}
for $\beta_n$ in Corollary  \ref{corDotPerfect}.
If, additionally, $u \sum_{k = 0}^n \wabs{x_k y_k} \geq \wfpa$
then
\[
\wabs{\, \wfc{\wrm{dot}_{\wflmt,\wrm{R}}}{\sum_{k = 0}^n x_k y_k} - \sum_{k = 0}^n x_k y_k
\,
}
\leq \frac{3}{2} \wlr{n + 1} u \sum_{k = 0}^n \wabs{x_k y_k}.
\pFullLink{corDotIEEE}
\]
\peCorolX{corDotIEEE}

\pbCorolBT{corDotMPFR}{MPFR dot prod.}
If $\wflmt = \wset{\wflmk{1},\dots,\wflmk{n}}$
and $\wrm{R} = \wset{\wrm{r}_0,\dots,\wrm{r}_n}$
round to nearest in a MPFR system,
$20 n u \leq 1$ and
$\wvec{x},\wvec{z} \in \wrn{n+1}$ then
\pbTClaim{thSharpDotNearMPFRB}
\wabs{
\, \wfc{\wrm{dot}_{\wflmt,\wrm{R}}}{\sum_{k = 0}^n x_k y_k} - \sum_{k = 0}^n x_k y_k
\,
} \leq
\frac{\wlr{2.05 n + 1.05}\wfpa}{2} + \beta_n \sum_{k = 0}^n \wabs{x_k y_k},
\peTClaim{thSharpDotNearMPFRB}
for $\beta_n$ in Corollary \ref{corDotPerfect}.
If,  additionally, $u \sum_{k = 0}^n \wabs{x_k y_k} \geq \wfpa$
then the last equation in Corollary \ref{corDotIEEE} is satisfied.
\peFullCorol{corDotMPFR}

Finally, in all bounds above we can use the Cauchy-Schwarz inequality
and replace the terms  $\sum_{k = 0}^n \wabs{x_k y_k}$ by
$\wvnorm{x}{}_2 \wvnorm{y}{}_2$.
With this replacement, we can compare our bounds to the ones in
\cite{DemmelA} and \cite{Neumaier}.

\section{Proofs}
\label{secProofs}
This section we prove our main results.
Section \ref{secAuxDefs} contains more definitions
and Section \ref{secAuxLemmas} presents more lemmas.
In Section \ref{secProp} we state basic results
about floating point systems and rounding to nearest. We call such results
by ``Propositions,'' because they are obvious and readers should
be able to deduce them with little effort.
Section \ref{secLemmas} begins with the proofs of the main lemmas,
and after that we prove some of the corollaries.
The extended version of the article contains the
proofs of the remaining lemmas and corollaries and the propositions,
and the verification of the examples.

\subsection{More definitions}
\label{secAuxDefs}
The proofs of our bounds on the errors in sums use the following definitions:

\pbLongDefBT{longDefTightFunction}{Tight function}
Let $\wcal{A}$ and $\wcal{B}$ be topological spaces and $\wcal{R}$ a set.
A function $f: \wcal{A} \times \wcal{R} \rightarrow \wcal{B}$ is
tight if for every sequence
$\wset{\wlr{a_k,r_k}, k \in \wn{}} \subset \wcal{A} \times \wcal{R}$  such that
$\lim_{k \rightarrow \infty} a_k$
there exists $r \in \wcal{R}$
and a subsequence $\wset{\wlr{a_{n_k},r_{n_k}}, k \in \wn{}}$
with
$\lim_{k \rightarrow \infty} \wfc{f}{a_{n_k},r_{n_k}} = \wfc{f}{a,r}$.
\peLongDef{longDefTightFunction}

\pbLongDefBT{longDefTightSet}{Tight set of functions}
Let $\wcal{A}$ and $\wcal{B}$ be topological spaces and let
$\wcal{R}$ be a set of functions from $\wcal{A}$ to $\wcal{B}$.
We say that $\wcal{R}$ is tight if the function
$f: \wcal{A} \times \wcal{R} \rightarrow \wcal{B}$ given
by $\wfc{f}{a,r} := \wfc{r}{a}$ is tight.
\peLongDef{longDefTightSet}

\subsection{More lemmas}
\label{secAuxLemmas}

\pbLemmaBT{lemUNearS}{Sharp epsilons}
Suppose $\wflm$ rounds to nearest in $\wfpf$ and $e$ is an exponent for $\wfpf$.
If $\wabs{z} = \beta^{e} \wlr{\beta^{\mu} + w}$ with
$w \in [0, \wlr{\beta - 1}\beta^{\mu}]$ then
\pbTClaim{rHat}
\wfl{z} = \wsign{z} \beta^{e}\wlr{\beta^{\mu} + r}
\hspace{0.3cm}
\wrm{for} \hspace{0.3cm} \ r \in [0, \wlr{\beta - 1} \beta^{\mu}] \cap \wz{}
\hspace{0.3cm} \wrm{and} \hspace{0.3cm}
\wabs{r - w} \leq 1/2.
\peTClaim{rHat}
Moreover,
\pbTClaim{uNearSA}
\wabs{\frac{\wfl{z} - z}{z}} \leq \frac{u}{1 + \max \wset{1, 2 w} u} \leq \frac{u}{1 + u}
\peTClaim{uNearSA}
and
\[
\wabs{\frac{\wfl{z} - z}{z}} \leq \frac{u}{1 + \wlr{2 r - 1} u}
\hspace{0.5cm} \wrm{and} \hspace{1cm}
\wabs{\frac{\wfl{z} - z}{\wfl{z}}} \leq \frac{u}{1 + 2 r u}.
\pLink{lemUNearS}
\]
\peLemmaX{lemUNearS}

%
%
%%%%%%%%%%%%%%%%%%%%%%%%%%%%%%%%%%%%%%%%%%%%%%%%%%%%%%%%%%%%%%%%%%%%%%%%%%%%%%%%%%%
% LemOpt
%%%%%%%%%%%%%%%%%%%%%%%%%%%%%%%%%%%%%%%%%%%%%%%%%%%%%%%%%%%%%%%%%%%%%%%%%%%%%%%%%%%

\pbLemmaBT{lemOpt}{Compactness}
Let $\wcal{R}$ be a set, $\wcal{L} \subset \wrone{} \setminus \wset{0}$,
$\wcal{A},\wcal{B} \subset \wrn{n}$, $\wcal{X} \subset \wrn{m}$  and
$\wcal{K} \subset \wcal{A}$.
Define $\wcal{Z} := \wcal{A} \cup \wcal{B}$.
If the functions
$f: \wcal{Z} \times \wcal{X} \rightarrow \wrone{}$,
$h: \wcal{Z} \times \wcal{R} \rightarrow \wcal{X}$,
and $\wfc{g}{\wvec{z},r} := \wfc{f}{\wvec{z}, \wfc{h}{\wvec{z},r}}$
and $\varphi \in \wrone{}$ are such that
\begin{itemize}
\item $\wcal{K}$ is compact and for $\wvec{z} \in \wcal{A}$ there exist $\lambda \in \wcal{L}$
such that $\lambda \wvec{z} \in \wcal{K}$.
\item If $\lambda \in \wcal{L}$, $\wvec{z} \in \wcal{A}$, $\lambda \wvec{z} \in \wcal{K}$ and
$r \in \wcal{R}$ then $\wfc{h}{\lambda \wvec{z},r'} = \lambda \wfc{h}{\wvec{z},r}$
for some $r' \in \wcal{R}$.
\item $f$ is upper semi-continuous and $\wfc{f}{\lambda \wvec{z}, \lambda \wvec{x}} \geq \wfc{f}{\wvec{z},\wvec{x}}$
for $\wvec{z} \in \wcal{A}$ and $\lambda \in \wcal{L}$.
\item $h$ is tight, in the sense of Definition \ref{longDefTightFunction}.
\item  $\wfc{g}{\wvec{z},r} \leq \varphi$ for $\wlr{\wvec{z},r} \in \wcal{B} \times \wcal{R}$.
\end{itemize}
then either $\wfc{g}{\wvec{z},r} \leq \varphi$ for all
$\wlr{\wvec{z},r} \in \wcal{Z} \times \wcal{R}$
or there exist $\wlr{\wvec{z}^*,r^*} \in  \wcal{K} \times \wcal{R}$ such that
$\wfc{g}{\wvec{z}^*,r^*} \geq \wfc{g}{\wvec{z},r}$ for all $\wlr{\wvec{z},r} \in \wcal{Z} \times \wcal{R}$.
\peLemma{lemOpt}

Lemma \ref{lemOpt} is a compactness argument. Its purpose
is to show that either there exists examples for which the
relative effects of rounding errors
are the worst possible or these errors are small.
It is necessary because  floating point systems are
infinite and we cannot take this existence for granted.
The intuition behind Lemma \ref{lemOpt} is simple.
The vector $\wvec{z}$  represents the input to
computation. The vector $\wfc{h}{\wvec{z},r}$ is obtained
by rounding functions of $\wvec{z}$ using the rounding functions
$r \in \wcal{R}$. The bad set $\wcal{B}$ represents situations like underflow or very poor scaling,
and its elements are handled separately.
For $\wvec{z}$ outside of the bad set, we can use scaling
by powers of $\beta$ (represented by $\lambda \in \wcal{L}$)
to reduce the analysis of $\wfc{f}{\wvec{z},\wvec{x}}$ to
real numbers $\wvec{z}$  in the compact set $\wcal{K}$. We can then deal
with the discontinuity in rounding by analyzing
all functions which round to nearest
(represented by $\wcal{R}$) instead of a single function. In the end,
as in the applications of the classic Banach-Alaoglu Theorem,
by using compactness and continuity in their full generality,
we can analyze the existence of maximizers for the relative effects of
rounding errors. We can then exploit the implications of maximality in order
to describe precisely such maximizers.

\subsection{Propositions}
\label{secProp}
In this section we present auxiliary results about floating point systems.
We believe readers will find most of them to be trivial,
and they are presented only to make our arguments more precise.
In all propositions $\beta$ is a base, $\mu$ is a positive integer,
$\wfpf$ is a floating point system associated to $\beta$ and $\mu$,
$z \in \wrone{}$, $x \in \wfpf$,
and $u$, $\wfpemin{}$, $\wfpa$ and $\wfpnu$ are the numbers
related to this system in Definitions \ref{longDefEps},
\ref{longDefPerfect}, \ref{longDefMPFR}, \ref{longDefIEEE},
\ref{longDefAlpha} and \ref{longDefNu}, Finally
the function $\wflm$ rounds to nearest in $\wfpf$.

%
%
%%%%%%%%%%%%%%%%%%%%%%%%%%%%%%%%%%%%%%%%%%%%%%%%%%%%%%%%%%%%%%%%%%%%%%%%%%%%%%%%%%%
% propOrder
%%%%%%%%%%%%%%%%%%%%%%%%%%%%%%%%%%%%%%%%%%%%%%%%%%%%%%%%%%%%%%%%%%%%%%%%%%%%%%%%%%%

\pbPropBT{propOrder}{Order by the exponent}
Let $d$ and $e$ be integers and $v,w \in \wrone{}$,
with $v <  \wlr{\beta - 1}\beta^{\wfpmu}$ and $w \geq 0$.
If $d < e$ then $\beta^{d} \wlr{\beta^{\wfpmu} + v} < \beta^{e} \wlr{\beta^{\wfpmu} + w}$.
\peFullProp{propOrder}

%
%
%%%%%%%%%%%%%%%%%%%%%%%%%%%%%%%%%%%%%%%%%%%%%%%%%%%%%%%%%%%%%%%%%%%%%%%%%%%%%%%%%%%
% propNormalForm
%%%%%%%%%%%%%%%%%%%%%%%%%%%%%%%%%%%%%%%%%%%%%%%%%%%%%%%%%%%%%%%%%%%%%%%%%%%%%%%%%%%

\pbPropBT{propNormalForm}{Normal form}
If $z \in \wrone{}$ is different from zero then there exist
unique $e \in \wz{}$ and $w \in [0,  \wlr{\beta - 1}\beta^{\wfpmu})$
such that $z = \wsign{z} \beta^{e} \wlr{\beta^{\wfpmu} + w}$.
\peFullProp{propNormalForm}

%
%
%%%%%%%%%%%%%%%%%%%%%%%%%%%%%%%%%%%%%%%%%%%%%%%%%%%%%%%%%%%%%%%%%%%%%%%%%%%%%%%%%%%
% propIntegerForm
%%%%%%%%%%%%%%%%%%%%%%%%%%%%%%%%%%%%%%%%%%%%%%%%%%%%%%%%%%%%%%%%%%%%%%%%%%%%%%%%%%%

\pbPropBT{propIntegerForm}{Integer form}
If $\wfpf$ is unperfect and $x \in \wfpf$ then
there exists $e,r \in \wz{}$ with $e \geq \wfpemin$
such that $x = \beta^e r$ and
\begin{itemize}
\item $r = 0$ if and only if $x = 0$.
\item $0 < \wabs{r} < \beta^{\mu}$ if and only if $x$ is subnormal
and $\wfpf$ is an IEEE system.
\item $\beta^{\mu} \leq \wabs{r} < \beta^{1 + \mu}$ if and only if $\wabs{x} \in \wfpbin_e$.
\end{itemize}
\peFullProp{propIntegerForm}

%
%
%%%%%%%%%%%%%%%%%%%%%%%%%%%%%%%%%%%%%%%%%%%%%%%%%%%%%%%%%%%%%%%%%%%%%%%%%%%%%%%%%%%
% propSymmetry
%%%%%%%%%%%%%%%%%%%%%%%%%%%%%%%%%%%%%%%%%%%%%%%%%%%%%%%%%%%%%%%%%%%%%%%%%%%%%%%%%%%

\pbPropBT{propSymmetry}{Symmetry}
$\wfpf$ is symmetric, that is, $x \in \wfpf \Leftrightarrow -x \in \wfpf \Leftrightarrow \wabs{x} \in \wfpf$.
\peFullProp{propSymmetry}

%
%
%%%%%%111%%%%%%%%%%%%%%%%%%%%%%%%%%%%%%%%%%%%%%%%%%%%%%%%%%%%%%%%%%%%%%%%%%%%%%%%%%%%%
% propNu
%%%%%%%%%%%%%%%%%%%%%%%%%%%%%%%%%%%%%%%%%%%%%%%%%%%%%%%%%%%%%%%%%%%%%%%%%%%%%%%%%%%

\pbPropBT{propNu}{The minimality of nu}
Let $x$ be a floating point number.
If $\wabs{x} \geq \nu$ and $x \neq 0$ then $x$ is normal, that is, there exists
an exponent $e$ for $\wfpf$ such that $\wabs{x} \in \wfpbin_{e}$.
If $0 < \wabs{x} < \nu$ then $\wfpf$ is an IEEE system $\wfpsi_{\wfpemin}$ and $x$ is subnormal,
that is, $\wabs{x} \in \wfps_{\wfpemin}$.
Conversely, if $e,r \in \wz{}$ and $z = \beta^{e} r$ with
$\wabs{r} \leq \beta^{1 + \mu}$ and $\wabs{z} \geq \wfpnu$ then $z \in \wfpf$.
\peFullProp{propNu}

%
%
%%%%%%%%%%%%%%%%%%%%%%%%%%%%%%%%%%%%%%%%%%%%%%%%%%%%%%%%%%%%%%%%%%%%%%%%%%%%%%%%%%%
% propSubnormalSum
%%%%%%%%%%%%%%%%%%%%%%%%%%%%%%%%%%%%%%%%%%%%%%%%%%%%%%%%%%%%%%%%%%%%%%%%%%%%%%%%%%%

\pbPropBT{propSubnormalSum}{Subnormal sum}
Let $\wfpsi$ be an IEEE system.
If $x,y \in \wfpsi$ are subnormal then $x + y \in \wfpsi$.
\peFullProp{propSubnormalSum}

%
%
%%%%%%%%%%%%%%%%%%%%%%%%%%%%%%%%%%%%%%%%%%%%%%%%%%%%%%%%%%%%%%%%%%%%%%%%%%%%%%%%%%%
% propCriticalSum
%%%%%%%%%%%%%%%%%%%%%%%%%%%%%%%%%%%%%%%%%%%%%%%%%%%%%%%%%%%%%%%%%%%%%%%%%%%%%%%%%%%

\pbPropBT{propCriticalSum}{Critical sum}
If $e$ is an exponent for $\wfpf$ and
$x \in \wfpf$ and $z \in \wrone{}$ are such that
$\wabs{x + z} = \beta^{e} \wlr{\beta^{\mu} + r + 1/2}$ with
$r \in [0,\wlr{\beta - 1} \beta^{\mu}) \cap \wz{}$
then $\wabs{z} \geq \beta^{e}/2$.
\peFullProp{propCriticalSum}

%
%
%%%%%%%%%%%%%%%%%%%%%%%%%%%%%%%%%%%%%%%%%%%%%%%%%%%%%%%%%%%%%%%%%%%%%%%%%%%%%%%%%%%
% propRoundIdentity
%%%%%%%%%%%%%%%%%%%%%%%%%%%%%%%%%%%%%%%%%%%%%%%%%%%%%%%%%%%%%%%%%%%%%%%%%%%%%%%%%%%

\pbPropBT{propRoundIdentity}{Identity}
If $x \in \wfpf$ then $\wfl{x} = x$.
\peFullProp{propRoundIdentity}

%
%
%%%%%%%%%%%%%%%%%%%%%%%%%%%%%%%%%%%%%%%%%%%%%%%%%%%%%%%%%%%%%%%%%%%%%%%%%%%%%%%%%%%
% propRoundMonotone
%%%%%%%%%%%%%%%%%%%%%%%%%%%%%%%%%%%%%%%%%%%%%%%%%%%%%%%%%%%%%%%%%%%%%%%%%%%%%%%%%%%

\pbPropBT{propRoundMonotone}{Monotonicity}
If $z \leq w$ then $\wfl{z} \leq \wfl{w}$, and if $x \in \wfpf$ then
\begin{itemize}
\item $x > \wfl{z} \Rightarrow x > z$,
\item $x < \wfl{z} \Rightarrow x < z$,
\item $\wabs{x} > \wabs{\wfl{z}} \Rightarrow \wabs{x} > \wabs{z}$,
\item $\wabs{x} < \wabs{\wfl{z}} \Rightarrow \wabs{x} < \wabs{z}$.
\end{itemize}
\peFullProp{propRoundMonotone}

%
%
%%%%%%%%%%%%%%%%%%%%%%%%%%%%%%%%%%%%%%%%%%%%%%%%%%%%%%%%%%%%%%%%%%%%%%%%%%%%%%%%%%%
% propRoundMinus
%%%%%%%%%%%%%%%%%%%%%%%%%%%%%%%%%%%%%%%%%%%%%%%%%%%%%%%%%%%%%%%%%%%%%%%%%%%%%%%%%%%

\pbPropBT{propRoundMinus}{Symmetric rounding}
$\wflr{m}{z} := - \wfl{-z}$ rounds to nearest in $\wfpf$.
\peFullProp{propRoundMinus}

%
%
%%%%%%%%%%%%%%%%%%%%%%%%%%%%%%%%%%%%%%%%%%%%%%%%%%%%%%%%%%%%%%%%%%%%%%%%%%%%%%%%%%%
% propRoundNormal
%%%%%%%%%%%%%%%%%%%%%%%%%%%%%%%%%%%%%%%%%%%%%%%%%%%%%%%%%%%%%%%%%%%%%%%%%%%%%%%%%%%

\pbPropBT{propRoundNormal}{Normal rounding}
Let $e$ be an exponent for $\wfpf$. If $\wabs{z} = \beta^{e} \wlr{\beta^{\wfpmu} + w}$,
with $w \in [0, \wlr{\beta - 1} \beta^{\wfpmu})$, then $\wfl{z} \in \wset{a,b}$ for
\[
a := \wsign{z} \beta^{e} \wlr{\beta^{\wfpmu} + \wfloor{w}} \in \wfpf
\hspace{0.7cm}
\wrm{and}
\hspace{0.7cm}
b := \wsign{z} \beta^{e} \wlr{\beta^{\wfpmu} + \wceil{w}}  \in \wfpf,
\]
and
\pbDef{thRoundNormal}
\wabs{\wfl{z} - z} = \min \wset{\wabs{z - a}, \, \wabs{z - b}} \leq \frac{\wabs{b - a}}{2} \leq \beta^e/2.
\peDef{thRoundNormal}
If $z < m := \wlr{a + b}/2$ then
$\wfl{z} = \min \wset{a,b}$, and if $z > m$ then $\wfl{z} = \max \wset{a,b}$.

In particular, if $r \in \wz{}$ and $\wabs{r - w} < 1/2$ then
$\wfl{z} = \wsign{z} \beta^e \wlr{\beta^{\mu} + r}$.
\peFullProp{propRoundNormal}

%
%
%%%%%%%%%%%%%%%%%%%%%%%%%%%%%%%%%%%%%%%%%%%%%%%%%%%%%%%%%%%%%%%%%%%%%%%%%%%%%%%%%%%
% propRoundSubNormal
%%%%%%%%%%%%%%%%%%%%%%%%%%%%%%%%%%%%%%%%%%%%%%%%%%%%%%%%%%%%%%%%%%%%%%%%%%%%%%%%%%%

\pbPropBT{propRoundSubnormal}{Subnormal rounding}
Let $\wfpsi = \wfpsi_{\wfpemin}$ be an IEEE system.
If $\wabs{z} \leq \nu$  then $\wfl{z} \in \wset{a,b}$ for
\[
a := \beta^{\wfpemin} \wfloor{\beta^{-\wfpemin} z} \in \wfpsi
\hspace{0.7cm}
\wrm{and}
\hspace{0.7cm}
b :=  \beta^{\wfpemin} \wceil{ \beta^{-\wfpemin} z} \in \wfpsi
\]
and
\[
\wabs{\wfl{z} - z} = \min \wset{z - a, b - z} \leq \frac{b - a}{2} \leq \wfpa/2.
\]
If $z < m:= \wlr{a + b}/2$ then
$\wfl{z} = a$, and if $z > m$ then $\wfl{z} = b$.

If $r \in [-\beta^{\mu}, \beta^{\mu}) \cap \wz{}$ then $\wfl{z} = \beta^{\wfpemin}{r}$
for $\beta^{\wfpemin} \wlr{r - 1/2} < z <  \beta^{\wfpemin} \wlr{r + 1/2}$ and
$\wfl{\beta^{\wfpemin}\wlr{ r + 1/2}} \in \wset{ \beta^{\wfpemin} r, \beta^{\wfpemin}\wlr{r+1}}$.
\peFullProp{propRoundSubnormal}

%
%
%%%%%%%%%%%%%%%%%%%%%%%%%%%%%%%%%%%%%%%%%%%%%%%%%%%%%%%%%%%%%%%%%%%%%%%%%%%%%%%%%%%
% propRoundBelowAlpha
%%%%%%%%%%%%%%%%%%%%%%%%%%%%%%%%%%%%%%%%%%%%%%%%%%%%%%%%%%%%%%%%%%%%%%%%%%%%%%%%%%%

\pbPropBT{propRoundBelowAlpha}{Rounding below alpha}
If $\wabs{z} < \wfpa/2$ then $\wfl{z} = 0$.
If $\wabs{z} = \wfpa/2$ then $\wfl{z} \in \wset{0, \wsign{z} \wfpa}$
and if $\wfpa/2 < \wabs{z} \leq \wfpa$ then $\wfl{z} = \wsign{z} \wfpa$.
In particular, if $\wabs{z} \leq \wfpa$ then $\wabs{\wfl{z} - z} \leq \wfpa / 2$.
\peFullProp{propRoundBelowAlpha}

%
%
%%%%%%%%%%%%%%%%%%%%%%%%%%%%%%%%%%%%%%%%%%%%%%%%%%%%%%%%%%%%%%%%%%%%%%%%%%%%%%%%%%%
% propRoundAdapt
%%%%%%%%%%%%%%%%%%%%%%%%%%%%%%%%%%%%%%%%%%%%%%%%%%%%%%%%%%%%%%%%%%%%%%%%%%%%%%%%%%%

\pbPropBT{propRoundAdapt}{Perfect adapter}
Let $\wfpsc_{\beta,\mu}$ be a perfect system and
$\wfpf_{\wfpemin,\beta,\mu}$ an unperfect one.
If $\wflm$ rounds to nearest in $\wfpf$ then
there exists $\wflmx$ which rounds to nearest in $\wfpsc$
and is such that $\wflx{z} = \wfl{z}$ for
$z$ with $\wabs{z} \geq \nu_{\wfpf}$.
\peFullProp{propRoundAdapt}

%
%
%%%%%%%%%%%%%%%%%%%%%%%%%%%%%%%%%%%%%%%%%%%%%%%%%%%%%%%%%%%%%%%%%%%%%%%%%%%%%%%%%%%
% propRoundIEEEX
%%%%%%%%%%%%%%%%%%%%%%%%%%%%%%%%%%%%%%%%%%%%%%%%%%%%%%%%%%%%%%%%%%%%%%%%%%%%%%%%%%%

\pbPropBT{propRoundIEEEX}{IEEE adapter}
If $\wflmt = \wset{\wflmk{1},\dots,\wflmk{n}}$ rounds to nearest in
the IEEE system $\wfpsi_{\wfpemin,\beta,\mu}$ and $\wfpsc_{\beta,\mu}$ is a perfect system
then there exists $\wflmtx = \wset{\wflmxk{1},\dots,\wflmxk{n}}$
which rounds to nearest in $\wfpsc$ and is such that
$\wflt{\sum_{k=0}^n x_k} = \wflmtxf{\sum_{k=0}^n x_k}$ for
all $\wvec{x} \in \wfpsi^{n+1}$. In particular, $\wfltx{\sum_{k=0}^n x_k} \in \wfpsi$.
\peFullProp{propRoundIEEEX}

%
%
%%%%%%%%%%%%%%%%%%%%%%%%%%%%%%%%%%%%%%%%%%%%%%%%%%%%%%%%%%%%%%%%%%%%%%%%%%%%%%%%%%%
% propRoundMPFRX
%%%%%%%%%%%%%%%%%%%%%%%%%%%%%%%%%%%%%%%%%%%%%%%%%%%%%%%%%%%%%%%%%%%%%%%%%%%%%%%%%%%

\pbPropBT{propRoundMPFRX}{MPFR adapter}
If $\wflmt = \wset{\wflmk{1},\dots,\wflmk{n}}$ rounds to nearest in
the MPFR system $\wfpsm_{\wfpemin,\beta,\mu}$ and $\wfpsc_{\beta,\mu}$ is a perfect system
then there exists $\wflmtx = \wset{\wflmxk{1},\dots,\wflmxk{n}}$
which rounds to nearest in $\wfpsc$ and is such that
$\wflt{\sum_{k=0}^n x_k} = \wflmtxf{\sum_{k=0}^n x_k}$ for
$\wvec{x} \in \wfpsm^{n+1}$ with $\wflt{\sum_{i = 0}^{k-1} x_i} + x_k \geq 0$ for $k = 0,\dots,n$.
\peFullProp{propRoundMPFRX}

%
%
%%%%%%%%%%%%%%%%%%%%%%%%%%%%%%%%%%%%%%%%%%%%%%%%%%%%%%%%%%%%%%%%%%%%%%%%%%%%%%%%%%%
% propFlat
%%%%%%%%%%%%%%%%%%%%%%%%%%%%%%%%%%%%%%%%%%%%%%%%%%%%%%%%%%%%%%%%%%%%%%%%%%%%%%%%%%%

\pbPropBT{propFlat}{Flatness}
Let $e$ be an exponent for $\wfpf$ and $z$ with $\wabs{z} = \beta^{e} \wlr{\beta^{\mu} + w}$
for $w \in [0,\wlr{\beta - 1} \beta^{\mu})$. On the one hand, if
$w = \wfloor{w} + 1/2$ then
\[
\wabs{w - y} < \beta^{e}/2 \ \Rightarrow \ \wfl{y}
\in \wset{ \wsign{z} \beta^{e} \wlr{\beta^{\mu} + \wfloor{w}}, \,
           \wsign{z} \beta^{e} \wlr{\beta^{\mu} + \wceil{w}}}.
\]
On the other hand, if $w - \wfloor{w} \neq 1/2$
then there exists $\delta > 0$ such that if
$\wflmk{1}$ and $\wflmk{2}$ round to nearest in $\wfpf$
and $\wabs{y - z} < \delta$ then $\wflk{1}{y} = \wflk{2}{z}$.
\peFullProp{propFlat}

%
%
%%%%%%%%%%%%%%%%%%%%%%%%%%%%%%%%%%%%%%%%%%%%%%%%%%%%%%%%%%%%%%%%%%%%%%%%%%%%%%%%%%%
% propSumScale
%%%%%%%%%%%%%%%%%%%%%%%%%%%%%%%%%%%%%%%%%%%%%%%%%%%%%%%%%%%%%%%%%%%%%%%%%%%%%%%%%%%

\pbPropBT{propSumScale}{Scaled sums}
Suppose $\wflmt = \wset{\wflmk{1},\dots,\wflmk{n}}$ rounds
to nearest in a perfect system $\wfpsc$ and $\wfpsumk{k}$ is the sum in Definition \ref{longDefPSum}.
If $\wvec{z} \in \wrn{n}$, $\sigma \in \wset{-1,1}$ and $m \in \wz{}$ then there exist
$\wflmtx = \wset{\wflmxk{1},\dots,\wflmxk{n}}$ which round
to nearest in $\wfpsc$ such that
$\wfpsumkf{k}{\sigma \beta^m \, \wvec{z}, \, \wflmt} \, =
              \, \sigma \beta^m \, \wfpsumkf{k}{\wvec{z}, \, \wflmtx}$
for $k = 1,\dots,n$.
\peFullProp{propSumScale}

%
%
%%%%%%%%%%%%%%%%%%%%%%%%%%%%%%%%%%%%%%%%%%%%%%%%%%%%%%%%%%%%%%%%%%%%%%%%%%%%%%%%%%%
% Proof Whole is Tight
%%%%%%%%%%%%%%%%%%%%%%%%%%%%%%%%%%%%%%%%%%%%%%%%%%%%%%%%%%%%%%%%%%%%%%%%%%%%%%%%%%%

\pbPropBT{propWholeIsTight}{Whole is tight}
The set of all functions which round to nearest in $\wfpf$ is tight.
\peFullProp{propWholeIsTight}

%
%
%%%%%%%%%%%%%%%%%%%%%%%%%%%%%%%%%%%%%%%%%%%%%%%%%%%%%%%%%%%%%%%%%%%%%%%%%%%%%%%%%%%
% Proof Sums are tight
%%%%%%%%%%%%%%%%%%%%%%%%%%%%%%%%%%%%%%%%%%%%%%%%%%%%%%%%%%%%%%%%%%%%%%%%%%%%%%%%%%%

\pbPropBT{propSumsAreTight}{Sums are tight}
Let $\wcal{R}$ be a tight set of functions
which round to nearest in $\wfpf$ and $\wfpsumk{k}$ the sum in Definition \ref{longDefPSum}.
The function $T_n: \wrone{}^{n} \times \wcal{R}^n \rightarrow \wrn{n+1}$ given by
\[
\wfc{T_n}{\wvec{z},\wflmt} :=
\wlr{\wfpsumkf{0}{\wvec{z},\wflmt}, \,
     \wfpsumkf{1}{\wvec{z},\wflmt}, \,
     \wfpsumkf{2}{\wvec{z},\wflmt}, \,
     \dots,
     \wfpsumkf{n}{\wvec{z},\wflmt}}
\]
is tight.
\peFullProp{propSumsAreTight}

\subsection{Lemmas}
\label{secLemmas}
This section presents the proofs of the Lemmas other than
\ref{lemSterbenz} and \ref{lemConvexity}, which are proved
in the extended version of the article.\\

%
%
%%%%%%%%%%%%%%%%%%%%%%%%%%%%%%%%%%%%%%%%%%%%%%%%%%%%%%%%%%%%%%%%%%%%%%%%%%%%%%%%%%%
% Proof of lemma lemUNear
%%%%%%%%%%%%%%%%%%%%%%%%%%%%%%%%%%%%%%%%%%%%%%%%%%%%%%%%%%%%%%%%%%%%%%%%%%%%%%%%%%%

\pbProofB{Lemma}{lemUNear}
If $z= 0$ then $\wfl{z} = z = 0$ by Prop. \ref{propRoundIdentity} and
Equation \pRef{rhoNear} holds.
If $z \neq 0$ then, by Prop. \ref{propNormalForm},
$z = \wsign{z} \beta^{e} \wlr{\beta^{\wfpmu{}} + w}$,
with $e \in \wz{}$ and $r \in \wrone{}$ with $w \in [0, \wlr{\beta - 1} \beta^{\wfpmu{}})$.
When $\wfpf$ is unperfect
$\wfpnu = \beta^{\wfpemin + \mu}$  and, by Prop. \ref{propOrder}, $e \geq \wfpemin$ because $\wabs{z} \geq \wfpnu$.
Therefore, $e$ is an exponent for $\wfpf$
and Lemma \ref{lemUNear} follows from Lemma \ref{lemUNearS}.
Finally, Lemma \ref{lemUNear} applies to all $z$ when $\wfpf$
is perfect because $\nu = 0$ in this case.
\peProof{Lemma}{lemUNear}\\

%%%%%%%%%%%%%%%%%%%%%%%%%%%%%%%%%%%%%%%%%%%%%%%%%%%%%%%%%%%%%%%%%%%%%%%%%%%%%%%%%%%
% End of proof of lemma lemUNear
%%%%%%%%%%%%%%%%%%%%%%%%%%%%%%%%%%%%%%%%%%%%%%%%%%%%%%%%%%%%%%%%%%%%%%%%%%%%%%%%%%%
%
%

%
%
%%%%%%%%%%%%%%%%%%%%%%%%%%%%%%%%%%%%%%%%%%%%%%%%%%%%%%%%%%%%%%%%%%%%%%%%%%%%%%%%%%%
% Proof of lemma lemSmallSum
%%%%%%%%%%%%%%%%%%%%%%%%%%%%%%%%%%%%%%%%%%%%%%%%%%%%%%%%%%%%%%%%%%%%%%%%%%%%%%%%%%%

\pbProofB{Lemma}{lemSmallSum} If $\wfpf$
is perfect then $\wfpa = \wfpnu = 0$ and Lemma \ref{lemSmallSum} holds
because $0 \in \wfpf$. It is clear that the Lemma also holds when $x = 0$ or $y =0$,
and from now on we suppose that $x,y \neq 0$ and $\wfpf$ is unperfect.
In this case $\wfpnu = \beta^{\wfpemin + \wfpmu}$, and
we can assume that $\wabs{y} \geq \wabs{x}$ because $x + y = y + x$.
Moreover, $x + y \in \wfpf \Leftrightarrow -(x + y) \in \wfpf$
by Prop. \ref{propSymmetry} and we
can also assume that
\pbDef{SNSA}
\wfpa \leq x + y \leq \beta \nu = \beta^{1 + \wfpemin + \wfpmu}
\hspace{1cm} \wrm{and} \hspace{1cm}
y > 0.
\peDef{SNSA}
If $y$ is subnormal then $0 < \wabs{x} \leq y < \nu$,
$x$ is also subnormal by Prop. \ref{propNu} and
Lemma \ref{lemSmallSum} follows from Prop. \ref{propSubnormalSum}.
Therefore, we can assume that $y$ is normal, that is,
\pbDef{SNSY}
y = \beta^{\wfpemin + e} \wlr{\beta^{\wfpmu} +  r_y}
\hspace{0.5cm} \wrm{with} \hspace{0.5cm}
e \geq 0
\hspace{0.5cm} \wrm{and} \hspace{0.5cm}
r_y \in [0,\wlr{\beta - 1} \beta^{\wfpmu}) \cap \wz{}.
\peDef{SNSY}

On the one hand, if $x > 0$ then Equation \pRef{SNSA} leads to $y < \beta^{1 + \wfpemin{} + \wfpmu}$
and Equation \pRef{SNSY} yields $e = 0$.
Prop. \ref{propIntegerForm} and the assumption $0 < \wabs{x} \leq y$ lead to
\[
x = \beta^{\wfpemin} r_x
\hspace{1cm} \wrm{with} \hspace{1cm}
r_x \in \wz{}
\hspace{0.5cm} \wrm{and} \hspace{0.5cm}
1 \leq r_x \leq \beta^{\wfpmu} + r_y,
\]
and Equations \pRef{SNSA} and \pRef{SNSY} imply that
\[
x + y = \beta^{\wfpemin} \wlr{\beta^{\wfpmu} + r_x + r_y}
\leq \beta^{1 + \wfpemin + \wfpmu} \ \
\Rightarrow x + y \geq \nu \ \ \wrm{and} \ \ r_x + r_y \leq \wlr{\beta - 1} \beta^{\wfpmu},
\]
and Prop. \ref{propNu} with $r = \beta^{\mu} + r_x + r_y$ shows that $x + y \in \wfpf$.

On the other hand, if $x < 0$ then Prop. \ref{propIntegerForm} and the assumption
$0 < \wabs{x} \leq y$ lead to
\[
x = -\beta^{\wfpemin + d} r_x
\hspace{0.5cm} \wrm{with} \hspace{0.5cm}
0 \leq d \leq e,
\hspace{0.5cm}
r_x \in \wz{}
\hspace{0.5cm} \wrm{and} \hspace{0.5cm}
1 \leq r_x < \beta^{1 + \wfpmu}.
\]
It follows that
$
x + y = \beta^{\wfpemin + d} \wlr{\beta^{e - d} \wlr{\beta^{\wfpmu} + r_y} - r_x}
= \beta^{\wfpemin} r
$
for
\[
r := \beta^d \wlr{\beta^{e - d} \wlr{\beta^{\wfpmu} + r_y} - r_x} \in \wz{}.
\]
Since $x < 0$,  using \pRef{SNSA} and the identity $\wfpnu = \beta^{\wfpemin + \mu}$ we deduce that
\[
0 < r = \beta^{-\wfpemin} \wlr{x + y} < \beta^{-\wfpemin} \beta \wfpnu = \beta^{1 + \mu}.
\]
When $\wfpf$ is a MPFR system we have that $\wfpa = \wfpnu$, Equation
\pRef{SNSA} implies that $x + y \geq \wfpnu$ and the
equation above and Prop. \ref{propNu} show that $x + y \in \wfpf$.
Finally, when $\wfpf$ is an IEEE system we either have (i)
$r \geq \beta^{\mu}$, in which case $x + y \geq \beta^{\wfpemin + \mu} = \wfpnu$
and $x + y \in \wfpf$ by Prop. \ref{propNu}, or
(ii) $r < \beta^{\mu}$, and $x + y \in \wfps_{\wfpemin}$
is a subnormal number, which belongs to $\wfpf$.
Therefore, $x + y \in \wfpf$ in all cases and we are done.
\peProof{Lemma}{lemSmallSum}\\

%%%%%%%%%%%%%%%%%%%%%%%%%%%%%%%%%%%%%%%%%%%%%%%%%%%%%%%%%%%%%%%%%%%%%%%%%%%%%%%%%%%
% End of the Proof of lemma lemSmallSum
%%%%%%%%%%%%%%%%%%%%%%%%%%%%%%%%%%%%%%%%%%%%%%%%%%%%%%%%%%%%%%%%%%%%%%%%%%%%%%%%%%%
%
%

%
%
%%%%%%%%%%%%%%%%%%%%%%%%%%%%%%%%%%%%%%%%%%%%%%%%%%%%%%%%%%%%%%%%%%%%%%%%%%%%%%%%%%%
% Proof of lemma lemIEEESum
%%%%%%%%%%%%%%%%%%%%%%%%%%%%%%%%%%%%%%%%%%%%%%%%%%%%%%%%%%%%%%%%%%%%%%%%%%%%%%%%%%%

\pbProofB{Lemma}{lemIEEESum}
By Prop. \ref{propIntegerForm},
$x = \beta^{d} r$ and $y = \beta^{e} s$ for $d,e,r,s \in \wz{}$ such that
$d, e \geq \wfpemin$. It follows that
$z = \beta^{\wfpemin} t$ for
$t := \beta^{d - \wfpemin} r + \beta^{e - \wfpemin} s \in \wz{}$.
We have that
$\wabs{t} \geq 1$ because $t \in \wz{} \setminus \wset{0}$ and
$\wabs{z} = \beta^{\wfpemin} \wabs{t} \geq \beta^{\wfpemin} = \wfpa$,
and $z \in \wfpf$ by Lemma \ref{lemSmallSum}.
\peProof{Lemma}{lemIEEESum}\\

%%%%%%%%%%%%%%%%%%%%%%%%%%%%%%%%%%%%%%%%%%%%%%%%%%%%%%%%%%%%%%%%%%%%%%%%%%%%%%%%%%%
% end of proof of lemma lemIEEESum
%%%%%%%%%%%%%%%%%%%%%%%%%%%%%%%%%%%%%%%%%%%%%%%%%%%%%%%%%%%%%%%%%%%%%%%%%%%%%%%%%%%
%
%

%
%
%%%%%%%%%%%%%%%%%%%%%%%%%%%%%%%%%%%%%%%%%%%%%%%%%%%%%%%%%%%%%%%%%%%%%%%%%%%%%%%%%%%
% Proof of lemma lemNormOneBound
%%%%%%%%%%%%%%%%%%%%%%%%%%%%%%%%%%%%%%%%%%%%%%%%%%%%%%%%%%%%%%%%%%%%%%%%%%%%%%%%%%%

\pbProofB{Lemma}{lemNormOneBound}
This proof illustrates the use of optimization to
bound rounding errors. We define
$z_1 := y_0 + y_1$, $z_{k} := y_k$ for $k > 1$
and use the sums $\wfpsumk{k}$ in Definition \ref{longDefPSum},
the set $\wcal{R}$ of all $n-$tuples which round to nearest
and the function
\pbDef{defEta}
\wfc{\eta}{\wvec{z},\wflmt} := \sum_{k = 1}^n \wabs{\wfpsumkf{k}{\wvec{z},\wflmt} - \wlr{\wfpsumkf{k-1}{\wvec{z},\wflmt} + z_k}},
\peDef{defEta}
from $\wrn{n} \times \wcal{R}$ to $\wrone{}$.
We show that Example \ref{exSharpSum} is the worst case for the ratio
\pbDef{ssNQ}
\wfc{q_n}{\wvec{z},\wflmt} :=
\frac{\wfc{\eta}{\wvec{z},\wflmt}}{\sum_{k = 1}^n \wabs{z_k}}.
\peDef{ssNQ}
This ratio is related to Equation \pRef{normOneBound} because
\[
\wabs{\wflt{\sum_{k = 0}^n y_k} - \sum_{k = 0}^n y_k} =
\wabs{\sum_{k = 1}^n \wlr{\wfpsumkf{k}{\wvec{z},\wflmt} - \wfpsumkf{k-1}{\wvec{z},\wflmt} - z_k}} \leq
\wfc{\eta}{\wvec{z},\wflmt}
\]
and
\[
\wabs{\wflt{\sum_{k = 0}^n y_k} - \sum_{k = 0}^n y_k} \leq \wfc{q_n}{\wvec{z},\wflmt} \sum_{k=1}^n \wabs{z_k} \leq
\wfc{q_n}{\wvec{z},\wflmt} \sum_{k=0}^n \wabs{y_k}.
\]
Therefore, to prove Lemma \ref{lemNormOneBound} it suffices to show that
\pbTClaim{ssNT}
\sup_{\wlr{\wvec{z},\wflmt} \in \wlr{\wrn{n} \setminus \wset{0}} \times \wcal{R}} \wfc{q_n}{\wvec{z},\wflmt} = \theta_n u
\hspace{1cm} \wrm{for} \hspace{1cm}
\theta_n := \frac{n}{1 + n u}.
\peTClaim{ssNT}
The ratio $q_n$ can be written as
$\wfc{q_n}{\wvec{z},\wflmt} = \wfc{f}{\wvec{z},\wfc{h}{\wvec{z},\wflmt}}$
for
\[
\wfc{h}{\wvec{z},\wflmt} := \wlr{\wfpsumkf{0}{\wvec{z},\wflmt}, \wfpsumkf{1}{\wvec{z},\wflmt},\dots,\wfpsumkf{n}{\wvec{z},\wflmt}} \in \wrn{n+1}
\]
and
\[
\wfc{f}{\wvec{z},\wvec{x}} :=
\frac{\sum_{k= 1}^n \wabs{x_k - \wlr{x_{k-1} + z_k}}}{\sum_{k = 1}^n \wabs{z_k}}.
\]
The function $f$ is continuous for $\wvec{z} \neq 0$ and satisfies
$\wfc{f}{\lambda \wvec{z},\lambda \wvec{x}} = \wfc{f}{\wvec{z},\wvec{x}}$ for $\lambda \neq 0$,
and Prop. \ref{propSumScale}, \ref{propWholeIsTight} and \ref{propSumsAreTight}  show that
$h$ satisfies the requirements of Lemma \ref{lemOpt}.
We can then apply this Lemma to prove that either (i) $\sup q_n \leq \theta_n u$ or (ii) $q_n$
has a maximizer $\wlr{\wvec{z}^*,\wflmt^*}$. In case (ii) we use the properties
of this maximizer to prove that it is no worse than what is described in
Example \ref{exSharpSum}.
For instance, this example tells us that
$\wfc{q_n}{\wvec{z}^*,\wflmt^*} \geq \theta_n u$
and if $q$ has a partial derivative with respect
to $z_k$ at $\wvec{z}^*$ then this derivative is zero.

We prove Equation \pRef{ssNT} by induction.
For $n = 1$, this equation follows from Lemma \ref{lemUNear}.
Let us then assume that $n > 1$ and Equation \pRef{ssNT} is valid for $\wvec{z} \in \wrn{m}$ and
rounding tuples $\wflmt = \wset{\wflmk{1},\dots,\wflmk{m}}$ when $m < n$ and
show that it also holds for $n$.
To apply Lemma \ref{lemOpt}, let us
define the numbers
\pbDef{ssA}
a:=  \frac{1 + \wlr{1 + 2u} \theta_{n-1} - \wlr{1 + u} \theta_n}{1 + u}
= \frac{\wlr{n - 1} \wlr{3 + 2 n u} u }{\wlr{1 + n u} \wlr{1 + \wlr{n-1} u} \wlr{1 + u}}  > 0
\peDef{ssA}
(recall that $n \geq 2$) and
\pbDef{ssB}
b := \theta_n - \theta_{n-1} = \frac{1}{\wlr{1 + n u} \wlr{1 + \wlr{n-1} u}} > 0,
\peDef{ssB}
and split $\wrn{n} \setminus \wset{\wvec{0}}$ as the union of the set
\pbDef{ssDefD}
\wcal{B} := \wset{ \wvec{z} \in \wrn{n} \ \ \wrm{with} \ \ b \sum_{k = 2}^n \wabs{z_n} > a \wabs{z_1} \ }
\peDef{ssDefD}
and the cone
$\wcal{A} := \wset{\lambda \wvec{z}, \ \wrm{with} \ \wvec{z} \in \wcal{K}, \ \lambda \in \wrone{} \setminus \wset{0}}$,
for
\pbDef{ssDefK}
\wcal{K} := \wset{ \wvec{z} \in \wrn{n} \ \ \wrm{with} \ \ 2/3 \leq z_1 \leq 2 \beta / 3
\ \ \wrm{and} \ \ b \sum_{k = 2}^n \wabs{z_k} \leq a z_1 \ }.
\peDef{ssDefK}
We claim that $\wfc{q_n}{\wvec{z},\wflmt} \leq \theta_n u$ for $\wvec{z} \in \wcal{B}$
and $\wflmt \in \wcal{R}$.
In fact, writing
$\hat{s}_k := \wfc{S_k}{\wvec{z},\wflmt}$ and $s_k := \sum_{i = 1}^k z_i$
for $k = 0,\dots,n$
and using Equation \pRef{ssNT} with
$\wflmtx := \wset{\wflmk{2}, \dots, \wflmk{n}}$
and $\tilde{\wvec{z}} := \wlr{\hat{s}_1 + z_2, z_3, \dots,z_n}$
we obtain by induction that
\pbDef{ssTemp}
\sum_{k = 2}^n \wabs{\hat{s}_k - \hat{s}_{k-1} - z_k} \leq
\theta_{n-1} u \wlr{ \wabs{\hat{s}_1 + z_2} + \sum_{k = 3}^n \wabs{z_k}}.
\peDef{ssTemp}
Keeping in mind that $z_1 = s_1$, we have that
\[
\wabs{\hat{s}_1 - s_1} + \sum_{k = 2}^n \wabs{\hat{s}_k - \hat{s}_{k-1} - z_k} \leq
\wlr{\wabs{\hat{s}_1 - s_1} + \theta_{n-1} u \wabs{\hat{s}_1}} +
\theta_{n-1} u \sum_{k = 2}^n \wabs{z_k},
\]
and Lemma \ref{lemUNear}, the definitions \pRef{defEta}, \pRef{ssA} and \pRef{ssB} of $\eta$, $a$ and $b$
and $\hat{s}_0 = 0$ yield
\[
\wfc{\eta}{\wvec{z},\wflmt} = \sum_{k = 1}^n \wabs{\hat{s}_k - \hat{s}_{k-1} - z_k} \leq
\wlr{1 + \theta_{n-1} \wlr{1 + 2 u} } \frac{u \wabs{s_1}}{1+ u} + \theta_{n-1} u \sum_{k=2}^n \wabs{z_k}
\]
\[
=
u \wlr{\wlr{1 + \wlr{1 + 2 u} \theta_{n-1} - \wlr{1 + u} \theta_n} \frac{\wabs{z_1}}{1 + u}
- \wlr{\theta_n - \theta_{n-1}} \sum_{z = 2}^n \wabs{z_k}}
+ \theta_n u \sum_{k = 1}^n \wabs{z_k}
\]
\[
= u \wlr{a \wabs{z_1} - b \sum_{k = 2}^n \wabs{z_k}} +
\theta_n u \sum_{k = 1}^n \wabs{z_k}.
\]
The definitions \pRef{ssNQ} and \pRef{ssDefD} of $q$ and $\wcal{B}$ and this equation imply that
$\wfc{q_n}{\wvec{z},\wflmt} \leq \theta_n u$,
and, indeed, $\wfc{q_n}{\wvec{z},\wflmt} \leq \theta_n u$ for $\wvec{z} \in \wcal{B}$
and $\wflmt \in \wcal{R}$.
 As a result, Lemma \ref{lemOpt} shows that either (i) the supremum of
$q_n$ is at most $\theta_n u$ or (ii) there exists $\wvec{z}^* \in \wcal{K}$ and
$\wflmt{}^* \in \wcal{R}$ with
\[
\wfc{q_n}{\wvec{z}^*,\wflmt{}^*} =
\sup_{\wlr{\wvec{z},\wflmt{}} \in \wlr{\wrn{n} \setminus \wset{\wvec{0}}}  \times  \wcal{R}}
\wfc{q_n}{\wvec{z},\wflmt{}^*}.
\]
In case (i) we are done and we now analyze case (ii).
Let us define $\hat{s}^*_k := \wfpsumkf{k}{\wvec{z}^*,\wflmt^*}$, and $s^*_k := \sum_{i = 1}^k z^*_k$,
for $k = 0,\dots, n$. Since $\wvec{z}^* \in \wcal{K}$, the
definitions of $a$ and $b$ lead to
\[
\sum_{k=2}^n \wabs{z_k^*} \leq \frac{\wlr{n - 1} \wlr{3 + 2 n u}}{1 + u} u z^*_1.
\]
Using Lemma \ref{lemUNear}, the hypothesis $20 n u \leq 1$ and induction we deduce that
\[
\wabs{\hat{s}^*_k - z^*_1} \leq \wabs{\hat{s}^*_k - \wlr{ \wlr{\hat{s}^*_1 + z^*_2} + \sum_{i=3}^n z^*_i}}
+ \wabs{\hat{s}^*_1 - s^*_1} + \sum_{i = 2}^n \wabs{z^*_i}
\]
\[
\leq \frac{\wlr{n - 1} u}{1 + \wlr{n-1} u} \wlr{\wabs{\hat{s}_1^* + z_2^*} + \sum_{i=3}^n \wabs{z_i^*}}
+ \frac{u}{1 + u} z^*_1 + \frac{\wlr{n - 1} \wlr{3 + 2 n u}}{1 + u} u z^*_1
\]
\[
\leq
\wlr{\frac{n - 1}{1 + \wlr{n-1} u} \wlr{1 + 2 u + \wlr{n - 1} \wlr{3 + 2 n u} u}
+ 1 + \wlr{n - 1} \wlr{3 + 2 n u}} \frac{u}{1 + u} z^*_1
\]
and, since $s^*_1 = z^*_1$ and $20 n u \leq 1$,
\pbDef{ssSk}
\wabs{\hat{s}_k^* - z_1^*} \leq \kappa n u z_1^* \leq \kappa z_1^*/20,
\peDef{ssSk}
for
\[
\kappa :=
\wlr{\frac{1}{1 + n u} \wlr{1 +  \wlr{3 + 2 n u} n u} + 3 + 2 n u} \frac{1}{1 + u}
\]
\pbDef{ssTheta}
\leq
 \frac{1}{1 + \frac{1}{20}} \wlr{1 + \frac{1}{20} \wlr{3 + \frac{1}{10}}}+ 3 + \frac{1}{10} = \frac{21}{5}.
\peDef{ssTheta}
Since $2/3 \leq z^*_1 \leq 2 \beta /3$, Equations \pRef{ssSk} and \pRef{ssTheta} lead to
\[
\frac{1}{\beta} \leq \frac{1}{2} < \frac{158}{300} \leq \frac{79}{100} z^*_1 \leq \hat{s}^*_k \leq \frac{121}{100} z^*_1 \leq
\frac{121}{150}\beta < \beta
\]
for $1 < k \leq n$, and since $\hat{s}_1^* = \wflk{1}{z^*_1}$ and
$2/3 \leq z^*_1 \leq 2/3 \beta$ this equation also holds for $k = 1$.
Monotonicity (Prop. \ref{propRoundMonotone})
and the fact that $\hat{s}^*_k = \wflk{k}{\hat{s}_{k-1} + z_k}$ lead to
\pbDef{ssNZS}
1/\beta < \hat{s}_{k-1} + s_k^* < \beta \ \ \wrm{for} \ \ 1 \leq k \leq n.
\peDef{ssNZS}

We now explore the implications of $\wlr{\wvec{z}^*,\wflmt^*}$ being a
maximizer of $q_n$. Example \ref{exSharpSum} shows
that $\wfc{q_n}{\wvec{z}^*,\wflmt{}^*} \geq \theta_n u$
and this implies that $z_k \neq 0$ for all $k$, because if $z_k = 0$ for some
$k$ then we would have
$\wfc{q_n}{\wvec{z}^*,\wflmt{}^*} = \wfc{q_{n-1}}{\tilde{\wvec{z}},\wflmtx}$
for $\tilde{\wvec{z}} \in \wrn{n-1}$ and $\wflmtx$ obtained by removing the $k$th
coordinate of $\wvec{z}^*$  and $\wflmk{k}$ from $\wflmt^*$, and
$\wfc{q_{n-1}}{\wflmtx,\tilde{\wvec{z}}} \leq \theta_{n-1} u < \theta_n u$,
contradicting the maximality of $\wlr{\wvec{z}^*,\wflmt^*}$.
Therefore, $z^*_k \neq 0$ for $k = 1,\dots, n$,
and the denominator of $q_n$ has non zero partial derivatives at $\wvec{z}^*$.
Equation \pRef{ssNZS} shows that $\hat{s}^*_{k-1} + z^*_k \neq 0$,
and Prop. \ref{propFlat} implies that
the numerator of $q_n$ will have a zero
partial derivative with respect to $z_k$ if
$\hat{s}^*_{k-1} + z^*_k$ is not of the form
\pbDef{ssMid}
\hat{s}^*_{k-1} + z^*_k = \beta^{e_k} \wlr{\beta^{\mu} + r_k + 1/2}
\hspace{0.7cm} \wrm{with} \hspace{0.7cm} e_k \in \wz{}
\ \ \wrm{and} \ \
r_k \in [0, \wlr{\beta - 1} \beta^{\mu}),
\peDef{ssMid}
and this would imply that the derivative of $q_n$
is well defined and different from zero.
By the maximality of $\wlr{\wvec{z}^*,\wflmt^*}$, we conclude that
Equation \pRef{ssMid} is valid.
Combining this equation with Equation \pRef{ssNZS}
we conclude that
we can write $\wset{1,2,\dots,n} = \wcal{L} \cup \cal{H}$
(for low and high) so that
the exponents in $e_k$ Equation \pRef{ssMid} are
$e_k = - \mu -1$ for $k \in \wcal{L}$ and
$e_k = -\mu$ for $k \in \wcal{H}$. Since $\beta^{-\mu}/2 = u$, this leads to
\begin{eqnarray}
\nonumber
k \in \wcal{L}\ & \Rightarrow & \  \frac{1 + u}{\beta} \leq \hat{s}^*_{k-1} + z^*_k \leq \frac{\beta - u}{\beta},\\
\nonumber
k \in \wcal{U}\ & \Rightarrow & \ 1 + u \leq \hat{s}^*_{k-1} + z^*_k \leq \beta - u,
\end{eqnarray}
As a result, Prop. \ref{propCriticalSum} implies that
\pbDef{ssZlh}
k \in \wcal{L} \Rightarrow \wabs{z^*_k}\geq u/\beta%
\hspace{1cm} \wrm{and} \hspace{1cm}
k \in \wcal{H} \Rightarrow \wabs{z^*_k}\geq u,
\peDef{ssZlh}
and Prop. \ref{propRoundNormal} yields
\pbDef{ssEk}
k \in \wcal{L} \Rightarrow \wabs{\hat{s}_k^* - \wlr{\hat{s}_{k-1}^* + z^*_k}} = u/\beta
\hspace{0.5cm} \wrm{and} \hspace{0.5cm}
k \in \wcal{H} \Rightarrow \wabs{\hat{s}_k^* - \wlr{\hat{s}_{k-1}^* + z^*_k}} = u.
\peDef{ssEk}

We now show that if $1 \in \wcal{L}$ then we
obtain a contradiction to the maximality of $\wlr{\wvec{z}^*,\wflmt^*}$. Indeed,
 let $m \in [1,n]$
be the last index such that $k \in \wcal{L}$ for $1 \leq k \leq m$.
If $m = n$ then $k \in \wcal{L}$ for all $k \in [1,n]$ and
the inequality $z_1^* \geq 2 / 3$ and Equations \pRef{ssZlh} and \pRef{ssEk}
and the fact that $2 \beta / 3 - u > 1$ imply that
\[
\wfc{q_n}{\wvec{z}^*,\wflmt^*} / \wlr{\theta_n u} =
\frac{\frac{n u / \beta}{2/3 + \wlr{n - 1} u / \beta}}{\frac{n u}{1 + n u}}
=
\frac{1 + n u}{\wlr{\frac{2 \beta}{3} - u} + n u} < 1,
\]
and this contradicts the maximality of $\wlr{\wvec{z}^*,\wflmt^*}$.
For  $m < n$ we have
\[
\sum_{k = 1}^{m} \wabs{z^*_k} \geq
\sum_{k = 1}^{m} z^*_k  =
\wlr{\hat{s}^*_{m} + z^*_{m+1}} -
\wlr{\sum_{k=1}^{m} \wlr{\hat{s}^*_{k} - \wlr{\hat{s}^*_{k-1} + z^*_k}}} - z^*_{m+1}
\]
\[
\geq \wlr{1 + u} - \wlr{m u / \beta} - \wabs{z^*_{m+1}}.
\]
Let $\ell$ be the size of $\wcal{L}$ and $h$ the size of $\wcal{H}$.
Equations \pRef{ssZlh} and \pRef{ssEk}, the identity $n = \ell + h$ and the hypothesis $20 n u \leq 1$ lead to
\[
\wfc{q_n}{\wvec{z}^*,\wflmt^*} - \theta_n u
\leq \frac{\ell u / \beta + h u}{1 + u - m u / \beta - \wabs{z^*_{m+1}} + (\ell - m) u / \beta + \wabs{z^*_{m+1}} + \wlr{h-1} u}
- \frac{n u}{1 + n u}
\]
\[
= - u \frac{\xi} {\wlr{1 + n u}\wlr{\beta - 2 m u + \ell u + \beta h u}},
\]
for
\[
\xi := \wlr{\beta - 1} \ell - 2 h m u - 2 \ell m u
= \ell \wlr{\wlr{\beta - 1} - \wlr{\frac{m}{\ell}} \wlr{2 h u} - 2 m u}
\geq 0.8 \ell > 0,
\]
and, again, $\wfc{q_n}{\wvec{z}^*,\wflmt^*} < \theta_n u$.
Therefore, by the maximality of $\wlr{\wflmt^*,\wvec{z}^*}$ we must have $z_1^* \geq 1$, and Equation
\pRef{ssMid} shows that $z_1^* \geq 1 + u$ and Equations \pRef{ssZlh} and \pRef{ssEk}
lead to
\pbDef{nobL}
\wfc{q_n}{\wvec{z}^*,\wflmt^*} \leq \frac{\ell u / \beta + h u}{1 + u + \ell u / \beta +\wlr{h - 1} u}
= \frac{\ell u / \beta + h u}{1 + \ell u / \beta + h u}.
\peDef{nobL}
Since $n = \ell + h$, $\theta_n = \wlr{\ell + h}/\wlr{1 + \wlr{\ell + h} u}$ and
\[
\frac{\ell + h}{1 + \wlr{\ell + h} u} - \frac{\ell / \beta + h}{1 + \ell u / \beta + h u}
=  \frac{\wlr{\beta - 1} \ell}{\wlr{1 + \wlr{\ell + h} u}\wlr{\beta + \ell u + \beta h u}}  \geq 0,
\]
Equation \pRef{nobL} implies
that $\wfc{q_n}{\wvec{z}^*,\wflmt^*} \leq \theta_n u$ and we are done.
\peProof{Lemma}{lemNormOneBound}\\

%%%%%%%%%%%%%%%%%%%%%%%%%%%%%%%%%%%%%%%%%%%%%%%%%%%%%%%%%%%%%%%%%%%%%%%%%%%%%%%%%%%
% End of the Proof of lemma lemNormOneBound
%%%%%%%%%%%%%%%%%%%%%%%%%%%%%%%%%%%%%%%%%%%%%%%%%%%%%%%%%%%%%%%%%%%%%%%%%%%%%%%%%%%
%
%

%
%
%%%%%%%%%%%%%%%%%%%%%%%%%%%%%%%%%%%%%%%%%%%%%%%%%%%%%%%%%%%%%%%%%%%%%%%%%%%%%%%%%%%
% Proof of lemma lemPositiveBound
%%%%%%%%%%%%%%%%%%%%%%%%%%%%%%%%%%%%%%%%%%%%%%%%%%%%%%%%%%%%%%%%%%%%%%%%%%%%%%%%%%%

\pbProofB{Lemma}{lemPositiveBound}
Let us define $z_1 := y_0 + y_1$ and $z_{k} := y_k$ for $k > 1$.
Using Lemma \ref{lemUNear} and induction in $n$ we can show that
\[
\wfpsumkf{n}{\wvec{z},\wflmt} \geq \sum_{k=1}^n \wlr{1 + u}^{-\wlr{n - k + 1}} z_k
= \frac{1}{1 + u} \sum_{k = 1}^{n} \wlr{1 + u}^{-\wlr{n - k}} z_k.
\]
The convexity of the
functions $\wlr{1 + u}^{-\wlr{n - k}}$, which have value $1$ and
derivative $- \wlr{n - k}$ at $u = 0$, lead to
\[
\wfpsumkf{n}{\wvec{z},\wflmt} \geq
 \frac{1}{1 + u} \wlr{\sum_{k=1}^n z_k - u \sum_{k = 1}^{n} \wlr{n - k} z_k}
\]
\[
= \frac{1}{1 + u} \wlr{\wlr{1 + u} \sum_{k=1}^n z_k - u \sum_{k = 1}^{n} \wlr{n - k + 1} z_k}
= \sum_{k=1}^n z_k  -  \frac{u}{1+u} \sum_{k=1}^n \wlr{n - k + 1} z_k,
\]
and the lower bound in Equation \pRef{thPositiveBound} follows from the identities
\[
\sum_{i = 0}^k y_i = \sum_{i=1}^k z_i,
\hspace{0.6cm}
\wflt{\sum_{k=0}^{n} y_k} = \wfpsumkf{n}{\wvec{z},\wflmt}
\hspace{0.6cm} \wrm{and} \hspace{0.6cm}
\sum_{k = 1}^n \sum_{i = 0}^k y_i = \sum_{k = 1}^n \wlr{n - k + 1} z_k.
\]

In order to prove the second inequality in Equation \pRef{thPositiveBound},
we proceed as in the proof of Lemma \ref{lemNormOneBound}
(We ask the reader to look at the first two paragraphs
of that proof.) This time we consider only the rounding tuple
$\wflmt := \wset{\wflm,\dots,\wflm}$ where $\wflm$ rounds to nearest
and breaks all ties upward, because our function
\pbDef{snbQ}
\wfc{q_n}{\wvec{z}} :=
\frac{\wfc{\eta}{\wvec{z}}}{\sum_{k = 1}^n \wlr{n - k + 1} z_k}
\peDef{snbQ}
for
\pbDef{ssnbNum}
\wfc{\eta}{\wvec{z}} := \wfpsumkf{n}{\wvec{z},\wflmt} - \sum_{k=1}^n z_k
= \sum_{k=1}^n \wlr{\wfpsumkf{k}{\wvec{z},\wflmt} - \wfpsumkf{k-1}{\wvec{z},\wflmt} - z_k}
\peDef{ssnbNum}
is clearly maximized by the rounding tuple $\wflmt$ for which all ties
are broken upward.

We prove by induction that
\pbDef{snbUpper}
\wfc{q_n}{\wvec{z}} \leq \tau_n u := \frac{u}{1 + u \wlr{\frac{\beta - 2}{\beta - 1} + \frac{n}{\beta^{n} - 1}}}.
\peDef{snbUpper}
For $n = 1$ Equation \pRef{snbUpper} follows
from Lemma \ref{lemUNear}. Let us then assume it holds for $n - 1$ and prove
it for $n$ using Lemma \ref{lemOpt} to show that either Equation \pRef{snbUpper}
holds or there exists a maximizer for $q_n$, which we then analyze.
With this purpose, define
\pbDef{snbCAB}
a := \frac{1 + \wlr{n - 1} \wlr{1 + 2 u} \tau_{n - 1} - n \wlr{1 + u} \tau_{n}}{1 + u}
\hspace{0.5cm} \wrm{and} \hspace{0.5cm}
b := \tau_{n} - \tau_{n-1}.
\peDef{snbCAB}
In order to prove that $a$ and $b$ are positive, note that
\[
\tau_n = \frac{1}{1 + u \phi_n}
\hspace{0.5cm} \wrm{and} \hspace{0.5cm}
\tau_{n-1} := \frac{1}{1 + u \wlr{\phi_n + \delta_n}}
\hspace{0.5cm} \wrm{for} \hspace{0.5cm}
\phi_n := \frac{\beta - 2}{\beta - 1} + \frac{n}{\beta^{n} - 1}
\]
and
\[
\delta_n := \frac{n-1}{\beta^{n-1} - 1} - \frac{n}{\beta^n - 1}
= \frac{n \wlr{\beta - 1} - \wlr{\beta - \beta^{1 - n}}}{\beta^{n} \wlr{1 - \beta^{1 - n}}\wlr{1 - \beta^{-n}}} > 0.
\]
For $\beta,n \geq 2$ we have that $\delta_n > 0$, and
the positivity of $\delta_n$ implies that
\[
b = \tau_n - \tau_{n-1} = u \delta_n \tau_{n-1} \tau_n > 0.
\]
For $n = 2$ the software Mathematica shows that
\[
a = u \frac{\beta - 1 + u \wlr{\beta - 2}}{\wlr{1 + u}^2 \wlr{\beta + 1 + \beta u}} > 0.
\]
Mathematica also shows that when $\beta = 2$
\[
a = u \frac{\wlr{2^{n} \wlr{n - 2} + 2}\wlr{2^n - n - 1}}{\wlr{1 + u} \wlr{2^n - 1 + n u} \wlr{2^n - 2 + 2 \wlr{n - 1} u}},
\]
which is positive for $n \geq 3$.
For $\beta = 3$ we have
\[
a = \frac{u \wlr{u + 2} \wlr{3^n - 2 n - 1} \wlr{\wlr{2 n - 3} 3^{n} + 3}}
{\wlr{1 + u} \wlr{2 \wlr{3^n - 1} + u \wlr{3^n + 2 n -1}}\wlr{2 \times 3^{n} - 6 + u \wlr{3^n + 6 n - 9 } }},
\]
which is also positive for $n \geq 3$.
Finally, for $\beta \geq 4$ and $n \geq 3$
\[
n \delta_n \leq  n \frac{\wlr{n - 1} \wlr{\beta - 1}}{\wlr{1 - \beta^{1 - n}} \wlr{1 - \beta^{-n}}} \beta^{-n}
\leq \frac{3 \times 2 \times 4^{-3}}{\wlr{1 - 4^{-2}}\wlr{1 - 4^{-3}}} = \frac{32}{315} < 0.2,
\]
and the software Mathematica also shows that
\[
a = \frac{\wlr{n - 3} + \wlr{1 - \wlr{n - 1 + n u} \delta_n} + \phi_n \wlr{1 + \wlr{\delta_n + \phi_n + n  - 2}u} }
{\wlr{1 + u} \wlr{1 + u \phi_n} \wlr{1 + u \wlr{\phi_n + \delta_n}}}
\]
and this number is positive for $n \geq 3$ because $\wlr{n - 1 + n u} \delta_n \leq n \delta_n \leq 0.2$.
Therefore, $a$ and $b$ are positive and the set
\[
\wcal{K} := \wset{\wvec{z} \in \wrn{n}\setminus \wset{0} \ \wrm{with} \ 2 / 3 \leq z_1 \leq 2 \beta / 3, \ \
z_k \geq 0 \ \ \wrm{and} \
b \sum_{k = 2}^{n} \wlr{n - k + 1} z_k \leq a \, z_1}
\]
is compact. We now split $\wset{\wvec{x} \in \wrn{n} \setminus \wset{\wvec{0}} \ \wrm{with} \ x_k \geq 0 }$
as the union of the set
\pbDef{snbD}
\wcal{B} := \wset{\wvec{z} \in \wrn{n} \ \wrm{with} \ \ z_k \geq 0 \ \
\ \wrm{and} \ \ b \sum_{k = 2}^{n} \wlr{n - k + 1} z_k > a\, z_1}
\peDef{snbD}
and the cone
\[
\wcal{A} :=
\wset{\lambda \wvec{x}  \ \wrm{with} \ \wvec{x} \in \wcal{K} \ \wrm{and} \ \lambda \in \wrone{}, \ \lambda > 0}
\]
and show that
$\wfc{q_n}{\wvec{z}} \leq \tau_n u $ for $\wvec{z} \in \wcal{B}$. In fact, for
such $\wvec{z}$, let us write
 $\hat{s}_k := \wfc{S_k}{\wvec{z},\wflmt}$
for $k = 0,\dots,n$.
Using induction, Lemma \ref{lemUNear},
and the definitions of $a$, $b$,
and keeping in mind that $s_1 = z_1$, we deduce that
\[
\sum_{k = 1}^n \wlr{\hat{s}_k - \wlr{\hat{s}_{k-1} + z_k}}
= \wlr{\hat{s}_{1} - s_{1}} +
\wlr{\wlr{\hat{s}_2 - \wlr{\hat{s}_1 + z_2}} + \sum_{k = 3}^n \wlr{\hat{s}_k - \wlr{s_{k-1} - z_k}}}
\]
\[
\leq
\frac{u}{1 + u} z_1 +
 \tau_{n-1} u \wlr{ \wlr{n-1} \wlr{\hat{s}_1 + z_2} + \sum_{k = 3}^n \wlr{n - k + 1} z_k}
\]
\[
= \frac{u}{1 + u} z_1 +
\tau_{n-1} u \wlr{ \wlr{n-1} \hat{s}_1 + \sum_{k = 2}^n \wlr{n - k + 1} z_k}
\]
\[
\leq
\wlr{1 + \wlr{1 + 2 u} \wlr{n-1} \tau_{n-1}} \frac{u z_1}{1 + u}
+  \tau_{n-1} u \sum_{k = 2}^n \wlr{n - k + 1} z_k
\]
\[
= \wlr{1 + \wlr{1 + 2 u} \wlr{n-1}  \tau_{n-1} - n \wlr{1 + u} \tau_{n}}
\frac{u z_1}{1 + u}
\]
\[
- \wlr{\tau_{n} - \tau_{n-1}} u \sum_{k = 2}^n \wlr{n - k + 1 } z_k
+  \tau_{n} u \wlr{n z_1 + \sum_{k = 2}^n \wlr{n - k + 1 } z_k},
\]
and it follows that
\[
\wfc{\eta}{\wvec{z}} \leq
\wlr{a z_1 - b
\sum_{k = 2}^n \wlr{n - k + 1} z_k} u + \tau_{n} u \sum_{k = 1}^n \wlr{n - k + 1 } z_k.
\]
By the definition of $\wcal{B}$ the term in parenthesis above is
negative and this equation shows that $\wfc{q_n}{\wvec{z}} \leq \tau_{n} u$ for
$\wvec{z} \in \wcal{B}$.
According to Lemma \ref{lemOpt} we have that either
(i) Equation \pRef{snbUpper} holds or (ii) $q_n$ has
a maximizer $\wvec{z}^* \in \wcal{K}$.
In case (i) we are done and we
now suppose that there exists such $\wvec{z}^{*}$.
Define $\hat{s}^*_k := \wfc{S_k}{\wvec{z}^*,\wflmt}$ for $k = 0,\dots,n$.
The same argument used in the proof of Lemma
\ref{lemNormOneBound} to deduce that $z_k^* \neq 0$ and Equation
\pRef{ssMid} shows that $z^*_k \neq 0$  for $k = 1,\dots, n$,
and
\pbDef{snbDecomp}
\hat{s}^*_{k-1} + z^*_k = \beta^{d_k} \wlr{\beta^{\mu} + r_k + 1/2}
\hspace{0.3cm} \wrm{with} \hspace{0.3cm}
d_k \in \wz{}
\hspace{0.3cm} \wrm{and} \hspace{0.3cm}
r_k \in [0,\wlr{\beta - 1} \beta^{\mu}) \cap \wz{}.
\peDef{snbDecomp}
Since $\wflm$ break ties upward, we have that
\pbDef{snbRUp}
\hat{s}^*_{k} = \beta^{d_k} \wlr{\beta^{\mu} + r_k + 1},
\peDef{snbRUp}

If the $r_k$ in Equation \pRef{snbDecomp} were all zero then,
since $\hat{s}_{0}^* = 0$ and
\[
\frac{1}{\beta} < 2/3 \leq z_1^* \leq \frac{2 \beta}{3} < \beta,
\]
Equation \pRef{snbDecomp} would yield $z_1^* = \beta^{-\mu}\wlr{\beta^\mu + 1/2} = 1 + u$ and,
for $k > 1$, Equations \pRef{snbDecomp} and \pRef{snbRUp} would
lead to $\hat{z}^*_k = \beta^{d_k + \mu}\wlr{1 + u} - \beta^{d_{k-1} + \mu} \wlr{1 + 2 u}$
and the $z_k^*$ would correspond to the $x_k$ in Example \ref{exSharpSumNearC} with
$e_k = d_k + \mu$ (take $x_0 = 0$ and $x_1 = z^*_1$). Therefore, by
the last line in the statement of Example \ref{exSharpSumNearC}, in order to
complete this proof it suffices to show that $r_k = 0$ for all $k$, and
this is what we do next.

We start with $k < n$ and after that we
handle the case $k = n$. Let us define $r_{0} := 0$, assume that
$r_{i} = 0$ for $i < k < n$ and show that $r_k = 0$.
Take  $\delta_k := \min\wset{1,r_k}$ and $\wvec{z}' \in \wrn{n}$
given by $z_i' := z_i^*$ for $i < k$ or $i > k + 1$ and
\[
z_k' := z_k^* - \beta^{d_k} \delta_k
\hspace{1cm} \wrm{and} \hspace{1cm}
z_{k+1}' := z_{k+1}^* + \beta^{d_k} \delta_k.
\]
We now prove that $\delta_k = 0$ by showing that $\wvec{z}' = \wvec{z}^*$.
If $\delta_k = 0$ then $\wvec{z}' = \wvec{z}^*$
and $\wvec{z}$ is in the domain of $q_n$.
If $\delta = 1$ then $z'_{k+1} > 0$ and showing that
$z_k' \geq 0$ suffices to prove that $\wvec{z}'$ is in the
domain of $q_n$. Indeed,
Equations \pRef{snbDecomp} and \pRef{snbRUp} and $r_{k-1} = 0$ lead to
\[
z_k' := \beta^{d_k} \wlr{\beta^\mu + r_k + 1/2} - \beta^{d_{k-1}} \wlr{\beta^\mu + 1} - \beta^{d_k} \delta_k
\]
\[
= \beta^{d_k} \wlr{\beta^\mu + \wlr{r_k - \delta_k} + 1/2} - \beta^{d_{k-1}} \wlr{\beta^\mu + 1}.
\]
Equations \pRef{snbDecomp}, \pRef{snbRUp} and $z_k^* \geq 0$ imply that
\[
\hat{s}_{k}^* = \wfl{\hat{s}_{k-1}^* + z_k} \geq \hat{s}_{k-1}^* + \beta^{d_k}/2 > \hat{s}_{k-1}^*,
\]
Prop. \ref{propOrder} leads to $d_k \geq d_{k-1}$. Moreover, $\delta_k \leq r_k$ by definition
and it follows that if $d_k > d_{k-1}$ then $\beta^{d_k}/2 \geq \beta^{d_{k-1}}$ and $z_k' \geq 0$.
If $d_k = d_{k-1}$ then $\hat{s}^*_{k} > \hat{s}^*_{k-1}$
implies that $\beta^{\mu} + r_k + 1 > \beta^{\mu} + 1$, $r_k > 1$,
and $r_k - \delta_k \geq 1$ and $z_k' \geq 0$. Therefore, $\wvec{z}'$
is on the domain of $q_n$.

We now analyze $\eta$ defined in Equation \pRef{ssnbNum} and
show that all parcels in $\wfc{\eta}{\wvec{z}^*}$ and
$\wfc{\eta}{\wvec{z}'}$ are equal.
Since we break ties upward, Equation \pRef{snbRUp} shows that
\[
\wfl{\hat{s}^*_{k-1} + z_k'}  = \wfl{\beta^{d_k} \wlr{\beta^{\mu} + \wlr{r_k - \delta_k} + 1/2}} =
\beta^{d_k} \wlr{\beta^{\mu} + \wlr{r_k - \delta_k} + 1}
\]
\pbDef{ssNBA}
= \beta^{d_k} \wlr{\beta^{\mu} + r_k + 1}  - \beta^{d_k} \delta_k
= \wfl{\hat{s}^*_{k-1} + z^*_k} - \beta^{d_k} \delta_k
= \hat{s}^*_k - \beta^{d_k} \delta_k.
\peDef{ssNBA}
It follows that
\[
\wfpsumkf{k}{\wvec{z}',\wflmt} + z_{k+1}' =
\wfl{\hat{s}^*_{k-1} + z_k'} + z_{k+1}' =
\]
\[
\wlr{\hat{s}^*_k - \beta^{d_k} \delta_k}  +  \wlr{z_{k+1}^* + \beta^{d_k} \delta_k} =
\hat{s}^*_k + z_{k+1}^* = \wfpsumkf{k}{\wvec{z}^*,\wflmt} + z_{k+1}^*.
\]
This equation leads to
\[
\wfpsumkf{k+1}{\wvec{z}',\wflmt} =
\wfl{ \wfpsumkf{k}{\wvec{z}',\wflmt} + z_{k+1}'}  =
\wfl{\hat{s}^*_k + z_{k+1}^*}  =
\wfpsumkf{k+1}{\wvec{z}^*,\wflmt},
\]
and
\[
\wfpsumkf{k+1}{\wvec{z}',\wflmt} - \wlr{\wfc{S_{k}}{\wvec{z}',\wflmt} + z_{k+1}'} =
\wfpsumkf{k+1}{\wvec{z}^*,\wflmt} - \wlr{\wfc{S_{k}}{\wvec{z}^*,\wflmt} + z_{k+1}^*}.
\]
Therefore, $\wfpsumkf{i}{\wvec{z}',\wflmt} = \wfpsumkf{i}{\wvec{z}^*,\wflmt}$
for $i < k$ and $i \geq k + 1$. It follows that
\[
\wfpsumkf{i}{\wvec{z}',\wflmt} - \wlr{\wfc{S_{i-1}}{\wvec{z}',\wflmt} + z_i'} =
\wfpsumkf{i}{\wvec{z}^*,\wflmt} - \wlr{\wfc{S_{i-1}}{\wvec{z}^*,\wflmt} + z^*_i}
\]
for $i <k$ and $i \geq k + 1$. For $i = k$, the definition $z_k' := z_k^* - \beta^{d_k} \delta_k$
and Equation \pRef{ssNBA} yield
\[
\wfpsumkf{k}{\wvec{z}',\wflmt} - \wlr{\wfpsumkf{k-1}{\wvec{z}',\wflmt} + z_k'} =
\wfl{\hat{s}^*_{k-1} + z_k'} - \wlr{\hat{s}^*_{k-1} + z_k'} =
\]
\[
= \wlr{\wfl{\hat{s}^*_{k-1} + z_k^*} - \beta^{d_k} \delta_k}
- \wlr{\hat{s}^*_{k-1} + z_k^* - \beta^{d_k} \delta_k}
\]
\[
= \wfpsumkf{k}{\wvec{z}^*,\wflmt} - \wlr{\wfpsumkf{k-1}{\wvec{z}^*,\wflmt} + z^*_k},
\]
Therefore, all parcels in the numerators $\eta$ in Equation \pRef{ssnbNum}
are equal for $\wvec{z}^*$ and $\wvec{z}'$.

Let us now analyze the denominator $D_n$ of $q_n$. Note that
\[
\wlr{n - k + 1} z_k' + \wlr{n - k} z_{k+1}' =
\]
\[
\wlr{n - k + 1} \wlr{z_k^* - \beta^{d_k} \delta_k}
+ \wlr{n - k} \wlr{z_{k+1}^* + z_k^* + \beta^{d_k} \delta_k}
\]
\[
= \wlr{n - k + 1} z_k^* + \wlr{n - k} z^*_{k+1} - \beta^{d_k} \delta_k.
\]
Moreover, $z_i' = z_i^*$ for $i \not \in \wset{k,k+1}$ and
\[
\wfc{D_n}{\wvec{z}'} -
\wfc{D_n}{\wvec{z}^*} =
\wlr{\sum_{i = 1}^n \wlr{n - i - 1} z_i'} - \wlr{\sum_{i = 1}^n \wlr{n - i - 1} z_i^*} =
\]
\[
\wlr{\wlr{n - k - 1} z_k' + \wlr{n - k} z_{k+1}'}
-
\wlr{\wlr{n - k - 1} z_k^* + \wlr{n - k} z_{k+1}^*} = -\beta^{d_k} \delta_k.
\]
Since the numerators of $\wfc{q_n}{\wvec{z}'}$ and
$\wfc{q_n}{\wvec{z}^*}$ are equal and $\wvec{z}^*$ is maximal
this equation implies that $\beta^{d_k} \delta_k \leq 0$.
Therefore, $\delta_k = \min \wset{1,r_k} = 0$, and $r_k = 0$.

Finally, for $k = n$, define
$\wvec{z}'$ with $z_k' = z_k^*$ for $k < n$ and
$z_n' = z_n^* - \beta^{d_n} r_n$. As before,
$\wvec{z}'$ is in the domain of $q_n$
and
$\wfpsumkf{k}{\wvec{z}',\wflmt} = \wfpsumkf{k}{\wvec{z}^*,\wflmt}$
for $k < n$. For $k = n$, Equation \pRef{snbDecomp} leads to
\[
\wfpsumkf{n-1}{\wvec{z}',\wflmt}  + z_n' =
\hat{s}^*_{n-1} + z_n^* - \beta^{d_n} r_n
= \beta^{d_n} \wlr{\beta^\mu + 1/2}.
\]
We break ties upward, $\wfpsumkf{n}{\wvec{z}',\wflmt} = \wfl{\wfpsumkf{n-1}{\wvec{z}',\wflmt}  + z_n'} =
\beta^{d_n} \wlr{\beta^\mu + 1}$ and
\[
\wfpsumkf{n}{\wvec{z}',\wflmt} - \wlr{\wfpsumkf{n- 1}{\wvec{z}',\wflmt} +  z_n'} =
\beta^{d_n} \wlr{\beta^\mu + 1} - \beta^{d_n} \wlr{\beta^\mu + 1/2} =
\]
\[
\beta^{d_n} / 2 = \beta^{d_n} \wlr{\beta^\mu + r_n + 1} - \beta^{d_n} \wlr{\beta^\mu + r_n + 1/2}
\]
\[
= \wfpsumkf{n}{\wvec{z}^*,\wflmt} - \wlr{\wfpsumkf{n- 1}{\wvec{z}^*,\wflmt} +  z_n^*},
\]
and the numerator of $q_n$ in \pRef{ssnbNum} would not change if were to
replace $\wvec{z}^*$ by $\wvec{z}'$. However, the denominator would be reduced by $\beta^{d_n} r_n$,
and this would contradict the maximality of $\wvec{z}^*$. Therefore $r_n = 0$.
In summary, $r_k = 0$ for all $k$, the $z_k^*$ are as the $x_k$
in Example \ref{exSharpSumNearC} and we are done.
\peProof{Lemma}{lemPositiveBound} \\

%%%%%%%%%%%%%%%%%%%%%%%%%%%%%%%%%%%%%%%%%%%%%%%%%%%%%%%%%%%%%%%%%%%%%%%%%%%%%%%%%%%
% End of the Proof of lemma lemPositiveBound
%%%%%%%%%%%%%%%%%%%%%%%%%%%%%%%%%%%%%%%%%%%%%%%%%%%%%%%%%%%%%%%%%%%%%%%%%%%%%%%%%%%
%
%

%
%
%%%%%%%%%%%%%%%%%%%%%%%%%%%%%%%%%%%%%%%%%%%%%%%%%%%%%%%%%%%%%%%%%%%%%%%%%
% Proof of lemma lemSignedSum
%%%%%%%%%%%%%%%%%%%%%%%%%%%%%%%%%%%%%%%%%%%%%%%%%%%%%%%%%%%%%%%%%%%%%%%%%

\pbProofB{Lemma}{lemSignedSum}
Let us write $z_1 := y_0 + y_1$, $z_k := y_k$ for $k > 1$,
$s_k := \sum_{i = 1}^k z_i$ and $\hat{s}_k = \wfpsumkf{k}{\wvec{z},\wflmt}$ for $k = 0,\dots,n$.
We prove by induction that
\pbDef{sspA}
\wabs{\hat{s}_n - s_n} \leq \frac{u}{1 - \wlr{n - 2} u} \sum_{k = 1}^n  \wabs{\sum_{i =1}^n z_i},
\peDef{sspA}
which is equivalent to Equation \pRef{thSignedSum}.
For $n = 1$, Equation \pRef{sspA} follows from Lemma \ref{lemUNear}.
We now prove Equation \pRef{sspA} for $n \geq 2$, assuming that it holds for $n - 1$.
For $\wvec{w} \in \wrn{n-1}$ with
$w_1 = \hat{s}_1 + z_2$ and $w_k = y_{k+1}$ for $k > 1$, we obtain by induction that
$\wfpsumkf{k}{\wvec{w},\wflmtx} = \hat{s}_{k+1}$ for $\wflmtx = \wset{\wflmk{2},\dots,\wflmk{n}}$,
\[
\wabs{\hat{s}_n - \wlr{\hat{s}_1 + z_2} - \sum_{k=3}^n z_k} \leq
\frac{u}{1 - \wlr{n - 3} u} \wlr{\sum_{k = 2}^n \wabs{ \wlr{\hat{s}_1 + z_2} + \sum_{i = 3}^k z_i}}
\]
and
\[
\wabs{\hat{s}_n - s_n} - \wabs{\hat{s}_1 - z_1} \leq
\frac{u}{1 - \wlr{n - 3} u} \wlr{
\wlr{n - 1} \wabs{\hat{s}_1 - z_1} + \sum_{k = 2}^n \wabs{\sum_{i = 1}^k  z_i}}.
\]
Since $\hat{s}_1 = \wflk{1}{z_1}$, Lemma \ref{lemUNear} leads to
\[
\wabs{\hat{s}_n - s_n} \leq
\frac{u}{1 + u} \wabs{z_1} +
\frac{u}{1 - \wlr{n - 3} u} \wlr{
\wlr{n - 1} \frac{u}{1 + u} \wabs{z_1} + \sum_{k = 2}^n \wabs{\sum_{i = 1}^k  z_i}}.
\]
\[
= \frac{u}{1 + u} \wlr{1 + \frac{\wlr{n - 1}u}{1 - \wlr{n - 3} u}} \wabs{z_1} +
\frac{u}{1 - \wlr{n - 3} u} \sum_{k = 2}^n \wabs{\sum_{i = 1}^k  z_i}.
\]
\[
\leq
\frac{u}{1 - \wlr{n - 2} u} \sum_{k = 1}^n \wabs{\sum_{i = 1}^k  z_i}
\]
\pbDef{lssB}
+ \wlr{\frac{1}{1 + u} \wlr{1 + \frac{\wlr{n - 1} u}{1 - \wlr{n - 3}u}} - \frac{1}{1 - \wlr{n - 2} u}} u \wabs{z_1}
\peDef{lssB}
The software Mathematica shows that
\[
 \frac{1}{1 + u} \wlr{ 1 + \frac{\wlr{n - 1} u}{1 - \wlr{n - 3} u}} - \frac{1}{1 - \wlr{n - 2} u}
= - \frac{\wlr{n -1} u^2}{\wlr{1 + u}\wlr{1 - \wlr{n - 2} u} \wlr{1 - \wlr{n - 3} u}},
\]
and this number is negative for $n \geq 2$ because $n u < 1$. As a result,
Equation \pRef{lssB} implies Equation \pRef{sspA} and we are done.
\peProof{Lemma}{lemSignedSum}\\

%%%%%%%%%%%%%%%%%%%%%%%%%%%%%%%%%%%%%%%%%%%%%%%%%%%%%%%%%%%%%%%%%%%%%%%%%
% End of proof of lemma lemSignedSum
%%%%%%%%%%%%%%%%%%%%%%%%%%%%%%%%%%%%%%%%%%%%%%%%%%%%%%%%%%%%%%%%%%%%%%%%%
%
%

%
%
%%%%%%%%%%%%%%%%%%%%%%%%%%%%%%%%%%%%%%%%%%%%%%%%%%%%%%%%%%%%%%%%%%%%%%%%%%%%%%%%%%%
% Proof of lemma lemUNearS
%%%%%%%%%%%%%%%%%%%%%%%%%%%%%%%%%%%%%%%%%%%%%%%%%%%%%%%%%%%%%%%%%%%%%%%%%%%%%%%%%%%

\pbProofB{Lemma}{lemUNearS} Let us start with $z > 0$
and define $m := \wlr{\wfloor{w} + \wceil{w}}/2$.
By Prop. \ref{propRoundNormal}, there are three possibilities :
\begin{itemize}
\item If $w < m$ then $r = \wfloor{w}$ satisfies Equation \pRef{rHat}.
\item If $w > m$ then $r = \wceil{w}$ satisfies Equation \pRef{rHat}.
\item If $w = m$ then $r_1 := \wfloor{w}$ and  $r_2 := \wceil{w}$
satisfy $r_i \in [0,\wlr{\beta - 1} \beta^{\wfpmu{}})$,
$\wabs{r_i - w} \leq 1/2$ and
$\wfl{z} = \beta^{e} \wlr{\beta^{\wfpmu{}} + r}$
for $r \in \wset{r_1,r_2}$. Therefore, Equation \pRef{rHat} is also
satisfied.
\end{itemize}
According to Definition \ref{longDefEps}, $2 u \times \beta^{\wfpmu{}} = 1$ and
Equation \pRef{rHat} yields
\pbDef{unsA}
\wabs{\frac{\wfl{z} - z}{z}} =
\frac{\wabs{r - w}}{\beta^{\wfpmu{}} + w} =
\frac{2 u \wabs{r -  w}}{1 + 2 w u} \leq
\frac{u}{1 + 2 w u}.
\peDef{unsA}
When $w \geq 1/2$, this equation implies that
\[
\wabs{\frac{\wfl{z} - z}{z}} \leq \frac{u}{1 + \max \wset{1,2 w} u},
\]
and when $w < 1/2$, Equation \pRef{rHat} and the fact that $r$ is integer
imply that $r = 0$ and
\[
\wabs{\frac{\wfl{z} - z}{z}} = \frac{w}{\beta^{\wfpmu{}} + w} = \frac{2 w u}{1 + 2 w u}
<
\frac{u}{1 + u} =
\frac{u}{1 + \max \wset{1, 2 w} u},
\]
and we have verified Equation \pRef{uNearSA}.
Equation \pRef{unsA} also leads  to
\[
\wabs{\frac{\wfl{z} - z}{z}} \leq
\frac{u}{1 + 2 w u} \leq
\frac{u}{1 + 2 \wlr{r - 1/2} u}
=
\frac{u}{1 + \wlr{2 r - 1} u}
\]
and
\[
\wabs{\frac{\wfl{z} - z}{\wfl{z}}} =
\frac{\wabs{r - w}}{\beta^{\wfpmu} + r} =
\frac{2 u \wabs{r -  w}}{1 + 2 r u}
\leq \frac{u}{1 + 2 r u}.
\]
This proves the last equation in Lemma \ref{lemUNearS} and we are done with $z > 0$.
To prove Lemma \ref{lemUNearS} for $z < 0$,
use the argument above for $z' = -z$ and the function $\wrm{m}$ in
Prop. \ref{propRoundMinus}.
\peProof{Lemma}{lemUNearS}\\

%%%%%%%%%%%%%%%%%%%%%%%%%%%%%%%%%%%%%%%%%%%%%%%%%%%%%%%%%%%%%%%%%%%%%%%%%%%%%%%%%%%
% End of proof of lemma lemUNearS
%%%%%%%%%%%%%%%%%%%%%%%%%%%%%%%%%%%%%%%%%%%%%%%%%%%%%%%%%%%%%%%%%%%%%%%%%%%%%%%%%%%
%
%

%
%
%%%%%%%%%%%%%%%%%%%%%%%%%%%%%%%%%%%%%%%%%%%%%%%%%%%%%%%%%%%%%%%%%%%%%%
% Proof of Lemma LemOpt
%%%%%%%%%%%%%%%%%%%%%%%%%%%%%%%%%%%%%%%%%%%%%%%%%%%%%%%%%%%%%%%%%%%%%%

\pbProofB{Lemma}{lemOpt}
Let us define
$\psi := \sup_{\wlr{\wvec{z},r} \in \wcal{Z} \times \wcal{R}} \wfc{g}{\wvec{z},r}$.
If $\psi \leq \varphi$ then $\wfc{g}{\wvec{z},r} \leq \varphi$ for
all $\wlr{\wvec{z},r} \in \wcal{Z} \times \wcal{R}$ and we
are done. Let us then assume that $\varphi < \psi$
and let
$\wset{\wlr{\wvec{z}_k, r_k}, k \in \wn{} } \subset \wcal{Z} \times  \wcal{R}$
be a sequence such that $\lim_{k \rightarrow \infty} \wfc{g}{\wvec{z}_k, r_k} = \psi$
and $\wfc{g}{\wvec{z}_k,r_k} > \varphi$. It follows that
$\wvec{z}_k \in \wcal{A}$ for each $k$ and
there exists $\lambda_k \in \wcal{L}$ and $r_k' \in \wcal{R}$
for which $\wvec{z}_k' := \lambda_k \wvec{z}_k \in \wcal{K}$
satisfies $\wfc{h}{\wvec{z}_k', r_k'} = \lambda_k \wfc{h}{\wvec{z}_k, r_k}$.
Since the sequence $\wvec{z}_k'$ is contained in the compact set $\wcal{K}$, it
has a subsequence which converges to $\wvec{z}^* \in \wcal{K}$,
and we may assume that this subsequence is $\wvec{z}_k'$ itself.
The scaling properties of $f$ lead to
\[
\wfc{f}{\wvec{z}_k',\wfc{h}{\wvec{z}_k',r_k'}}
= \wfc{f}{\lambda_k \wvec{z}_k,  \lambda_k  \wfc{h}{\wvec{z}_k,r_k}}
\geq
\wfc{f}{\wvec{z}_k,\wfc{h}{\wvec{z}_k,r_k}} = \wfc{g}{\wvec{x}_k, r_k}
\]
and
\[
\liminf_{k \rightarrow \infty} \wfc{f}{ \wvec{z}_{k}', \ \wfc{h}{\wvec{z}_{k}', r_{k}'}}
\geq \liminf_{k \rightarrow \infty } \wfc{g}{\wvec{z}_k, r_k} =
\lim_{k \rightarrow \infty } \wfc{g}{\wvec{z}_k, r_k} = \psi.
\]
Since $h$ is tight, there exists $r^* \in \wcal{R}$
and a subsequence $\wvec{z}_{n_k}'$
such that $\lim_{k \rightarrow \infty} \wfc{h}{\wvec{z}_{n_k}',r_{n_k}'} = \wfc{h}{\wvec{z}^*,r^*}$.
By the upper semi-continuity of $f$ and the maximality of $\psi$ we have
\[
\psi \geq \wfc{g}{\wvec{z}^*,r^*} =
 \wfc{f}{\wvec{z}^*, \wfc{h}{\wvec{z}^*, r^*}}
\]
\[
 \geq
 \limsup_{k \rightarrow \infty} \wfc{f}{ \wvec{z}_{n_k}', \ \wfc{h}{\wvec{z}_{n_k}', r_{n_k}'}}
 \geq
 \liminf_{k \rightarrow \infty} \wfc{f}{ \wvec{z}_{k}', \ \wfc{h}{\wvec{z}_{k}', r_{k}'}}
 \geq \psi.
\]
Therefore, $\wfc{g}{\wvec{z}^*, r^*} = \psi$ and we are done.
\peProof{Lemma}{lemOpt}\\

%%%%%%%%%%%%%%%%%%%%%%%%%%%%%%%%%%%%%%%%%%%%%%%%%%%%%%%%%%%%%%%%%%%%%%
% End of Proof of Lemma LemOpt
%%%%%%%%%%%%%%%%%%%%%%%%%%%%%%%%%%%%%%%%%%%%%%%%%%%%%%%%%%%%%%%%%%%%%%
%
%

\subsection{Corollaries}
\label{secCorol}
In this section we prove some of the corollaries stated in the article.
The remaining corollaries are proved in the extended version.\\

%
%
%%%%%%%%%%%%%%%%%%%%%%%%%%%%%%%%%%%%%%%%%%%%%%%%%%%%%%%%%%%%%%%%%%%%%%%%%%%%%%%%%%%
% Proof of Corollary CorSqrt
%%%%%%%%%%%%%%%%%%%%%%%%%%%%%%%%%%%%%%%%%%%%%%%%%%%%%%%%%%%%%%%%%%%%%%%%%%%%%%%%%%%

\pbProofB{Corollary}{corSqrt}
$\wfl{x^2} \geq \nu$ by Monotonicity (Prop. \ref{propRoundMonotone}),
and Lemma \ref{lemUNear} yield
\pbDef{sqrtA}
\wabs{z - \wabs{x}} \, \wlr{z + \wabs{x}} = \wabs{z^2 - x^2} =
\wabs{\wfl{x^2} - x^2} \leq \frac{\wabs{x}^2 u}{1 + u}
\peDef{sqrtA}
for $z := \sqrt{\wfl{x^2}} > 0$. It follows that $\delta := \wabs{z - \wabs{x}}/\wabs{x}$
satisfies
\[
\delta \leq \frac{u}{1 + u} \frac{\wabs{x}}{\wabs{x} + z} \leq \frac{u}{1 + u} < u = \beta^{-\mu}/2 \leq \frac{1}{4}
\ \ \ \Rightarrow \ \ 1 - \delta > 0.
\]
Equation \pRef{sqrtA} leads to
\[
\frac{u}{1 + u} \geq \delta
\frac{z + \wabs{x}}{\wabs{x}} \geq \delta \frac{2 \wabs{x} - \wabs{z - \wabs{x}}}{\wabs{x}} = \delta \wlr{2 - \delta} > 0,
\]
and
\[
1 - \delta =
\sqrt{\wlr{1 - \delta}^2} = \sqrt{1 - \delta \wlr{2 - \delta}} \geq \sqrt{1 - \frac{u}{1 + u}} = \frac{1}{\sqrt{1 + u}},
\]
and
\pbDef{sqrtPsi}
\delta \leq 1 - \frac{1}{\sqrt{1 + u}} = \frac{u}{2} \psi
\hspace{0.5cm} \wrm{for} \hspace{0.5cm}
\psi := \frac{2}{u} \frac{\sqrt{1 + u} - 1}{\sqrt{1 + u}} = \frac{2}{1 + u + \sqrt{1 + u}} < 1.
\peDef{sqrtPsi}

Let $\wfpsc$ be the complete system
with the same $\beta$ and $\mu$ as $\wfpf$. By Prop. \ref{propRoundAdapt}
there exists $\wflmx$ which rounds to nearest in $\wfpsc$ and is such that
$\wflx{w} = \wfl{w}$ for $w$ with $\wabs{w} \geq \nu_{\wfpf}$. In particular,
$\wfl{x^2} = \wflx{x^2}$.
Since $\nu < 1$ and $x^2 \geq \nu$ we have that
$\wabs{x} \geq \nu$ and by Prop. \ref{propNu}
there exists an exponent $e$ for $\wfpf$ and
$r \in [0,\wlr{\beta - 1} \beta^{\mu}) \cap \wz{}$ such that
$\wabs{x} = \beta^{e} \wlr{\beta^{\mu} + r}$.
This implies that $\beta^{e + \mu} \leq \wabs{x} < \beta^{e + \mu + 1}$.
The numbers $\beta^{2 e + 2 \mu}$ and $\beta^{2 e + 2 \mu + 2}$
are in $\wfpsc$  (although $\beta^{2 e + 2 \mu}$ may not be in $\wfpf$)
and, by the monotonicity of $\wflmx$,
\[
\beta^{2 e + 2 \mu} \leq \wfl{x^2} = \wflx{x^2} \leq \beta^{2 e + 2 \mu + 2},
\]
and $\beta^{e + \mu} \leq \sqrt{\wfl{x^2}} = z = \sqrt{\wfl{x^2}} \leq \beta^{e + \mu + 1}$.
By Prop. \ref{propOrder} and Prop. \ref{propNormalForm},
$z = \beta^{e} \wlr{\beta^{\mu} + w}$ with $0 \leq w \leq \wlr{\beta - 1} \beta^{\mu}$. As a result,
\[
\delta = \frac{\wabs{\beta^{e} \wlr{\beta^{\mu} + w} - \beta^{e} \wlr{\beta^{\mu} + r}}}
{\beta^{e} \wlr{\beta^{\mu} + r}} = \frac{\wabs{w - r}}{\beta^{\mu} + r},
\]
and recalling that $2 u \beta^{\mu} = 1$ and using Equation \pRef{sqrtPsi} we obtain
\pbDef{sqrtDR}
\wabs{w - r} \leq \frac{1}{4} \psi \wlr{1 + \beta^{-\mu} r} = \frac{1}{4} \psi \wlr{1 + 2 r u}.
\peDef{sqrtDR}
There are two possibilities: either
\pbDef{sqrtCRA}
\frac{1}{4} \psi \wlr{1 + 2 r u} < \frac{1}{2}
\peDef{sqrtCRA}
or
\pbDef{sqrtCRB}
\frac{1}{4} \psi \wlr{1 + 2 r u} \geq \frac{1}{2}.
\peDef{sqrtCRB}
In case \pRef{sqrtCRA} $\wabs{w - r} < 1/2$ by Equation \pRef{sqrtDR},
Prop. \ref{propRoundNormal} shows that $\wfl{z} = \wabs{x}$
and Corollary \ref{corSqrt} holds for $x$. For instance, if
$\beta = 2$ then $2 r u < 2 \wlr{2 - 1} 2^{\mu} u = 1$ and $r$ satisfies
Equation \pRef{sqrtCRA} because $\psi < 1$.
Therefore, we have proved Corollary \ref{corSqrt} for $\beta = 2$.

In order to complete the proof for the
cases in which Equation \pRef{sqrtCRB} is valid,
it suffices to show that
\pbDef{sqrtCond}
\frac{\wabs{x}}{\wfl{z}} = \frac{\wabs{x}}{\wfl{\sqrt{\wfl{x^2}}}} < 1 + u = \beta^{-\mu} \wlr{\beta^{\mu} + 1/2},
\peDef{sqrtCond}
because this equation implies that
$\wfl{\wabs{x}/\sqrt{\wfl{z}}} \leq 1$ by Prop. \ref{propRoundNormal} and monotonicity.

We first show that Equation \pRef{sqrtCond} is valid when
\pbDef{sqrtH}
\zeta := 1 + 2 r u > \wlr{1 + u}^{3/2} + 1 + u.
\peDef{sqrtH}
In fact, for $\psi$ in Equation \pRef{sqrtPsi}, Equation \pRef{sqrtH} is equivalent to
\[
\zeta > \frac{1 + u}{1 - \frac{\psi}{2} \wlr{1 + u}},
\hspace{0.7cm}
\frac{\zeta}{1 + u} - \frac{\psi}{2} \zeta - 1 > 0
\hspace{0.7cm} \wrm{and} \hspace{0.7cm}
\zeta - \frac{\psi}{2} \zeta u - u > \frac{\zeta}{1+u},
\]
and can also be written as
\pbDef{sqrtHB}
1 + u > \frac{\zeta}{\zeta - \frac{\psi}{2} \zeta u - u}
\hspace{0.7cm} \wrm{or} \hspace{0.7cm}
\frac{1 + 2 r u}{1 + 2 r u - \frac{\psi}{2}\wlr{1 + 2 r u} u  - u} < 1 + u.
\peDef{sqrtHB}
Since $w \in [0,\wlr{\beta - 1} \beta^{\mu}]$, Prop. \ref{propRoundNormal}
implies that $\wfl{z} \geq \beta^{e} \wlr{\beta^{\mu} + w - 1/2}$ and
\[
\frac{\wabs{x}}{\wfl{z}} \leq \frac{\beta^{e} \wlr{\beta^{\mu} + r}}{\beta^{e} \wlr{\beta^{\mu} + w - 1/2}}
= \frac{1 + 2 r u}{1 + 2  w u - u},
\]
because $2 u \beta^{\mu} = 1$. Equations \pRef{sqrtDR} shows that $w \geq r - \psi \wlr{1 + 2 r u} / 4$
and
\pbDef{sqrtHC}
\frac{\wabs{x}}{\wfl{z}} \leq \frac{1 + 2 r u}{1 + 2 r u - \frac{\psi}{2} \wlr{1 + 2 r u} u  - u}.
\peDef{sqrtHC}
Equations \pRef{sqrtHB}  and \pRef{sqrtHC} lead to Equation \pRef{sqrtCond}.
Therefore, Equation \pRef{sqrtH} implies Equation \pRef{sqrtCond} and
Corollary \ref{corSqrt} is valid when Equation \pRef{sqrtH} is satisfied.

In the case opposite to Equation \pRef{sqrtH} we have that
\pbDef{sqrtTRPA}
2 r u \leq \wlr{1 + u}^{3/2} + u = 1 + \frac{5}{2} u + \frac{3}{8 \sqrt{1 + \xi_1}} u^2
\peDef{sqrtTRPA}
for some $\xi_1 \in [0,u]$. Since $r$ is integer and $2 u = \beta^{-\mu}$, Equation \pRef{sqrtTRPA} implies that
\pbDef{sqrtTRA}
r < \beta^{\mu} + \frac{5}{4} + \frac{3}{16} u < \beta^{\mu} + 2 \Rightarrow r \leq \beta^{\mu} + 1.
\peDef{sqrtTRA}
Moreover, Equation \pRef{sqrtCRB} leads to
\[
r \geq \beta^{\mu} \frac{2 - \psi}{\psi} = \beta^{\mu} \wlr{u + \sqrt{1 + u}} =
\beta^{\mu} \wlr{1 + \frac{3}{2} u - \frac{1}{8 \wlr{1 + \xi_2}^{3/2}} u^2}
\]
for some $\xi_2 \in [0,u]$, and since $r$ is integer and $2 u = \beta^{-\mu}$,
we have that
\pbDef{sqrtTR}
r \geq \beta^{\mu} + \frac{3}{4}  - \frac{1}{16 \wlr{1 + \xi_2}^{3/2}} u
\Rightarrow r \geq \beta^{\mu} + 1.
\peDef{sqrtTR}
Equations \pRef{sqrtTRA} and \pRef{sqrtTR} show that
there is just one $r$ left: $r = \beta^{\mu} + 1$, which corresponds
to $\wabs{x} = \beta^{e} \wlr{2 \beta^{\mu} + 1}$. It follows that
\[
x^2 = \beta^{2 e}\wlr{4 \beta^{2 \mu} + 4 \beta^{\mu} + 1} =
\beta^{2 e + \mu}\wlr{\beta^{\mu} + \wlr{3 \beta^{\mu} + 4 + \beta^{-\mu}}}.
\]
If $\beta \geq 5$ then $3 \beta^{\mu} + 4 + \beta^{-\mu} < \wlr{\beta - 1} \beta^{\mu}$
and Prop. \ref{propRoundNormal} implies that
\[
 \wfl{x^2} = 4 \beta^{2 e + \mu}\wlr{\beta^{\mu} + 1}
 \Rightarrow z = \sqrt{\wfl{x^2}} = 2 \beta^{e + \mu} \sqrt{1 + \beta^{-\mu}}
  \]
 \[
 = 2 \beta^{e + \mu} \wlr{1 + \frac{1}{2} \beta^{-\mu} - \frac{\theta_5}{2} \beta^{-\mu}},
\]
where, for  some $\xi_5 \in [0,\beta^{-\mu}]$,
\[
0 \leq \theta_5 := \frac{1}{4 \wlr{1 + \xi_5}^{3/2}} \beta^{-\mu} \leq \frac{1}{4} \times \frac{1}{5} = \frac{1}{20}.
\]
Therefore,
$z := \sqrt{\wfl{x^2}} = \beta^{e} \wlr{2 \beta^{\mu} + 1 - \theta_5}$
and the bound $\wabs{\theta_5} \leq 1/20$ and Prop. \ref{propRoundNormal} imply that
$\wfl{z} = \beta^{e} \wlr{2 \beta^{\mu} + 1} = \wabs{x}$
and we are done with the case $\beta \geq 5$.

For $\beta = 3$, the critical $x$ is $3^{e} \wlr{2 \times 3^{\mu} + 1}$ and
\[
x^2 = 3^{2 e}\wlr{4 \times 3^{2 \mu} + 4 \times 3^{\mu} + 1}
=
3^{2 e + \mu + 1}\wlr{3^{\mu} + 3^{\mu - 1} + 1 + \wlr{\frac{1}{3} + 3^{-\mu-1}}}
\]
The bound
\[
\frac{1}{3} + 3^{-\mu-1} \leq \frac{1}{3} + \frac{1}{9} = \frac{4}{9} < 1/2
\]
and Prop. \ref{propRoundNormal} lead to
\[
\wfl{x^2} = 3^{2 e + \mu + 1}\wlr{3^{\mu} + 3^{\mu - 1} + 1}
= 4 \times 3^{2 e + 2 \mu}\wlr{1 + \frac{3}{4} \times 3^{-\mu}}
\]
and
\[
z := \sqrt{\wfl{x^2}} = 2 \times 3^{e + \mu}\wlr{1 + \frac{3}{8} \times 3^{-\mu} - \frac{\theta_3}{2} \times 3^{-\mu}}
= 3^{e}\wlr{2 \times 3^{\mu} + \frac{3}{4} - \theta_3}
\]
where, for some $\xi_3 \in [0,1/3]$,
\[
0 \leq
\theta_3 := \frac{1}{4 \wlr{1 + \xi_3}^{3/2}} \times \frac{9}{16} \times 3^{-\mu} \leq \frac{3}{64}.
\]
Since $3/4 - 3 / 64 = 45 / 64 > 1/2$, Prop. \ref{propRoundNormal} shows that
$\wfl{z} = \wabs{x}$ when $\beta = 3$.

Finally, for $\beta = 4$, we care about $x = 4^{e} \wlr{2 \times 4^{\mu} + 1}$ and
\[
x^2 =
4^{2 e}\wlr{4 \times 4^{2 \mu} + 4 \times 4^{\mu} + 1}
=
4^{2 e + 1 + \mu}\wlr{4^{\mu} + 1 + 4^{-\mu - 1}},
\]
$4^{-\mu - 1} < 1/2$ and Prop. \ref{propRoundNormal} yields
\[
\wfl{x^2} = 4^{2 e + 1 + \mu}\wlr{4^{\mu} + 1} =
4^{2 e + 1 + 2 \mu}\wlr{1 + 4^{-\mu}}.
\]
It follows that
\[
z := \sqrt{\wfl{x^2}} = 2 \times 4^{e + \mu}\sqrt{1 + 4^{-\mu}} =
 2 \times 4^{e + \mu} \wlr{1 + \frac{1}{2} \times 4^{-\mu} - \frac{\theta_4}{2} \times 4^{-\mu}}
\]
where, for some $\xi_4 \in [0,1/4]$,
\[
0 < \theta_4 := \frac{1}{4 \sqrt{1 + \xi_4}} \, 4^{-\mu} < \frac{1}{16}.
\]
Therefore, $z = 4^{e + 1} \wlr{2 \times 4^{\mu} + 1 - \theta_4}$,
$\wfl{z} = \wabs{x}$ and we are done.
\peProof{Corollary}{corSqrt} \\

%%%%%%%%%%%%%%%%%%%%%%%%%%%%%%%%%%%%%%%%%%%%%%%%%%%%%%%%%%%%%%%%%%%%%%%%%
% End of Proof of Corollary corSqrt
%%%%%%%%%%%%%%%%%%%%%%%%%%%%%%%%%%%%%%%%%%%%%%%%%%%%%%%%%%%%%%%%%%%%%%%%%
%
%

%
%
%%%%%%%%%%%%%%%%%%%%%%%%%%%%%%%%%%%%%%%%%%%%%%%%%%%%%%%%%%%%%%%%%%%%%%%%%
% Proof of corollary SharpSumIEEE
%%%%%%%%%%%%%%%%%%%%%%%%%%%%%%%%%%%%%%%%%%%%%%%%%%%%%%%%%%%%%%%%%%%%%%%%%

\pbProofB{Corollary}{corNormOneUnperfect}
Let $\wfpsc$ be the perfect system corresponding to $\beta$ and $\mu$
and $\wflmtx$ the rounding tuple
in Prop. \ref{propRoundIEEEX} or \ref{propRoundMPFRX}, depending
on whether $\wfpf$ is an IEEE system or a MPFR system.
As in the proof of Lemma \ref{lemNormOneBound}, we define
$z_1 := y_0 + y_1$, $z_k := y_k$ for $2 \leq k \leq n$,
$s_k := \sum_{i = 1}^k z_i$ and $\hat{s}_k := \wfpsumkf{k}{\wvec{x},\wflmt}$
for $k = 0,\dots,n$. We also
use the set $\wcal{T}$ of indexes $k$ in $[1,n]$ such that
$\wabs{\wfpsumkf{k-1}{\wvec{z},\wflmt} + z_k} < \tau$
for
\[
\tau := \beta^{\wfpemin} \wlr{\beta^{\mu} + r}
\hspace{1cm} \wrm{and} \hspace{1cm}
 r := \beta^{\mu} \frac{\beta - 1}{2}.
\]
Note that $\tau \in \wfpbin_{\wfpemin} \subset \wfpf$ because
$r$ is integer and $r < \wlr{\beta - 1} \beta^{\mu}$.
The threshold $\tau$ was chosen
because $\wfpnu = \beta^{\wfpemin + \wfpmu}$,
\pbDef{nouTau}
\tau =  \frac{\beta + 1}{2} \wfpnu < \beta \nu
\peDef{nouTau}
and Prop. \ref{propRoundBelowAlpha} shows that
\pbDef{sseeeA}
\wabs{z} \leq \beta \wfpnu \Rightarrow \wabs{\wfl{z} - z} \leq \wfpa / 2,
\peDef{sseeeA}
where $\wfpa = \beta^{\wfpemin}$ for IEEE systems and
$\wfpa = \wfpnu = \beta^{\wfpemin + \mu}$ for MPFR systems.

Let $m \in [0,n]$ be the size of $\wcal{T}$.
We prove by induction that
\[
\wfc{\eta}{\wvec{z},\wflmt}
:= \sum_{k=1}^n \wabs{\wfpsumkf{k}{\wvec{z},\wflmt} - \wlr{\wfpsumkf{k-1}{\wvec{z},\wflmt} + z_k}}
\]
satisfies
\pbTClaim{thSharpSumUnperfectPV}
\wfc{\eta}{\wvec{z},\wflmt} \leq
\frac{m \alpha}{2} +
\frac{\wlr{n - m} u}{1 + \wlr{n - m} u} \wlr{\frac{m \alpha}{2} + \sum_{k = 1}^n \wabs{z_k}}.
\peTClaim{thSharpSumUnperfectPV}
If $m = 0$ then $\wfpsumkf{n}{\wvec{z},\wflmt} = \wfpsumkf{n}{\wvec{z},\wflmtx}$
and Equation \pRef{thSharpSumUnperfectPV} follows from Lemma \ref{lemNormOneBound}.
Assuming that Equation \pRef{thSharpSumUnperfectPV} holds  for $m - 1$,
let us show that it holds for $m$. If $\wabs{s_1} < \tau$ then
the sum $\wlr{\hat{s}_1 + z_2} + \sum_{k = 3}^n z_k$ has $n - 1$
parcels and there are $m - 1$ indices in $[2,n] \cap \wcal{T}$.
As a result $(n-1) - (m-1) = n - m$, Equation \pRef{sseeeA}, the identity $s_1 = z_1$ and induction yield
\[
\wfc{\eta}{\wvec{z},\wflmt}  =
\wabs{\hat{s}_1 - s_1} +
\wlr{\sum_{k = 2}^n \wabs{\hat{s}_k - \wlr{\hat{s}_{k-1} + z_k}}}
\]
\[
\leq
\frac{\alpha}{2} +  \wlr{\frac{\wlr{m - 1}\alpha}{2} +
\frac{\wlr{n - m} u}{1 + \wlr{n - m} u} \wlr{
 \frac{\wlr{m - 1}\alpha}{2} +
 \wabs{\hat{s}_1 + z_2} + \sum_{k = 3}^n \wabs{z_k}}}
\]
\[
\leq
\frac{m \alpha}{2} +
\frac{\wlr{n - m} u}{1 + \wlr{n - m} u} \wlr{
\wlr{ \wabs{\hat{s}_1 - s_1} - \frac{\alpha}{2}} +
\frac{m \wfpa}{2} +
\wabs{s_1} + \sum_{k = 2}^n \wabs{z_k}}
\]
\[
\leq
\frac{m \alpha}{2} +
\frac{\wlr{n - m} u}{1 + \wlr{n - m} u} \wlr{ \frac{m \wfpa}{2} + \sum_{k = 1}^n \wabs{z_k}}.
\]
Therefore, Equation \pRef{thSharpSumUnperfectPV} holds when $\wabs{s_1} < \tau$.
Let us then assume that $\wabs{s_1} \geq \tau$ and define
$\ell \in [2,n]$ as the first index such that
$\wabs{\hat{s}_{\ell - 1} + z_{\ell}} < \tau$,
\pbDef{IEEESpq}
S := \sum_{k = 1}^{\ell - 1} \wabs{z_k},
\hspace{1cm}
p := \ell - 1
\hspace{1cm} \wrm{and} \hspace{1cm}
q := n - m - \ell + 1.
\peDef{IEEESpq}
Monotonicity and $\tau \in \wfpf$ implies that
$\wabs{\hat{s}_\ell} = \wabs{\wflk{\ell}{\hat{s}_{\ell - 1} + z_{\ell}}} \leq \tau$
and the proof of Lemma \ref{lemNormOneBound},  Equation \pRef{sseeeA} and induction yield
\[
\wfc{\eta}{\wvec{z},\wflmt} =
\sum_{k = 1}^{\ell - 1} \wabs{\hat{s}_{k} - \wlr{s_{k - 1} + z_k}} +
\wabs{\hat{s}_{\ell} - \hat{s}_{\ell - 1} - z_{\ell}} +
\sum_{k = \ell + 1}^{n} \wabs{\hat{s}_{k} - \wlr{\hat{s}_{k-1} + z_k}}
\]
\[
\leq
\frac{p u}{1 + p u} S +
\frac{\wfpa}{2} +
\wlr{  \frac{\wlr{m - 1}\alpha}{2} + \frac{q u}{1 + q u} \wlr{
\frac{ \wlr{m - 1} \alpha }{2} +
\wabs{\hat{s}_{\ell} + z_{\ell + 1}} + \sum_{k = \ell + 2}^n \wabs{z_k}}}
\]
\pbDef{ssIEEEA}
\leq
\frac{p u}{1 + p u} S +
\frac{m \alpha}{2} +
\frac{q u}{1 + q u} \wlr{\frac{\wlr{m - 1}\alpha}{2} + \tau + \sum_{k = \ell + 1}^n \wabs{z_k}}.
\peDef{ssIEEEA}
If $S \geq 7 \tau / 6$ then
\[
\frac{p}{1 + p u} S + \frac{q}{1 + q u} \tau
\leq \wlr{\frac{p}{1 + p u} + \frac{6}{7}  \frac{q}{1 + q u}} S
\leq \frac{\wlr{p + q} u}{1 + \wlr{p + q} u} S - \Delta S
\]
for
\[
\Delta := \frac{p + q}{1 + \wlr{p + q} u} - \wlr{\frac{p}{1 + p u} + \frac{6}{7} \frac{q}{1 + q u}}.
\]
The software Mathematica shows that
\[
\Delta = q \frac{1 + q u - 6 \wlr{2 + q u + p u} p u }{\wlr{1 + p u} \wlr{1 + q u} \wlr{1 + \wlr{p + q} u}}
\]
and the hypothesis $20 n u \leq 1$ implies that $\Delta \geq 0$.
Therefore, if $S \geq 7 \tau / 6$ then Equation \pRef{ssIEEEA} leads to
\[
\wfc{\eta}{\wvec{z},\wflmt} \leq  \frac{m \alpha}{2}  +
\frac{\wlr{p + q} u}{1 + \wlr{p + q} u} \wlr{\frac{m \alpha}{2} + S +\sum_{k = \ell + 1}^n \wabs{z_k}},
\]
and Equation \pRef{thSharpSumUnperfectPV} follows from Equation \pRef{IEEESpq}. We can
 then assume that $S < 7 \tau / 6$ and, for $1 \leq k < \ell$,
Lemma \ref{lemNormOneBound} leads to
\[
\wabs{\hat{s}_k} \leq \wabs{s_k} + \wabs{\hat{s}_k - s_k}  \leq \wlr{1 + \frac{k u}{1 + k u}} S <
1.05 \times \frac{7}{6} \frac{\beta + 1}{2} \nu
< \frac{2}{3} \wlr{\beta + 1} \nu \leq \beta \wfpnu,
\]
and Equation \pRef{sseeeA} implies that $\wabs{\hat{s}_{k} - \wlr{\hat{s}_{k} + z_k}} \leq \wfpa/2$
for $1 \leq k < \ell$. It follows that
\pbDef{ssIEEEBA}
\sum_{k=1}^{\ell - 1} \wabs{\hat{s}_{k} - \wlr{s_{k - 1} - z_k}}
\leq \wlr{\ell - 1} \wfpa / 2 = p \wfpa / 2.
\peDef{ssIEEEBA}
The identity $u \nu = \wfpa/2$ for IEEE systems
and the inequality $u \nu = u \wfpa \leq \wfpa/4$ for MPFR systems,
the hypothesis $20 n u \leq 1$ and
the fact that
\[
\sum_{k = 1}^{\ell - 1} \wabs{z_k} \geq \wabs{z_1} = \wabs{s_1} \geq \tau = \wlr{\beta + 1} \nu / 2 \geq \frac{3}{2} \nu
\]
imply that
\[
\frac{p \wfpa}{2} \leq \frac{p \wfpnu u}{2} \leq
\frac{2}{3} \wlr{1 + p u} \frac{p u}{1 + p u} \sum_{k = 1}^{\ell -1} \wabs{z_k}
\]
\[
\leq \frac{2}{3} \times \frac{21}{20} \times \frac{p u}{1 + p u} \sum_{k = 1}^{\ell - 1} \wabs{z_k} =
\frac{7}{10} \frac{p u}{1 + p u} \sum_{k = 1}^{\ell - 1} \wabs{z_k}.
\]
Using induction as in Equation \pRef{ssIEEEA}
and the bounds in the previous equation and in Equation \pRef{ssIEEEBA},
and recalling that $\wabs{z_1} \geq \tau$,  we obtain
\[
\wfc{\eta}{\wvec{z},\wflmt} \leq \frac{p \wfpa}{2} + \frac{\alpha}{2} +
\frac{\wlr{m - 1} \alpha}{2} + \frac{qu}{1 + q u} \wlr{\frac{\wlr{m - 1} \wfpa}{2} + \tau + \sum_{k = \ell + 1}^n \wabs{z_k}}
\]
\pbDef{IEEEUfa}
\leq
\frac{7}{10} \frac{p u}{1 + p u}\sum_{k = 1}^n \wabs{z_k} +
\frac{m \alpha}{2} + \frac{q u}{1 + q u} \wlr{\frac{m \wfpa}{2}  + \sum_{k = 1}^n \wabs{z_k}}.
\peDef{IEEEUfa}
According to the software Mathematica,
\[
\frac{p + q}{1 + \wlr{p + q} u} - \wlr{\frac{q}{1 + q u} + \frac{7}{10} \frac{p}{1 + p u}}
= p \frac{3 + 3 p u - 7 \wlr{2 + q u + p u} q u}{10 \wlr{1 + p u} \wlr{1 + q u} \wlr{1 + \wlr{p + q} u}},
\]
and this number is positive due to the hypothesis $20 n u \leq 1$. As a result,
Equation \pRef{IEEEUfa} implies Equation \pRef{thSharpSumUnperfectPV}
and this concludes the inductive proof of Equation \pRef{thSharpSumUnperfectPV}.
This equation leads to
\pbDef{nouF}
\wabs{
\, \wflt{\sum_{k = 0}^n y_k} - \sum_{k = 0}^n y_k \,} \leq
\frac{m \alpha}{2} + \frac{\wlr{n - m} u}{1 + \wlr{n - m} u} \wlr{\frac{m \wfpa}{2} + \sum_{k = 0}^n \wabs{y_k}},
\peDef{nouF}
and implies Equation \pRef{thSharpSumUnperfect} because $0 \leq m \leq n$.

Finally, when the additional condition in Corollary \ref{corNormOneUnperfect} holds
we have that
\[
\sum_{k = 0}^n \wabs{y_k} \geq \theta \frac{\wlr{1 + n u}^2}{u} \wfpa
\]
for $1 > \theta := \frac{1}{\wlr{1 + n u}^2} \geq \wlr{20/21}^2 > 0.9$ and
the software Mathematica shows that
\[
\frac{n u}{1 + n u} \theta \frac{\wlr{1 + n u}^2}{u} \wfpa
- \wlr{\frac{m \wfpa}{2} + \frac{\wlr{ n - m} u}{1 + \wlr{n - mu} u}
\wlr{\frac{m \wfpa}{2} + \theta \frac{\wlr{1 + n u}^2}{u} \wfpa}}
\]
\[
= \wfpa m \frac{2 \theta - 1 + 2 u \wlr{ \wlr{\theta - 1} n + m}}{1 + \wlr{n - m} u}
\geq \wfpa m \frac{0.8 - 2 u \wlr{0.1 n - m}}{1 + \wlr{n - m} u}
> 0,
\]
and Equation \pRef{normOneBound} follows from Equation \pRef{nouF}.
\peProof{Corollary}{corNormOneUnperfect}\\

%%%%%%%%%%%%%%%%%%%%%%%%%%%%%%%%%%%%%%%%%%%%%%%%%%%%%%%%%%%%%%%%%%%%%%%%%
% End of proof of corollary SharpSumIEEE
%%%%%%%%%%%%%%%%%%%%%%%%%%%%%%%%%%%%%%%%%%%%%%%%%%%%%%%%%%%%%%%%%%%%%%%%%
%
%

%
%
%%%%%%%%%%%%%%%%%%%%%%%%%%%%%%%%%%%%%%%%%%%%%%%%%%%%%%%%%%%%%%%%%%%%%%%%%
% Proof of Corollary corSignedSum
%%%%%%%%%%%%%%%%%%%%%%%%%%%%%%%%%%%%%%%%%%%%%%%%%%%%%%%%%%%%%%%%%%%%%%%%%

\pbProofB{Corollary}{corSignedSum}
Define $z_1 := y_0 + y_1$ and $z_k := y_k$ for $2 \leq k \leq n$,
$s_k := \sum_{i = 1}^k z_i$  and
$\hat{s}_{k} := \wfpsumkf{k}{\wvec{x},\wflmt}$ for $k = 0, \dots, n$.
Let $\wfpsc$ be the perfect system corresponding to
$\beta$ and $\mu$ and $\wflmtx$
the rounding tuple in Props. \ref{propRoundIEEEX} or \ref{propRoundMPFRX},
depending on whether $\wfpf$ is an IEEE system or a MPFR system.
By definition of $\wflmtx$, we have that
$\wflk{k}{s_{k-1} + z_k} = \wflxkf{k}{s_{k-1} + z_k}$
when $\wabs{s_{k-1} + z_k} \geq \wfpnu$.
Let $\wcal{T}$ be the set
of indexes $k$ in $[1,n]$ such that
$\wabs{\wfpsumkf{k-1}{\wvec{z},\wflmt} + z_k} < \wfpnu$
and $m \in [0,n]$ its size.
We prove by induction that
\pbTClaim{csi}
\wabs{
\, \wflt{\sum_{k = 1}^n z_k} - \sum_{k = 1}^n z_k
\,
} \leq
\wlr{1 + 2 \wlr{n - m} u} m \frac{\alpha}{2} +
\frac{1 - m u / 2}{1 - \wlr{n - 2} u} u \sum_{k = 1}^n \wabs{\sum_{i = 1}^k z_i}.
\peTClaim{csi}
When $m = 0$ we have that $\hat{s}_n = \wfpsumkf{n}{\wvec{z},\wflmtx}$ and
Equation \pRef{csi} follows from Lemma \ref{lemSignedSum}.
Assuming that Equation \pRef{csi} holds for $m - 1$, let us
prove it for $m$. Let $\ell$ be the last element of $\wcal{T}$
(Note that $\ell \geq m$.)
It follows that $\wabs{\hat{s}_{k - 1} - z_k} \geq \nu$
for $k > \ell$ and
$\hat{s}_k =  \wflxkf{k}{\hat{s}_{\ell} + \sum_{i = \ell + 1}^k z_i}$ for $k > \ell$.
The proof of Lemma \ref{lemSignedSum} shows that
\[
\wabs{\hat{s}_{n} - \wlr{\hat{s}_{\ell} + \sum_{k = \ell + 1}^n z_k}} \leq
\frac{u}{1 - \wlr{\wlr{n - \ell} - 2} u}
\sum_{k = \ell + 1}^n \wabs{\hat{s}_{\ell} + \sum_{i = \ell + 1}^k z_i}
\]
\[
\leq \frac{u}{1 - \wlr{n - \ell - 2} u}
\sum_{k = \ell + 1}^n  \wlr{\wabs{\hat{s}_{\ell} - s_{\ell}} + \wabs{\sum_{i = 1}^k z_i}}
\]
\pbTClaim{csib}
= A u \wlr{\wlr{n - \ell} \wabs{\hat{s}_{\ell} - s_{\ell}} + \sum_{k= \ell + 1}^n \wabs{\sum_{i = 1}^k z_i}},
\peTClaim{csib}
for
\[
A := \frac{1}{1 - \wlr{n - \ell - 2} u}.
\]
Moreover, $\wabs{\hat{s}_{\ell - 1} + z_\ell} < \wfpnu$ and, by induction
and Prop. \ref{propRoundBelowAlpha},
\[
\wabs{\hat{s}_{\ell} - s_{\ell}} \leq
\wabs{\hat{s}_{\ell} - \hat{s}_{\ell - 1} - z_{\ell}} +
\wabs{\hat{s}_{\ell - 1} - s_{\ell - 1}}
\]
\[
\leq
\frac{\alpha}{2} + \wlr{1 + 2  \wlr{\ell - m} u} \wlr{m - 1 } \frac{\alpha}{2} +
\frac{\wlr{1 - \wlr{m - 1} u/2} u}{1 - \wlr{\ell - 3} u} \sum_{k = 1}^{\ell - 1} \wabs{\sum_{i = 1}^k z_i}
\]
\pbTClaim{cisc}
= \wlr{m + 2 \wlr{\ell - m} \wlr{m - 1} u} \frac{\alpha}{2} +
C u \sum_{k = 1}^{\ell - 1} \wabs{\sum_{j = 1}^k z_j}.
\peTClaim{cisc}
for
\[
C := \frac{1 - \wlr{m - 1} u / 2}{1 - \wlr{\ell - 3} u}.
\]
Combining Equations \pRef{csib} and \pRef{cisc} we obtain
\[
\wabs{\hat{s}_n - s_n} \leq
\wabs{\hat{s}_n - \wlr{\hat{s}_{\ell} + \sum_{k=\ell + 1}^n z_k}}
+ \wabs{\hat{s}_\ell - s_{\ell}} \leq
\]
\[
\leq
\wlr{1 + A \wlr{n - \ell} u} \wabs{\hat{s}_\ell - s_{\ell}} + A u \sum_{k=\ell + 1}^n \wabs{\sum_{i =1}^k z_i}
\]
\pbTClaim{cisd}
\leq
 D \wlr{m + 2 \wlr{\ell - m} \wlr{m - 1} u} \frac{\alpha}{2}
+ D C u \sum_{k = 1}^{\ell - 1} \wabs{\sum_{i =1}^k z_i}
+ A u \sum_{k=\ell + 1}^n \wabs{\sum_{i =1}^k z_i},
\peTClaim{cisd}
for
\[
D := 1 +  A \wlr{n - \ell} u = \frac{1 + 2 u}{1 - \wlr{n - \ell - 2} u}.
\]
We now show that $Q < 1$ for
\[
 Q := \frac{D \wlr{m + 2 \wlr{\ell - m} \wlr{m - 1} u}}
 {
 \wlr{1 + 2 \wlr{n - m} u} m
 }
=
\frac{\wlr{1 + 2 u}\wlr{1 + 2 \wlr{\ell - m}\wlr{1 - 1/m}u}}
{\wlr{1 - \wlr{n - \ell - 2} u}{\wlr{1 + 2 \wlr{n - m} u}}} .
\]
It easy to see that $Q < 1$ when $\ell = n$. Since $20 n u \leq 1$,
when $\ell < n$ we have
\[
 Q < \frac{\wlr{1 + 2 u}\wlr{1 + 2 \wlr{\ell - m} u} }
{\wlr{1 - \wlr{n - \ell - 2} u}{\wlr{1 + 2 \wlr{n - m} u}}}
\]
\[
= \frac{1 + \wlr{2 \ell - 2 m + 2} u + 4 \wlr{\ell - m} u^2}
{
1 + \wlr{n + \ell - 2 m + 2} u - 2 \wlr{n - \ell - 2} \wlr{n - m} u^2
}
\]
\[
=
\frac{1 + \wlr{2 \ell - 2 m + 2} u + 4 \wlr{\ell - m} u^2}
{1 + \wlr{2 \ell - 2 m + 2} u + \wlr{n - \ell} \wlr{1 - 2
\frac{\wlr{n - \ell - 2}}{n - \ell} \wlr{n - m} u} u}
\]
\[
\leq  \frac{1 + \wlr{2 \ell - 2 m + 2} u + 0.2 u}
{
1 + \wlr{2 \ell - 2 m + 2} u + \wlr{1 - 0.1} u
} < 1.
\]
Therefore, $Q < 1$ and, equivalently,
\pbDef{csad}
D \wlr{m + 2 \wlr{\ell - m} \wlr{m - 1} u}  \leq \wlr{1 + 2 \wlr{n - m} u} m.
\peDef{csad}
Moreover,
\[
D C = \frac{1 + 2 u}{1 - \wlr{n - \ell - 2} u} \frac{1 - \wlr{m - 1} u / 2}{1 - \wlr{\ell - 3} u}
=  \frac{\wlr{1 + 2 u} \wlr{1 - \wlr{m - 1} u / 2}}{1 - \wlr{n - 5} u + \wlr{\ell - 3}\wlr{n - \ell - 2} u^2 }.
\]
Note that the function $\wf{h}{\ell} := \wlr{\ell - 3}\wlr{n - \ell - 2}$ is concave.
Therefore its minimum in the interval $[1,n]$ is at the endpoints. Since
$\wfc{h}{1} = \wf{h}{n} = - 2 \wlr{n - 3}$, we have
\[
D C \leq \frac{\wlr{1 + 2 u} \wlr{1 - \wlr{m - 1} u / 2}}{1 - \wlr{n - 5} u - 2  \wlr{n - 3} u^2 },
\]
and the software Mathematica shows that
\[
\frac{\wlr{1 + 2 u} \wlr{1 - \wlr{m - 1} u / 2}}{1 - \wlr{n - 5} u - 2  \wlr{n - 3} u^2 }
- \frac{1 - m u / 2}{1 - \wlr{n - 2} u }
= - u \frac{1 - 2 u - m u + n u}{2 \wlr{1 + 3 u - n u}\wlr{1 + 2 u - n u}} < 0,
\]
where the last inequality follows from the hypothesis $20 n u \leq 1$.
Therefore,
\[
D C \leq \frac{1 - m u / 2}{1 - \wlr{n - 2} u }.
\]
Not also that, since $\ell \geq m$  and $20 n u \leq 1$,
\[
A - \frac{1 - m u/2}{1 - \wlr{n - 2} u} =
\frac{1 - \wlr{n - 2} u - \wlr{1 - m u /2} \wlr{1 - \wlr{ n - \ell - 2} u}}{\wlr{1 - \wlr{n - 2} u} \wlr{1 - \wlr{n - \ell - 2} u}}
\]
\[
= - \frac{\ell - m / 2 \wlr{1 - \wlr{ n - \ell - 2} u}}{\wlr{1 - \wlr{n - 2} u} \wlr{1 - \wlr{n - \ell - 2} u}} u < 0,
\]
and
\[
A \leq \frac{1 - m u / 2}{1 - \wlr{n - 2} u }.
\]
The bounds on $DC$ and $A$ above, combined with Equations
\pRef{cisd} and \pRef{csad} imply Equation \pRef{csi}, and we
completed the inductive proof of this equation.

Finally, when $u \sum_{k = 1}^n \wabs{ \sum_{i = 0}^n y_i} \geq n \wfpa$
Equation \pRef{csi} leads to
\[
\wabs{
\, \wflt{\sum_{k = 0}^n y_k} - \sum_{k = 0}^n y_k
\,
} \leq  \theta_{m} u \sum_{k = 1}^n \wabs{\sum_{i = 0}^k y_i},
\]
for
\[
\theta_{m} := \wlr{1 + 2 \wlr{n - m} u} \frac{m}{2n } + \frac{1 - m u / 2}{1 - \wlr{n - 2} u}.
\]
The derivative of $\theta_m$ with respect to $m$ is
\[
\frac{1 - u \wlr{4 m - 2} - 2 u^2 \wlr{n^2 - 2 m n + 4 m - 2 n}}{2 n \wlr{1 + 2 u - n u}},
\]
and it is positive because $20 m u \leq 20 n u \leq 1$. Thus, $\theta_m$ is maximized for
$m = n$ and
\[
\theta_m \leq \frac{m}{2} + \frac{1 - n u / 2}{1 - \wlr{n - 2} u} =
\frac{3 - 2 \wlr{n - 1} u}{2 \wlr{1- \wlr{n -2} u}},
\]
and Equation \pRef{thCorSignedSumL} holds because
\[
\frac{3 - 2 \wlr{n - 1} u}{2 \wlr{1- \wlr{n -2} u}} - \frac{3}{2} \wlr{1 + \frac{n u}{2}}
= - u \frac{8 + n \wlr{1  - 3 \wlr{n - 2} u}}{1 - \wlr{n - 2} u} < 0
\]
when $20 n u \leq 1$.
\peProof{Corollary}{corSignedSum}\\

%%%%%%%%%%%%%%%%%%%%%%%%%%%%%%%%%%%%%%%%%%%%%%%%%%%%%%%%%%%%%%%%%%%%%%%%%
% End of Proof of Corollary corSignedSum
%%%%%%%%%%%%%%%%%%%%%%%%%%%%%%%%%%%%%%%%%%%%%%%%%%%%%%%%%%%%%%%%%%%%%%%%%
%
%

%%%%%%%%%%%%%%%%%%%%%%%%%%%%%%%%%%%%%%%%%%%%%%%%%%%%%%%%%%%%%%%
% The auxiliary results are proved only in the full version
%%%%%%%%%%%%%%%%%%%%%%%%%%%%%%%%%%%%%%%%%%%%%%%%%%%%%%%%%%%%%%%

\iftoggle{LatexFull}{

\section{Extended version}
\label{secExtend}
In this part of the article we prove Lemmas \ref{lemSterbenz}
and \ref{lemConvexity},
the corollaries which were not proved in the previous sections,
and the propositions. We try to prove every
assertion we make, no matter how trivial it may sound.
In all propositions $\wfpf$ is a floating point
system, $z \in \wrone{}$, $x \in \wfpf$,
$\wflm$ rounds to nearest in $\wfpf$,
and $u$, $\wfpemin{}$, $\mu$ $\wfpa$ and $\wfpnu$ are the numbers
related to this system in Definitions \ref{longDefEps},
\ref{longDefPerfect}, \ref{longDefMPFR}, \ref{longDefIEEE},
\ref{longDefAlpha} and \ref{longDefNu}.

\subsection{Proofs of Lemmas \ref{lemSterbenz} and \ref{lemConvexity}}
\label{secExtLem}

In this section we prove Lemmas \ref{lemSterbenz} and \ref{lemConvexity}.\\
%
%
%%%%%%%%%%%%%%%%%%%%%%%%%%%%%%%%%%%%%%%%%%%%%%%%%%%%%%%%%%%%%%%%%%%%%%%
% Proof of Sterbenz lemma
%%%%%%%%%%%%%%%%%%%%%%%%%%%%%%%%%%%%%%%%%%%%%%%%%%%%%%%%%%%%%%%%%%%%%%%

\pbProofB{Lemma}{lemSterbenz}
If $b - a < \beta \nu$ then Lemma \ref{lemSterbenz} follows from
Lemma \ref{lemSmallSum}. Therefore, we can assume that $b - a \geq \beta \nu$.
Prop. \ref{propNormalForm} implies that $a = \beta^{d}\wlr{\beta^{\mu} + r}$
and $b = \beta^{e} \wlr{\beta^{\mu} + s}$ with $d,e \in \wz{}$ and
$r,s \in [0,\wlr{\beta - 1} \beta^{\mu})$. Since $a \leq b \leq 2 a$
and $\beta \geq 2$,
\[
\beta^{d}\wlr{\beta^{\mu} + r} \leq \beta^{e} \wlr{\beta^{\mu} + s} \leq
2 \beta^{d}\wlr{\beta^{\mu} + r} \leq \beta^{d + 1}\wlr{\beta^{\mu} + r}.
\]
Prop. \ref{propOrder} shows that $d \leq e \leq d + 1$ and
either (i) $e = d$ or (ii) $e = d + 1$. In case (i)
$b - a = \beta^{e} \wlr{s - r} \geq \beta \nu$.
Since $0 \leq s - r < \wlr{\beta - 1} \beta^{\mu}$ and $b - a \geq \nu$,
Prop. \ref{propNu} implies that $b - a \in \wfpf$.
In case (ii) $0 < b - a = \beta^{d} t$
for $t := \wlr{ \wlr{\beta  - 1 } \beta^{\mu} + \beta s - r} > 0$ and
\[
b - a \leq a \Rightarrow t \leq \beta^{\mu} + r < \beta^{1 + \mu}.
\]
This bound, the assumption $b - a \geq \beta \nu$ and Prop.
\ref{propNu} imply that $z \in \wfpf$.
\peProof{Lemma}{lemSterbenz}\\

%%%%%%%%%%%%%%%%%%%%%%%%%%%%%%%%%%%%%%%%%%%%%%%%%%%%%%%%%%%%%%%%
% End of the Proof of lemma lemSterbenz
%%%%%%%%%%%%%%%%%%%%%%%%%%%%%%%%%%%%%%%%%%%%%%%%%%%%%%%%%%%%%%%%
%
%

%
%
%%%%%%%%%%%%%%%%%%%%%%%%%%%%%%%%%%%%%%%%%%%%%%%%%%%%%%%%%%%%%%%%%%%%%%%%%
% Proof of lemma lemConvexity
%%%%%%%%%%%%%%%%%%%%%%%%%%%%%%%%%%%%%%%%%%%%%%%%%%%%%%%%%%%%%%%%%%%%%%%%%
\pbProofB{Lemma}{lemConvexity}
The function $g_k$ has first derivative
\[
\wdf{g_k}{u} = -\wfc{g_k}{u} \sum_{i = 1}^k \frac{n_i}{1 + n_i u}
\]
and second derivative
\[
\wdsf{g_k}{u} = \wfc{g_k}{u} \wlr{ \wlr{\sum_{i = 1}^k \frac{n_i}{1 + n_i u}}^2 + \sum_{i = 1}^k \frac{n_i^2}{\wlr{1 + n_i u}^2}} > 0,
\]
and, therefore, it is convex. Similarly, the function $f_k$ has first derivative
\[
\wdf{f_k}{u} = \wfc{f_k}{u} \sum_{i = 1}^k \frac{n_i}{\wlr{1 + n_i u} \wlr{1 + 2 n_i u}}
\]
and second derivative
\[
\wdsf{f_k}{u} = \wfc{f_k}{u} \wlr{ \wlr{\sum_{i = 1}^k \frac{n_i}{\wlr{1 + n_i u} \wlr{1 + 2 n_i u}}}^2 -
\sum_{i=1}^k \frac{n_i^2 \wlr{3 + 4 n_i u}}{\wlr{\wlr{1 + n_i u} \wlr{1 + 2 n_i u}}^2}}.
\]
It follows that
\pbDef{dsf}
\wdsf{f_k}{u} =
-\wfc{f_k}{u} \wvec{v}^\wtr{} \wlr{3 \wvec{I} - \mathds{1} \mathds{1}^{\wtr{}}} \wvec{v}
- 4 \wfc{f_k}{u} u \sum_{i=1}^k \frac{n_i^3}{\wlr{\wlr{1 + n_i u} \wlr{1 + 2 n_i u}}^2},
\peDef{dsf}
where $\wvec{I}$ is the $k \times k$ identity matrix, $\mathds{1} \in \wrn{k}$ is the vector
with all entries equal to $1$ and $\wvec{v} \in \wrn{k}$ has entries
\[
v_i := \frac{n_i}{\wlr{1 + n_i u} \wlr{1 + 2 n_i u}}.
\]
The $k \times k$ symmetric matrix $\wvec{M} = 3 \wvec{I} - \mathds{1} \mathds{1}^{\wtr{}}$ has a $(k-1)$
dimensional eigenspace associated to
the eigenvalue $3$ which is orthogonal to $\mathds{1}$,
and $\mathds{1}$ is an eigenvector with
eigenvalue $3 - k$. Therefore, $\wvec{M}$ is positive semidefinite for $k \leq 3$,
Equation \pRef{dsf} implies that
\[
\wdsf{f_k}{u} \leq
- 4 \wfc{f_k}{u} u \sum_{i=1}^k \frac{n_i^3}{\wlr{\wlr{1 + n_i u} \wlr{1 + 2 n_i u}}^2} < 0
\]
and we are done.
\peProof{Lemma}{lemConvexity}\\

%%%%%%%%%%%%%%%%%%%%%%%%%%%%%%%%%%%%%%%%%%%%%%%%%%%%%%%%%%%%%%%%%%%%%%%
% End of Proof of lemma convexity
%%%%%%%%%%%%%%%%%%%%%%%%%%%%%%%%%%%%%%%%%%%%%%%%%%%%%%%%%%%%%%%%%%%%%%%
%
%

\subsection{Proofs of the remaining corollaries}
\label{secExtCor}
In this section we prove the corollaries which were
not proved in the previous sections.\\

%
%
%%%%%%%%%%%%%%%%%%%%%%%%%%%%%%%%%%%%%%%%%%%%%%%%%%%%%%%%%%%%%%%
% Proof of corollary corFourProd
%%%%%%%%%%%%%%%%%%%%%%%%%%%%%%%%%%%%%%%%%%%%%%%%%%%%%%%%%%%%%%%

\pbProofB{Corollary}{corFourProd} Corollary \ref{corFourProd} is
a consequence of the convexity of $\wlr{1 + u}^{-k}$ and
the concavity of $f^k$ for $k \leq 3$ and $f$ in \pRef{concF}, which yield
\[
1 - k u \leq \frac{1}{\wlr{1 + u}^k} \leq \wlr{\frac{1 + 2 u}{1 + u}}^k \leq 1 + k u
\]
for $k = 1$, $2$ and $3$.
\peProof{Corollary}{corFourProd}\\

%%%%%%%%%%%%%%%%%%%%%%%%%%%%%%%%%%%%%%%%%%%%%%%%%%%%%%%%%%%%%%%
% End of Proof of corollary corFourProd
%%%%%%%%%%%%%%%%%%%%%%%%%%%%%%%%%%%%%%%%%%%%%%%%%%%%%%%%%%%%%%%
%
%

%
%
%%%%%%%%%%%%%%%%%%%%%%%%%%%%%%%%%%%%%%%%%%%%%%%%%%%%%%%%%%%%%%%%%%%%%%%%%
% Proof of Corollary corNormOneIEEE
%%%%%%%%%%%%%%%%%%%%%%%%%%%%%%%%%%%%%%%%%%%%%%%%%%%%%%%%%%%%%%%%%%%%%%%%%

\pbProofB{Corollary}{corNormOneIEEE}
Let $\wflmtx$ be the rounding tuple in Prop. \ref{propRoundIEEEX}.
If the $y_k$ are floating point numbers
then $\wfpsumkf{k}{\wvec{y},\wflmt} = \wfpsumkf{k}{\wvec{y},\wflmtx}$ for all $k$
and Corollary \ref{corNormOneIEEE} follows from Lemma
\ref{lemNormOneBound}.
\peProof{Corollary}{corNormOneIEEE}\\

%%%%%%%%%%%%%%%%%%%%%%%%%%%%%%%%%%%%%%%%%%%%%%%%%%%%%%%%%%%%%%%%%%%%%%%%%
% End of proof of corollary SharpSumNearIEEE
%%%%%%%%%%%%%%%%%%%%%%%%%%%%%%%%%%%%%%%%%%%%%%%%%%%%%%%%%%%%%%%%%%%%%%%%%
%
%

%
%
%%%%%%%%%%%%%%%%%%%%%%%%%%%%%%%%%%%%%%%%%%%%%%%%%%%%%%%%%%%%%%%%%%%%%%%%%
% Proof of Corollary corNormOneMFPR
%%%%%%%%%%%%%%%%%%%%%%%%%%%%%%%%%%%%%%%%%%%%%%%%%%%%%%%%%%%%%%%%%%%%%%%%%

\pbProofB{Corollary}{corNormOneMPFR}
Let $\wflmtx$ be the rounding tuple in Prop. \ref{propRoundMPFRX}.
If all $y_k$ are non negative floating point numbers
then $\wfpsumkf{k}{\wvec{y},\wflmt} = \wfpsumkf{k}{\wvec{y},\wflmtx}$ for all $k$ and
Corollary \ref{corNormOneMPFR} follows from Lemma
\ref{lemNormOneBound}.
\peProof{Corollary}{corNormOneMPFR} \\

%%%%%%%%%%%%%%%%%%%%%%%%%%%%%%%%%%%%%%%%%%%%%%%%%%%%%%%%%%%%%%%%%%%%%%%%%
% end of Proof of Corollary corShaprSumNearMFPR
%%%%%%%%%%%%%%%%%%%%%%%%%%%%%%%%%%%%%%%%%%%%%%%%%%%%%%%%%%%%%%%%%%%%%%%%%
%
%

%
%
%%%%%%%%%%%%%%%%%%%%%%%%%%%%%%%%%%%%%%%%%%%%%%%%%%%%%%%%%%%%%%%%%%%%%%%%%
% Proof of Corollary corPositiveUnperfect
%%%%%%%%%%%%%%%%%%%%%%%%%%%%%%%%%%%%%%%%%%%%%%%%%%%%%%%%%%%%%%%%%%%%%%%%%

\pbProofB{Corollary}{corPositiveUnperfect}
If $\wfpf$ is a MPFR system, let  $\wflmtx$ be
the rounding tuple in Prop. \ref{propRoundMPFRX}.
Since all $y_k$ belong to $\wfpsm$ and are non negative
we have that $\wfpsumkf{k}{\wvec{y},\wflmt} = \wfpsumkf{k}{\wvec{y},\wflmtx}$ for all $k$
and Corollary \ref{corPositiveUnperfect} follows from Lemma
\ref{lemPositiveBound}.
If $\wfpf$ is an IEEE system,
let  $\wflmtx$ be rounding tuple in Prop. \ref{propRoundIEEEX}.
Since all $y_k$ are floating point numbers,
$\wfpsumkf{k}{\wvec{y},\wflmt} = \wfpsumkf{k}{\wvec{y},\wflmtx}$ for all $k$
and Corollary \ref{corPositiveUnperfect} follows from Lemma
\ref{lemPositiveBound}.
\peProof{Corollary}{corPositiveUnperfect}\\

%%%%%%%%%%%%%%%%%%%%%%%%%%%%%%%%%%%%%%%%%%%%%%%%%%%%%%%%%%%%%%%%%%%%%%%%%
% End of proof of corollary SharpSumNearIEEEC
%%%%%%%%%%%%%%%%%%%%%%%%%%%%%%%%%%%%%%%%%%%%%%%%%%%%%%%%%%%%%%%%%%%%%%%%%
%
%

%
%
%%%%%%%%%%%%%%%%%%%%%%%%%%%%%%%%%%%%%%%%%%%%%%%%%%%%%%%%%%%%%%%%%%%%%%%%%
% Proof of Corollary corFmaPerfect
%%%%%%%%%%%%%%%%%%%%%%%%%%%%%%%%%%%%%%%%%%%%%%%%%%%%%%%%%%%%%%%%%%%%%%%%%

\pbProofB{Corollary}{corFmaPerfect}
In a perfect system, the dot product of $n + 1$ numbers evaluated
using a fma, as in Definition \ref{longDefFmaDot},
is the floating point sum of the $(n + 2)$ real numbers
$p_0 := 0$ and $p_k := x_{k-1} y_{k-1}$ for $k > 0$,
 and Equation \pRef{thSharpFmaDotNear} follows from
 Lemma \ref{lemNormOneBound} applied to the $p_k$.
\peProof{Corollary}{corFmaPerfect}\\

%
%
%%%%%%%%%%%%%%%%%%%%%%%%%%%%%%%%%%%%%%%%%%%%%%%%%%%%%%%%%%%%%%%%%%%%%%%%%
% Proof of Corollary corFmaUnperfect
%%%%%%%%%%%%%%%%%%%%%%%%%%%%%%%%%%%%%%%%%%%%%%%%%%%%%%%%%%%%%%%%%%%%%%%%%

\pbProofB{Corollary}{corFmaUnperfect}
In an unperfect systems, the dot product of $n + 1$ numbers evaluated
using a fma, as in Definition \ref{longDefFmaDot},
is the floating point sum of the $(n + 2)$ real numbers
$p_0 := 0$ and $p_k := x_{k-1} y_{k-1}$ for $k > 0$,
 and Corollary \ref{corFmaUnperfect} follows from
 Corollary \ref{corNormOneUnperfect} applied to the $p_k$.
\peProof{Corollary}{corFmaUnperfect}\\

%%%%%%%%%%%%%%%%%%%%%%%%%%%%%%%%%%%%%%%%%%%%%%%%%%%%%%%%%%%%%%%%%%%%%%%%%
% End of Proof of Corollary corFmaUnperfect
%%%%%%%%%%%%%%%%%%%%%%%%%%%%%%%%%%%%%%%%%%%%%%%%%%%%%%%%%%%%%%%%%%%%%%%%%
%
%

%
%
%%%%%%%%%%%%%%%%%%%%%%%%%%%%%%%%%%%%%%%%%%%%%%%%%%%%%%%%%%%%%%%%%%%%%%%%%
% Proof of Corollary corDotPerfect
%%%%%%%%%%%%%%%%%%%%%%%%%%%%%%%%%%%%%%%%%%%%%%%%%%%%%%%%%%%%%%%%%%%%%%%%%

\pbProofB{Corollary}{corDotPerfect}
The dot product is the floating point sum of the floating point numbers
$p_k := \wflrk{r}{k}{x_k y_k}$. In a perfect system, Lemma
\ref{lemUNear} shows that
\[
p_k = x_k y_k + \theta_k \frac{u}{1 + u} x_k y_k
\hspace{1cm} \wrm{with} \hspace{1cm} \wabs{\theta_k} \leq 1,
\]
and Lemma \ref{lemNormOneBound} implies that
\[
\wabs{\wflt{\sum_{k = 0}^n p_k} - \sum_{k = 0}^n p_k}
= \frac{n u}{1 + n u} \sum_{k = 0}^n \wabs{p_k}.
\]
It follows that
\[
\wabs{\wflt{\sum_{k = 0}^n x_k y_k} - \sum_{k = 0}^n x_k y_k }
\leq
\sum_{k = 0}^n \wabs{p_k - x_k} +
\wabs{\wflt{\sum_{k = 0}^n x_k y_k} - \sum_{k = 0}^n p_k}
\leq \beta_n u \sum_{k = 0}^n \wabs{x_k y_k}
\]
for
\[
\beta_n := \frac{1}{1 + u} \wlr{1 + \frac{n}{1 + n u} \wlr{1 + 2 u}}
= \frac{n + 1 + 3 n u}{1 + \wlr{n    + 1} u + n u^2}.
\]
Finally, note that for $n \geq 1$ and $20 n u \leq 1$,
\[
\beta_n - \frac{n + 1}{1 + n u/2} =
- u \frac{\wlr{n - 2} \wlr{n - 1 - n u}}{\wlr{1 + n u/2} \wlr{1 + \wlr{n+1}u + n u^2}} \leq 0,
\]
and
\[
\beta_n - \frac{n + 1}{1 + \wlr{n - 3} u} =
- u \frac{n + 4 + 10 n u - 2 n^2 u}{\wlr{1 + \wlr{n + 1} u + n u^2}\wlr{1 + \wlr{n - 3} u}} < 0.
\]
\peProof{Corollary}{corDotPerfect}\\

%%%%%%%%%%%%%%%%%%%%%%%%%%%%%%%%%%%%%%%%%%%%%%%%%%%%%%%%%%%%%%%%%%%%%%%%%
% End of oroof of Corollary corDotPerfect
%%%%%%%%%%%%%%%%%%%%%%%%%%%%%%%%%%%%%%%%%%%%%%%%%%%%%%%%%%%%%%%%%%%%%%%%%
%
%

%
%
%%%%%%%%%%%%%%%%%%%%%%%%%%%%%%%%%%%%%%%%%%%%%%%%%%%%%%%%%%%%%%%%%%%%%%%%%
% Proof of Corollary corDotIEEE
%%%%%%%%%%%%%%%%%%%%%%%%%%%%%%%%%%%%%%%%%%%%%%%%%%%%%%%%%%%%%%%%%%%%%%%%%

\pbProofB{Corollary}{corDotIEEE}
The dot product is the sum of the
$n + 1$ floating point numbers $p_k := \wflrk{r}{k}{x_k y_k}$, and
 Corollary \ref{corNormOneIEEE} shows that
\[
\wabs{\wflt{\sum_{k = 0}^n p_k} - \sum_{k = 0}^n p_k} \leq \frac{n u}{1 + n u} \sum_{k = 0}^n \wabs{p_k}.
\]
We also have
\[
\wabs{p_k - x_k y_k} \leq \frac{u}{1 + u} \wabs{x_k y_k} + \frac{\wfpa}{2}
\]
and
\[
\wabs{\wflt{\sum_{k = 0}^n x_k y_k} - \sum_{k = 0}^n x_k y_k} \leq
\wabs{\wflt{\sum_{k = 0}^n p_k} - \sum_{k = 0}^n p_k} + \sum_{k = 0}^n \wabs{p_k - x_k y_k}
\]
\[
\leq \frac{n u}{1 + n u} \sum_{k = 0}^n \wabs{p_k} + \sum_{k = 0}^n \wabs{p_k - x_k y_k}
\]
\[
\leq \frac{n u}{1 + n u} \wlr{ \frac{1 + 2 u}{1 + u} \sum_{k = 0}^n \wabs{x_k y_k} + \wlr{n + 1} \frac{\alpha}{2}} +
\frac{\wlr{n+1} \alpha}{2} + \frac{u}{1 + u} \sum_{k = 0}^n \wabs{x_k y_k}
\]
\[
= \beta_n u \sum_{k = 0}^n \wabs{x_k y_k} + b \frac{\alpha}{2}
\]
for $\beta_n$ in Corollary \ref{corDotPerfect} and
\[
b := \wlr{n + 1} \wlr{1 + \frac{n u}{1 + n u}} = \wlr{n + 1} \frac{1 + 2 n u}{1 + n u} < 1.05 \wlr{n + 1},
\]
because $20 n u \leq 1$. Finally, if $ u \sum_{k = 0}^n \wabs{x_k y_k} \geq \wfpa$ then
\[
\beta_n u \sum_{k = 0}^n \wabs{x_k y_k} + b \frac{\alpha}{2} \leq
\theta_n u \sum_{k = 0}^n \wabs{x_k y_k}
\hspace{1cm} \wrm{for} \hspace{1cm}
\theta_n := \beta_n + \frac{n + 1}{2} \frac{1 + 2 n u}{1 + n u},
\]
and the software Mathematica shows that
\[
\theta_n - 3 \frac{n +1}{2} = - u \frac{n^2 - 3 n + 2 + n u \wlr{1 + n}}{2 \wlr{1 + u} \wlr{1 + n u}}
\]
which is negative for $n \geq 1$. This proves
the last equation in Corollary \ref{corDotIEEE}.
\peProof{Corollary}{corDotIEEE}\\

%
%
%%%%%%%%%%%%%%%%%%%%%%%%%%%%%%%%%%%%%%%%%%%%%%%%%%%%%%%%%%%%%%%%%%%%%%%%%
% End of Proof of Corollary corDotIEEE
%%%%%%%%%%%%%%%%%%%%%%%%%%%%%%%%%%%%%%%%%%%%%%%%%%%%%%%%%%%%%%%%%%%%%%%%%

%
%
%%%%%%%%%%%%%%%%%%%%%%%%%%%%%%%%%%%%%%%%%%%%%%%%%%%%%%%%%%%%%%%%%%%%%%%%%
% Proof of Corollary corDotMPFR
%%%%%%%%%%%%%%%%%%%%%%%%%%%%%%%%%%%%%%%%%%%%%%%%%%%%%%%%%%%%%%%%%%%%%%%%%

\pbProofB{Corollary}{corDotMPFR}
The dot product is the sum of the
$n + 1$ floating point numbers $p_k := \wflrk{r}{k}{x_k y_k}$, and
the proof of Corollary \ref{corNormOneUnperfect} shows that
\[
\wabs{\wflt{\sum_{k = 0}^n p_k} - \sum_{k = 0}^n p_k} \leq \frac{m  \wfpa}{2} +
\frac{\wlr{n - m} u}{1 + \wlr{n - m} u} \wlr{\frac{m  \wfpa}{2} + \sum_{k = 0}^n \wabs{p_k}},
\]
for some $m \in [0,n]$. We also have that
\[
\wabs{p_k - x_k y_k} \leq \frac{u}{1 + u} \wabs{x_k y_k} + \frac{\wfpa}{2}
\]
and
\[
\wabs{\wflt{\sum_{k = 0}^n x_k y_k} - \sum_{k = 0}^n x_k y_k} \leq
\wabs{\wflt{\sum_{k = 0}^n p_k} - \sum_{k = 0}^n p_k} + \sum_{k = 0}^n \wabs{p_k - x_k y_k}
\]
\[
\leq \frac{m \wfpa}{2}  +
\frac{\wlr{n - m} u}{1 + \wlr{n - m} u}
\wlr{\frac{m \wfpa}{2} + \sum_{k = 0}^n \wabs{p_k}} + \sum_{k = 0}^n \wabs{p_k - x_k y_k}
\]
\pbDef{MPRFBA}
\leq \frac{\wlr{n - m} u}{1 + \wlr{n - m} u}
\wlr{ \frac{1 + 2 u}{1 + u} \sum_{k = 0}^n \wabs{x_k y_k} + \wlr{m + n + 1} \frac{\wfpa}{2}} +
\peDef{MPRFBA}
\[
\frac{\wlr{m + n+1} \wfpa}{2} + \frac{u}{1 + u} \sum_{k = 0}^n \wabs{x_k y_k}
 \leq \beta_n u \sum_{k = 0}^n \wabs{x_k y_k} + b \frac{\wfpa}{2}
\]
for $\beta_n$ in Corollary \ref{corDotPerfect} and
\[
b := \frac{n^2 + n - m^2 -m}{1 + \wlr{n - m} u} u + \wlr{m + n + 1} \leq \frac{n \wlr{n + 1} u}{1 + n u} + 2 n + 1 \leq 2.05n + 1.05.
\]
Finally, if $ u \sum_{k = 0}^n \wabs{x_k y_k} \geq \wfpa$,
then
\[
\wabs{\wflt{\sum_{k = 0}^n x_k y_k} - \sum_{k = 0}^n x_k y_k} \leq
\gamma_n u \sum_{k = 0}^n \wabs{x_k y_k}
\]
for
\[
\gamma_n := \beta_n + \frac{1}{2} \wlr{\frac{n^2 + n - m^2 -m}{1 + \wlr{n - m} u} u + \wlr{m + n + 1}}.
\]
The derivative of $\gamma_n$ with respect to $m$ is
\[
- \frac{1 + 2 n u + u}{\wlr{1 + \wlr{n - mu} u}^2} < 0
\]
and $\gamma_n$ is maximized for $m = 0$, in which case
it is equal to the $\theta_n$ in the proof of Corollary \ref{corDotIEEE}.
This proves the last statement in Corollary \ref{corDotMPFR}.
\peProof{Corollary}{corDotMPFR}\\

%%%%%%%%%%%%%%%%%%%%%%%%%%%%%%%%%%%%%%%%%%%%%%%%%%%%%%%%%%%%%%%%%%%%%%%%%
% End of the Proof of Corollary corDotMPFR
%%%%%%%%%%%%%%%%%%%%%%%%%%%%%%%%%%%%%%%%%%%%%%%%%%%%%%%%%%%%%%%%%%%%%%%%%
%
%

\subsection{Numbers}
\label{secNumbers}
This section contains new propositions about real and integer numbers,
and the proofs of propositions related to these numbers stated
in the main part of the article.

\subsubsection{Propositions}
\label{secNumbersProps}
This sections presents more propositions regarding real and integer numbers.

%
%
%%%%%%%%%%%%%%%%%%%%%%%%%%%%%%%%%%%%%%%%%%%%%%%%%%%%%%%%%%%%%%%%%%%%%%%%%%%%%%%%%%%
% propNormalFormCont
%%%%%%%%%%%%%%%%%%%%%%%%%%%%%%%%%%%%%%%%%%%%%%%%%%%%%%%%%%%%%%%%%%%%%%%%%%%%%%%%%%%

\pbPropBT{propNormalFormCont}{Continuity of the normal form}
If $e$ is integer, $\wabs{z} = \beta^{e} \wlr{\beta^{\mu} + w}$ with
$0 < w < \wlr{\beta - 1} \beta^{\mu}$ and
\[
\wabs{y - z} < \beta^{e} \min \wset{w, \wlr{\beta - 1} \beta^{\mu} - w}
\]
then $y = \wsign{z} \beta^{e} \wlr{\beta^{\mu} + v}$ with $0 < v < \wlr{\beta - 1} \beta^{\mu}$
and $\wabs{v - w} = \beta^{-e} \wabs{y - z}$.
\peFullProp{propNormalFormCont}

%
%
%%%%%%%%%%%%%%%%%%%%%%%%%%%%%%%%%%%%%%%%%%%%%%%%%%%%%%%%%%%%%%%%%%%%%%%%%%%%%%%%%%%
% propNormalFormDis
%%%%%%%%%%%%%%%%%%%%%%%%%%%%%%%%%%%%%%%%%%%%%%%%%%%%%%%%%%%%%%%%%%%%%%%%%%%%%%%%%%%

\pbPropBT{propNormalFormDis}{Discontinuity of the normal form}
If $e$ is integer and $\wabs{z} = \beta^{e + \mu}$ with $\wabs{y - z} < \beta^{e + \mu - 1} \wlr{\beta - 1}$
then we have three possibilities:
\begin{itemize}
\item[(i)] $\wabs{y} < \wabs{z}$ and
$y = \wsign{z} \beta^{e - 1} \wlr{\beta^{\mu} + w}$ with
\[
0 < v = \wlr{\beta - 1}\beta^{\mu} - \beta^{1 - e} \wabs{y - z} < \wlr{\beta - 1}\beta^{\mu}.
\]
\item[(ii)] $\wabs{y} = \wabs{z}$ and $y = z$.
\item[(iii)] $\wabs{y} > \wabs{z}$ and
$y = \wsign{z} \beta^{e} \wlr{\beta^{\mu} + w}$ with $0 < w = \beta^{-e} \wabs{y - z} <  \beta^{\mu - 1} \wlr{\beta - 1}$.
\end{itemize}
\peProp{propNormalFormDis}

\subsubsection{Proofs}
\label{secNumbersProof}
In this section we prove the propositions regarding integer and real numbers.\\

%
%
%%%%%%%%%%%%%%%%%%%%%%%%%%%%%%%%%%%%%%%%%%%%%%%%%%%%%%%%%%%%%%%%%%%%%%%%%%%%%%%%%%%
% Proof of prop propOrder
%%%%%%%%%%%%%%%%%%%%%%%%%%%%%%%%%%%%%%%%%%%%%%%%%%%%%%%%%%%%%%%%%%%%%%%%%%%%%%%%%%%%

\pbProofB{Proposition}{propOrder} Since $d,e\in \wz{}$ and $d < e$ we have
that $e - d \geq 1$ and
\[
\beta^{e} \wlr{\beta^{\wfpmu} + w} - \beta^{d} \wlr{\beta^{\wfpmu} + v} \geq
\beta^{d} \wlr{\wlr{\beta^{e - d} - 1} \beta^{\wfpmu} - v} \geq
\beta^{d} \wlr{\wlr{\beta - 1} \beta^{\wfpmu} - v} > 0,
\]
and this shows that $\beta^{e} \wlr{\beta^{\wfpmu} + w} > \beta^{d} \wlr{\beta^{\wfpmu} + v}$.
\peProof{Proposition}{propOrder}\\

%%%%%%%%%%%%%%%%%%%%%%%%%%%%%%%%%%%%%%%%%%%%%%%%%%%%%%%%%%%%%%%%%%%%%%%%%%%%%%%%%%%
% End of proof of propOrder
%%%%%%%%%%%%%%%%%%%%%%%%%%%%%%%%%%%%%%%%%%%%%%%%%%%%%%%%%%%%%%%%%%%%%%%%%%%%%%%%%%%
%
%

%
%
%%%%%%%%%%%%%%%%%%%%%%%%%%%%%%%%%%%%%%%%%%%%%%%%%%%%%%%%%%%%%%%%%%%%%%%%%%%%%%%%%%%
% Proof of prop. propNormalForm
%%%%%%%%%%%%%%%%%%%%%%%%%%%%%%%%%%%%%%%%%%%%%%%%%%%%%%%%%%%%%%%%%%%%%%%%%%%%%%%%%%%

\pbProofB{Proposition}{propNormalForm}
The integer exponent $e := \wfloor{\wfc{\log_{\beta}}{\wabs{z}}} - \wfpmu$ satisfies
\[
\wfc{\log_{\beta}}{\wabs{z}}  - \wfpmu - 1 < e \leq \wfc{\log_{\beta}}{\wabs{z}} - \wfpmu
\hspace{1.0cm} \wrm{and} \hspace{1.0cm}
{\beta}^{- \wfpmu - 1} \wabs{z} < {\beta}^{e} \leq \wabs{z} {\beta}^{-\wfpmu}.
\]
The equation above shows that
$w :=  {\beta}^{-e} \wabs{z} - {\beta}^{\wfpmu}$  satisfies $0 \leq w < \wlr{\beta - 1}{\beta}^{\wfpmu}$
and $z = \wsign{z} {\beta}^{e}\wlr{{\beta}^{\wfpmu} + w}$. If
 $z = \wsign{z} {\beta}^{d} \wlr{{\beta}^{\wfpmu} + v}$
with $d \in \wz{}$ and $0 \leq v < \wlr{\beta - 1} {\beta}^{\wfpmu}$ then
\[
{\beta}^{e} \wlr{\beta^{\wfpmu} + w} = \wabs{z} = {\beta}^{d} \wlr{{\beta}^{\wfpmu} + v},
\]
and Prop. \ref{propOrder} implies that $d = e$, and the equation above implies that $v = w$.
\peProof{Proposition}{propNormalForm}\\

%%%%%%%%%%%%%%%%%%%%%%%%%%%%%%%%%%%%%%%%%%%%%%%%%%%%%%%%%%%%%%%%%%%%%%%%%%%%%%%%%%%
% End of proof of prop. propNormalForm
%%%%%%%%%%%%%%%%%%%%%%%%%%%%%%%%%%%%%%%%%%%%%%%%%%%%%%%%%%%%%%%%%%%%%%%%%%%%%%%%%%%
%
%

%%%%%%%%%%%%%%%%%%%%%%%%%%%%%%%%%%%%%%%%%%%%%%%%%%%%%%%%%%%%%%%%%%%%%%%%%%%%%%%%%%%
% End of proof of prop. propIntegerForm
%%%%%%%%%%%%%%%%%%%%%%%%%%%%%%%%%%%%%%%%%%%%%%%%%%%%%%%%%%%%%%%%%%%%%%%%%%%%%%%%%%%
%
%

%
%
%%%%%%%%%%%%%%%%%%%%%%%%%%%%%%%%%%%%%%%%%%%%%%%%%%%%%%%%%%%%%%%%%%%%%%%%%%%%%%%%%%%
% propNormalFormCont
%%%%%%%%%%%%%%%%%%%%%%%%%%%%%%%%%%%%%%%%%%%%%%%%%%%%%%%%%%%%%%%%%%%%%%%%%%%%%%%%%%%

\pbProofB{Proposition}{propNormalFormCont}
We have that
\[
\wabs{1 - \frac{y}{z}} < \frac{\beta^{e} w}{\wabs{z}} = \frac{w}{\beta^{\mu} + w} < 1
\ \ \Rightarrow \frac{y}{z} > 0 \Rightarrow y \neq 0,
\]
and Prop. \ref{propNormalForm} yield $d$ and $v \in [0,\wlr{\beta - 1} \beta^{\mu})$ such
that $y = \wsign{y} \beta^{d} \wlr{\beta^{\mu} + v}$. The inequality
\[
\frac{\wsign{y} \beta^{d}\wlr{\beta^{\mu} + v}}{\wsign{z} \beta^{e}\wlr{\beta^{\mu} + w}} = \frac{y}{z}  > 0
\]
implies that $\wsign{y} = \wsign{z}$. Moreover,
\[
\beta^{d} \wlr{\beta^{\mu} + v} = \wabs{y} \leq \wabs{z} + \wabs{y-z} < \wabs{z} + \beta^{e} \wlr{\wlr{\beta - 1}\beta^{\mu} - w} = \beta^{e + 1 +\mu}
\]
and Prop. \ref{propOrder} implies that $d \leq e$. Similarly,
\[
\beta^{d} \wlr{\beta^{\mu} + v} = \wabs{y}\geq \wabs{z} - \wabs{y-z} > \wabs{z} - \beta^{e} w = \beta^{e + \mu},
\]
and $d \geq e$. Therefore $d = e$,
$y = \wsign{y} \beta^{e} \wlr{\beta^{\mu} + v}$ and
$\wabs{y - z} =\beta^{e} \wabs{w - z}$.
\peProof{Proposition}{propNormalFormCont}\\

%%%%%%%%%%%%%%%%%%%%%%%%%%%%%%%%%%%%%%%%%%%%%%%%%%%%%%%%%%%%%%%%%%%%%%%%%%%%%%%%%%%
% End of proof of propNormalFormCont
%%%%%%%%%%%%%%%%%%%%%%%%%%%%%%%%%%%%%%%%%%%%%%%%%%%%%%%%%%%%%%%%%%%%%%%%%%%%%%%%%%%
%
%

%
%
%%%%%%%%%%%%%%%%%%%%%%%%%%%%%%%%%%%%%%%%%%%%%%%%%%%%%%%%%%%%%%%%%%%%%%%%%%%%%%%%%%%
% propNormalFormDis
%%%%%%%%%%%%%%%%%%%%%%%%%%%%%%%%%%%%%%%%%%%%%%%%%%%%%%%%%%%%%%%%%%%%%%%%%%%%%%%%%%%

\pbProofB{Proposition}{propNormalFormDis}
We have that
\[
\wabs{1 - \frac{y}{z}} < \frac{\beta - 1}{\beta} < 1
\ \ \Rightarrow \frac{y}{z} > 0 \Rightarrow y \neq 0,
\]
and Prop \ref{propNormalForm} yields $d \in \wz{}$ and $w \in [0,\wlr{\beta - 1} \beta^{\mu})$ such
that $y = \wsign{y} \beta^{d} \wlr{\beta^{\mu} + w}$. The inequality
\[
\frac{\wsign{y} \beta^{d}\wlr{\beta^{\mu} + w}}{\wsign{z} \beta^{e + \mu}} = \frac{y}{z}  > 0
\]
implies that $\wsign{y} = \wsign{z}$. We also have that
\[
\beta^{d} \wlr{\beta^{\mu} + w} = \wabs{y} \leq \wabs{z} + \wabs{y-z} < \wabs{z} + \beta^{e + \mu - 1} \wlr{\beta - 1}
= \beta^{e} \wlr{\beta^\mu + \beta^{\mu - 1} \wlr{\beta - 1} }
\]
and Prop. \ref{propOrder} implies that $d \leq e$. It $\wabs{y} \geq \wabs{z}$ then
Prop. \ref{propOrder} implies that $d \geq e$. It follows that $d = e$ and the conditions
in items (ii) and (iii) in Prop. \ref{propNormalFormDis} are satisfied. If $\wabs{y} < \wabs{z}$ then
Prop. \ref{propOrder} implies that $d < e$ and
\[
\beta^{d} \wlr{\beta^{\mu} + w} = \wabs{y} \geq \wabs{z} - \beta^{e + \mu - 1} \wlr{\beta - 1} =
\beta^{e - 1} \wlr{\beta^{\mu + 1} - \wlr{\beta^{\mu + 1} - \beta^{\mu}}} =
\beta^{e - 1 + \mu}
\]
and Prop. \ref{propOrder} imply that $d \geq e - 1$. Therefore $d = e - 1$
and the conditions in item (i) in Prop. \ref{propNormalFormDis} are satisfied.
\peProof{Proposition}{propNormalFormDis}\\

%%%%%%%%%%%%%%%%%%%%%%%%%%%%%%%%%%%%%%%%%%%%%%%%%%%%%%%%%%%%%%%%%%%%%%%%%%%%%%%%%%%
% End of proof of propNormalFormCont
%%%%%%%%%%%%%%%%%%%%%%%%%%%%%%%%%%%%%%%%%%%%%%%%%%%%%%%%%%%%%%%%%%%%%%%%%%%%%%%%%%%
%
%

\subsection{Floating point systems}
\label{secfloatSys}
In this section we present more definitions related to floating point systems
and more propositions about them. We prove the propositions regarding
floating point systems stated in the previous sections and the propositions
stated here. In most definitions, propositions and proofs in this section
$\wfpf$ is a floating point system,
$\wflm$ rounds to nearest in $\wfpf$, $z,w \in \wrone{}$ and $x,y \in \wfpf$,
and the numbers $\wfpa$ and $\wfpnu$ are as in Definitions
\ref{longDefAlpha} and \ref{longDefNu}, and the
exceptions are stated explicitly.

\subsubsection{Propositions}
\label{floatSysProps}
This section presents more propositions regarding floating point systems.
%
%
%%%%%%%%%%%%%%%%%%%%%%%%%%%%%%%%%%%%%%%%%%%%%%%%%%%%%%%%%%%%%%%%%%%%%%%%%%%%%%%%%%%
% propAlpha
%%%%%%%%%%%%%%%%%%%%%%%%%%%%%%%%%%%%%%%%%%%%%%%%%%%%%%%%%%%%%%%%%%%%%%%%%%%%%%%%%%%

\pbPropBT{propAlpha}{Minimality of alpha}
$\wfpa \in \wfpf$ and if $x \in \wfpf \setminus \wset{0}$ then $\wabs{x} \geq \wfpa$.
\peProp{propAlpha}

%
%%%%%%%%%%%%%%%%%%%%%%%%%%%%%%%%%%%%%%%%%%%%%%%%%%%%%%%%%%%%%%%%%%%%%%%%%%%%%%%%%%%
% propEmptyNormalRange
%%%%%%%%%%%%%%%%%%%%%%%%%%%%%%%%%%%%%%%%%%%%%%%%%%%%%%%%%%%%%%%%%%%%%%%%%%%%%%%%%%%

\pbPropBT{propEmptyNormalRange}{Empty normal range}
If $e$ is an exponent for $\wfpf$ and $r$ is an integer with $r \in [0,\wlr{\beta - 1} \beta^{\mu})$
then $\wfpf \cap \wlr{\beta^{e} \wlr{\beta^{\mu} + r}, \, \beta^{e} \wlr{\beta^{\mu} + r + 1}} = \emptyset$.
\peFullProp{propEmptyNormalRange}

%
%
%%%%%%%%%%%%%%%%%%%%%%%%%%%%%%%%%%%%%%%%%%%%%%%%%%%%%%%%%%%%%%%%%%%%%%%%%%%%%%%%%%%
% propEmptySubnormalRange
%%%%%%%%%%%%%%%%%%%%%%%%%%%%%%%%%%%%%%%%%%%%%%%%%%%%%%%%%%%%%%%%%%%%%%%%%%%%%%%%%%%

\pbPropBT{propEmptySubnormalRange}{Empty subnormal range}
Let $\wfpsi_{\wfpemin}$ be an IEEE system. If $r \in \wz{}$
and $-\beta^{\mu} \leq r < \beta^{\mu}$
then $\wlr{\beta^{\wfpemin} r, \, \beta^{\wfpemin} \wlr{r+ 1}} \cap \wfpsi_{\wfpemin} = \emptyset$.
\peFullProp{propEmptySubnormalRange}

%
%
%%%%%%%%%%%%%%%%%%%%%%%%%%%%%%%%%%%%%%%%%%%%%%%%%%%%%%%%%%%%%%%%%%%%%%%%%%%%%%%%%%%
% propRoundScaleSys
%%%%%%%%%%%%%%%%%%%%%%%%%%%%%%%%%%%%%%%%%%%%%%%%%%%%%%%%%%%%%%%%%%%%%%%%%%%%%%%%%%%

\pbPropBT{propSysScale}{Scale invariance}
If $\wfpf$ is perfect then  $x \in \wfpf$ if and only if $\beta x \in \wfpf$.
If $\wfpf$ is unperfect and $x \in \wfpf$ then $\beta x \in \wfpf$.
\peProp{propSysScale}

\subsubsection{Proofs}
\label{floatSysProofs}
In this section we prove the propositions regarding floating point systems.\\

%
%
%%%%%%%%%%%%%%%%%%%%%%%%%%%%%%%%%%%%%%%%%%%%%%%%%%%%%%%%%%%%%%%%%%%%%%%%%%%%%%%%%%%
% Proof of prop. propIntegerForm
%%%%%%%%%%%%%%%%%%%%%%%%%%%%%%%%%%%%%%%%%%%%%%%%%%%%%%%%%%%%%%%%%%%%%%%%%%%%%%%%%%%

\pbProofB{Proposition}{propIntegerForm}
According to Definitions \ref{longDefMPFR} and \ref{longDefIEEE} of MPFR system
and IEEE system, we have three possibilities: (i) $x = 0$, in which case
$x = \beta^{\wfpemin} r$ for $r = 0$,
(ii) $x$ is subnormal, and $\wfpf$ is an IEEE system
and $x = \beta^{\wfpemin} r$ with $\wabs{r} \in [1,\beta^{\mu}) \cap \wz{}$ and
(iii) $x \in \wfpbin_{e}$ for some $e \geq \wfpemin$,
and $x = \beta^{\wfpemin + e} r$ with $\wabs{r} \in [\beta^{\mu},\beta^{1 + \mu}) \cap \wz{}$.
\peProof{Proposition}{propIntegerForm}\\

%
%
%%%%%%%%%%%%%%%%%%%%%%%%%%%%%%%%%%%%%%%%%%%%%%%%%%%%%%%%%%%%%%%%%%%%%%%%%%%%%%%%%%%
% Proof of propSymmetry
%%%%%%%%%%%%%%%%%%%%%%%%%%%%%%%%%%%%%%%%%%%%%%%%%%%%%%%%%%%%%%%%%%%%%%%%%%%%%%%%%%%

\pbProofB{Proposition}{propSymmetry}
In the three possible cases, Definitions \ref{longDefPerfect}, \ref{longDefMPFR}
and \ref{longDefIEEE}, the floating point systems are clearly symmetric.
\peProof{Proposition}{propSymmetry}\\

%%%%%%%%%%%%%%%%%%%%%%%%%%%%%%%%%%%%%%%%%%%%%%%%%%%%%%%%%%%%%%%%%%%%%%%%%%%%%%%%%%%
% End of Proof of prop propSymmetry
%%%%%%%%%%%%%%%%%%%%%%%%%%%%%%%%%%%%%%%%%%%%%%%%%%%%%%%%%%%%%%%%%%%%%%%%%%%%%%%%%%%
%
%

%
%
%%%%%%%%%%%%%%%%%%%%%%%%%%%%%%%%%%%%%%%%%%%%%%%%%%%%%%%%%%%%%%%%%%%%%%%%%%%%%%%%%%%
% proof of propNu
%%%%%%%%%%%%%%%%%%%%%%%%%%%%%%%%%%%%%%%%%%%%%%%%%%%%%%%%%%%%%%%%%%%%%%%%%%%%%%%%%%%

\pbProofB{Proposition}{propNu}
If $\wfpf$ is a perfect system or a MPFR system and $x \in \wfpf \setminus \wset{0}$ then
$\wabs{x} \in \wfpbin_{e}$ for some exponent $e$ for $\wfpf$ by
Definitions \ref{longDefPerfect} and \ref{longDefMPFR}.
If $\wfpf$ is an IEEE system $\wfpsi_{\wfpemin}$ then
$\nu = \beta^{\wfpemin + \mu}$ and $x$ with $\wabs{x} \geq \nu$
is not subnormal. As a result, by definition of IEEE system, $\wabs{x} \in \wfpbin_{e}$ for
some exponent $e$ for $\wfpf$.  If $0 < \wabs{x} < \wfpnu$ then
$\wfpf$ is not perfect, because $\wfpnu = 0$ for perfect systems.
Moreover, $\wabs{x} \not \in \wfpbin_{e}$
for $e \geq \wfpemin$ and, by Definition \ref{longDefMPFR},
$\wfpf$ is not a MPFR system.
Therefore, $\wfpf$ is an IEEE system and $x$ is subnormal.

Regarding the converse part, if $r$ is a multiple of $\beta$
then we can replace $e$ by $e + 1$ and $r$ by $r/\beta$ and $z$ stays the same.
Therefore, we can assume that $r$ is not a multiple of $\beta$. In particular,
$\wabs{r} < \beta^{1 + \mu}$.
By symmetry (Prop. \ref{propSymmetry}), it suffices to show that $\wabs{z} \in \wfpf$
when $\wabs{z} \geq \wfpnu$. If $\wfpf$ is perfect then $\wabs{z} \in \wfpbin_{e}$
and Prop. \ref{propNu} holds. Therefore, we can assume that $\wfpf$ is unperfect.
In this case $\wfpnu = \beta^{\wfpemin + \mu}$ by Definition \ref{longDefNu} and
$\beta^{e} \wabs{r} \geq \beta^{\wfpemin + \mu}$;
actually $\beta^{e} \wabs{r} > \beta^{\wfpemin + \mu}$ because $r$ is not
a multiple of $\beta$. Since $0 < \wabs{r} < \beta^{1 + \mu}$,
there exists a first integer $d > 1$ such that $\beta^d \wabs{r} \geq \beta^{1 + \mu}$.
Dividing $\beta^{d-1} \wabs{r}$ by $\beta^{\mu}$ we obtain that
$\beta^{d-1} \wabs{r} = \beta^{\mu} q + p$ for $p,q \in \wz{}$ with $q \geq 0$ and $0 \leq p < \beta^{\mu}$.
The definition of $d$ yields $\beta^{1 + \mu} > \beta^{d - 1} \wabs{r} = \beta^{\mu} q + p$
and
\[
s := \wlr{q  - 1} \beta^\mu + p <  \wlr{\beta - 1} \beta^{\mu}.
\]
Moreover,
\[
\beta^{1 + \mu} q + \beta p = \beta^{d} \wabs{r} \geq \beta^{1 + \mu} \Rightarrow
q \geq 1 - p/\beta^\mu > 0 \ \ \Rightarrow q \geq 1 \ \ \Rightarrow s \geq 0.
\]
As a result, $\wabs{r} = \beta^{1 - d} \wlr{\beta^{\mu} + s}$
with $s \in [0, \wlr{\beta - 1} \beta^{\mu})$
and Prop. \ref{propOrder} leads to
\[
\wabs{z} \geq \wfpnu \Rightarrow \beta^{e + 1 - d} \wlr{\beta^{\mu} + s} \geq \beta^{\wfpemin} \wlr{\beta^\mu + 0}
\]
\[
\Rightarrow e + 1 - d \geq \wfpemin \Rightarrow \wabs{z} = \beta^{e + 1 - d} \wlr{\beta^{\mu} + s} \in \wfpbin_{e + 1 - d} \subset \wfpf.
\]

\peProof{Proposition}{propNu}\\

%%%%%%%%%%%%%%%%%%%%%%%%%%%%%%%%%%%%%%%%%%%%%%%%%%%%%%%%%%%%%%%%%%%%%%%%%%%%%%%%%%%
% End of proof of propNu
%%%%%%%%%%%%%%%%%%%%%%%%%%%%%%%%%%%%%%%%%%%%%%%%%%%%%%%%%%%%%%%%%%%%%%%%%%%%%%%%%%%
%
%

%
%
%%%%%%%%%%%%%%%%%%%%%%%%%%%%%%%%%%%%%%%%%%%%%%%%%%%%%%%%%%%%%%%%
% Proof of Prop propSubnormalSum
%%%%%%%%%%%%%%%%%%%%%%%%%%%%%%%%%%%%%%%%%%%%%%%%%%%%%%%%%%%%%%%%

\pbProofB{Proposition}{propSubnormalSum}
Prop. \ref{propSymmetry} states that
$x + y \in \wfpsi \Leftrightarrow -\wlr{x + y} \in \wfpsi$,
and it suffices to show that $\wabs{x + y} \in \wfpsi$.
Since $x$ and $y$ are subnormal,
$x = \wsign{x} \beta^{\wfpemin} r_x$ with
$r_x \in [1,\beta^{\wfpmu}) \cap \wz{}$  and
$y = \wsign{y} \beta^{\wfpemin} r_y$ with
$r_y \in [1,\beta^{\wfpmu}) \cap \wz{}$.
If $\wsign{x} = - \wsign{y}$ then $\wabs{x + y} = \beta^{\wfpemin} \wabs{r_x - r_y}$
and $\wabs{x + y}$ is either $0$ or subnormal, because
\[
\wabs{r_x - r_y} < \max \wset{r_x,r_y} < \beta^{\wfpmu}.
\]
If $\wsign{x} = \wsign{y}$ then $\wabs{x + y} = \beta^{\wfpemin} \wlr{r_x + r_y}$
with $1  < r_x + r_y < 2 \beta^{\wfpmu} \leq \beta^{1 + \mu}$. If $r_x + r_y < \beta^{\wfpmu}$ then
$x + y$ is subnormal, otherwise $\wabs{x + y} \geq \beta^{\wfpemin + \mu} = \wfpnu$
and Prop. \ref{propNu} implies that $\wabs{x + y} \in \wfpsi$.
\peProof{Proposition}{propSubnormalSum}\\

%%%%%%%%%%%%%%%%%%%%%%%%%%%%%%%%%%%%%%%%%%%%%%%%%%%%%%%%%%%%%%%%
% End of Proof of prop propSubnormalSum
%%%%%%%%%%%%%%%%%%%%%%%%%%%%%%%%%%%%%%%%%%%%%%%%%%%%%%%%%%%%%%%%
%
%

%
%
%%%%%%%%%%%%%%%%%%%%%%%%%%%%%%%%%%%%%%%%%%%%%%%%%%%%%%%%%%%%%%%%%%%%%%%%%%%%%%%%%%%
% Proof of prop. propCriticalSum
%%%%%%%%%%%%%%%%%%%%%%%%%%%%%%%%%%%%%%%%%%%%%%%%%%%%%%%%%%%%%%%%%%%%%%%%%%%%%%%%%%%

\pbProofB{Proposition}{propCriticalSum}
Let us start with $s := x + z > 0$.
If $x \leq \beta^{e + \mu}$ then Prop. \ref{propCriticalSum} holds
because $z = s - x \geq \beta^{e} \wlr{r + 1/2} \geq \beta^{e} / 2$.
If $x \geq \beta^{e + \mu + 1}$ then
\[
z := s - x \leq \beta^e \wlr{r + 1/2 + \beta^{\mu} - \beta^{\mu + 1}}
= - \beta^{e} / 2 - \wlr{ \wlr{\beta - 1} \beta^{\mu} - \wlr{r + 1}}
\leq -\beta^{e}/2,
\]
because $r \in [0,\wlr{\beta - 1} \beta^{\mu}) \cap \wz{}$,
and again $\wabs{z} \geq \beta^e/2$.
Therefore, we only need to analyze the case $\beta^{e + \mu} < x < \beta^{e + \mu + 1}$.
In this case, by Prop. \ref{propNormalForm}, $x = \beta^{d} \wlr{\beta^\mu + t}$ for
$d \in \wz{}$ and $t \in [0,\wlr{\beta - 1}\beta^{\mu}) \cap \wz{}$.
Prop. \ref{propOrder} implies that $d = e$ and
\[
z = s - x = \beta^{e} \wlr{r - t + 1/2} \ \ \Rightarrow \ \
\wabs{z} \geq \beta^{e} \wabs{r - t + 1/2} \geq \beta^{e}/2,
\]
because $r - t \in \wz{}$, and we are done with the case $x + z > 0$.
Finally, when $x + z  < 0$ the argument
above for $-x$ and $-z$ and leads to $\wabs{z} = \wabs{-z} \geq \beta^e/2$.
\peProof{Proposition}{propCriticalSum}\\

%%%%%%%%%%%%%%%%%%%%%%%%%%%%%%%%%%%%%%%%%%%%%%%%%%%%%%%%%%%%%%%%
% End of Proof of prop propCriticalSum
%%%%%%%%%%%%%%%%%%%%%%%%%%%%%%%%%%%%%%%%%%%%%%%%%%%%%%%%%%%%%%%%
%
%

%
%
%%%%%%%%%%%%%%%%%%%%%%%%%%%%%%%%%%%%%%%%%%%%%%%%%%%%%%%%%%%%%%%%%%%%%%%%%%%%%%%%%%%
% proof of propAlpha
%%%%%%%%%%%%%%%%%%%%%%%%%%%%%%%%%%%%%%%%%%%%%%%%%%%%%%%%%%%%%%%%%%%%%%%%%%%%%%%%%%%

\pbProofB{Proposition}{propAlpha}
If $\wfpf$ is perfect then $\alpha = 0 \in \wfpf$ and Prop. \ref{propAlpha} is trivial.
If $\wfpf$ is the MPFR system
$\wfpsm_{\wfpemin,\beta,\mu}$ then $\alpha = \beta^{\wfpemin + \mu} \in \wfpbin_{\wfpemin} \subset \wfpf$
and if $x \in \wfpf \setminus \wset{0}$ then $\wabs{x} \in \wfpbin_{e}$ for some $e \geq \wfpemin$.
By definition of $\wfpbin_{e}$, $\wabs{x} = \beta^{e} \wlr{\beta^{\mu} + r}$
with $r \geq 0$ and $\wabs{x} \geq \alpha$.
Finally, if $\wfpf$ is the IEEE system $\wfpsi_{\wfpemin,\beta,\mu}$ then
$\alpha = \beta^{\wfpemin} \in \wfpbin_{\wfpemin} \subset \wfpf$ and if $x \in \wfpf \setminus \wset{0}$ then
$\wabs{x} \in \wfpf \setminus \wset{0}$ and
either (i) $\wabs{x} \in \wfps_{\wfpemin}$ or (ii) $\wabs{x} \in \wfpbin_{e}$ with $e \geq \wfpemin$.
In case (i), $\wabs{x} = \beta^{\wfpemin} r$ for $r \in \wz{} \setminus {0}$ and $\wabs{x} \geq \alpha$.
As for the MPFR system, in case (ii) $\wabs{x} \geq \alpha$.
\peProof{Proposition}{propAlpha}\\

%%%%%%%%%%%%%%%%%%%%%%%%%%%%%%%%%%%%%%%%%%%%%%%%%%%%%%%%%%%%%%%%%%%%%%%%%%%%%%%%%%%
% End of proof of propAlpha
%%%%%%%%%%%%%%%%%%%%%%%%%%%%%%%%%%%%%%%%%%%%%%%%%%%%%%%%%%%%%%%%%%%%%%%%%%%%%%%%%%%
%
%

%
%
%%%%%%%%%%%%%%%%%%%%%%%%%%%%%%%%%%%%%%%%%%%%%%%%%%%%%%%%%%%%%%%%%%%%%%%%%%%%%%%%%%%
% Proof of prop. propEmptyNormalRange
%%%%%%%%%%%%%%%%%%%%%%%%%%%%%%%%%%%%%%%%%%%%%%%%%%%%%%%%%%%%%%%%%%%%%%%%%%%%%%%%%%%

\pbProofB{Proposition}{propEmptyNormalRange}
We show that if $z \in \wlr{\beta^{e}\wlr{\beta^{\mu} + r}, \ \beta^{e}\wlr{\beta^{\mu} + r + 1}}$
then $z \not \in \wfpf$.
By Prop. \ref{propOrder} and \ref{propNormalForm},
there exists $w$ with $r < w < r + 1$ such that $z = \beta^{e}\wlr{\beta^{\wfpmu} + w}$.
$z \not \in - \wfpbin_{d}$ because $z > 0$.
Prop. \ref{propOrder} implies that if $d > e$ then $z < y$ for
$y \in \wfpbin_{d}$ and if $d < e$ then $z > y$ for $y \in \wfpbin_{d}$.
Therefore, $z \not \in \bigcup_{d \neq e} \wfpbin_{d}$.
Moreover, if $y \in \wfpbin_{e}$ then
$y = \beta^{e} \wlr{\beta^{\wfpmu} + s}$ with $s \in \wz{}$ and
$y \neq z$ because $w \not \in \wz{}$. As a result,
$z \not \in \bigcup_{d \in \wz{}} \wlr{\wfpbin_{d} \cup -\wfpbin_{d}}$.
This proves that $z \not \in \wfpf$ when
$\wfpf$ is a perfect or MPFR system. Finally,
if $\wfpf$ is an IEEE system $\wfpsi_{\wfpemin}$ then
$e \geq \wfpemin$ because $e$ is an exponent for $\wfpf$,
and $z > y$ for all $y \in \wfps_{\wfpemin} \cup -\wfps_{\wfpemin}$.
This shows that $z \not \in \wfps_{\wfpemin} \cup -\wfps_{\wfpemin}$ and
$z \not \in \wfpf$.
\peProof{Proposition}{propEmptyNormalRange}\\

%%%%%%%%%%%%%%%%%%%%%%%%%%%%%%%%%%%%%%%%%%%%%%%%%%%%%%%%%%%%%%%%%%%%%%%%%%%%%%%%%%%
% End of proof of prop EmptyNormalRange
%%%%%%%%%%%%%%%%%%%%%%%%%%%%%%%%%%%%%%%%%%%%%%%%%%%%%%%%%%%%%%%%%%%%%%%%%%%%%%%%%%%
%
%

%
%
%%%%%%%%%%%%%%%%%%%%%%%%%%%%%%%%%%%%%%%%%%%%%%%%%%%%%%%%%%%%%%%%%%%%%%%%%%%%%%%%%%%
% Proof of prop. propEmptySubnormalRange
%%%%%%%%%%%%%%%%%%%%%%%%%%%%%%%%%%%%%%%%%%%%%%%%%%%%%%%%%%%%%%%%%%%%%%%%%%%%%%%%%%%

\pbProofB{Proposition}{propEmptySubnormalRange}
We show that if $z \in \wlr{\beta^{\wfpemin}r, \ \beta^{\wfpemin}\wlr{r + 1}}$
then $z \not \in \wfpsi$.
We have that $\wabs{z} < \beta^{\wfpemin} \max \wset{\wabs{r},\wabs{r + 1}} \leq \beta^{\wfpemin + \mu}$
and $\wabs{z} \not \in \bigcup_{e = \wfpemin}^{\infty} \wfpbin_{e}$.
Moreover, $w := \beta^{-\wfpemin} z$
is such that $r < w < r + 1$
and $z = \beta^{\wfpemin} w$. It follows that $w \not \in \wz{}$
and $\wabs{z} \not \in \wfps_{\wfpemin}$,
and combining the arguments above and symmetry (Prop. \ref{propSymmetry}) we conclude that $z \not \in \wfpsi$.
\peProof{Proposition}{propEmptySubnormalRange}\\

%%%%%%%%%%%%%%%%%%%%%%%%%%%%%%%%%%%%%%%%%%%%%%%%%%%%%%%%%%%%%%%%%%%%%%%%%%%%%%%%%%%
% End of proof of prop EmptySubNormalRange
%%%%%%%%%%%%%%%%%%%%%%%%%%%%%%%%%%%%%%%%%%%%%%%%%%%%%%%%%%%%%%%%%%%%%%%%%%%%%%%%%%%
%
%

%
%
%%%%%%%%%%%%%%%%%%%%%%%%%%%%%%%%%%%%%%%%%%%%%%%%%%%%%%%%%%%%%%%%%%%%%%%%%%%%%%%%%%%
% Proof of prop propSysScale
%%%%%%%%%%%%%%%%%%%%%%%%%%%%%%%%%%%%%%%%%%%%%%%%%%%%%%%%%%%%%%%%%%%%%%%%%%%%%%%%%%%

\pbProofB{Proposition}{propSysScale}
Since $\wcal{A}_e := \wset{\beta x, x \in \wfpbin_{e}} = \wfpbin_{e + 1}$,
the set $\wfpsc$ in Definition \ref{longDefPerfect} is such that
$x \in \wfpsc$ if and only if $\beta x \in \wfpsc$.
For the MPFR system $\wfpsm_{\wfpemin,\beta,\mu}$, if $x \in \wfpsm$ then
$\wabs{x} \in \wfpbin_{e}$ for some $e \geq \wfpemin$, $\wabs{\beta x} \in \wfpbin_{e +1} \subset \wfpsm$
and $\beta x \in \wfpsm$ by symmetry (Prop. \ref{propSymmetry}.)
For the IEEE system $\wfpsi_{\wfpemin,\beta,\mu}$, if $x \in \wfpsi$ then
either $x \in \wfpbin_{e}$ for some $e \geq \wfpemin$,
and the argument used in the MFPR case applies to $x$, or
$x = \wsign{x} \beta^{\wfpemin} r$ with $r \in [0,\beta^{\mu}) \cap \wz{}$.
If $\beta r < \mu$ then $\wabs{\beta x} = \beta^{\wfpemin} \wlr{\beta r} \in \wfps_{\wfpemin}$
and $\beta x \in \wfpsi$ by symmetry. If $\beta r \geq \beta^\mu$ then
$s = \beta r - \beta^{\mu} \in [0, \wlr{\beta - 1} \beta^{\mu}) \cap \wz{}$
and $\wabs{\beta x} = \beta^{\wfpemin} \wlr{\beta^{\mu} + s} \in \wfpbin_{\wfpemin} \subset \wfpsi$
and $s \in \wfpsi$ by symmetry.
\peProof{Proposition}{propSysScale}\\

%%%%%%%%%%%%%%%%%%%%%%%%%%%%%%%%%%%%%%%%%%%%%%%%%%%%%%%%%%%%%%%%%%%%%%%%%%%%%%%%%%%
% Proof of prop propSysScale
%%%%%%%%%%%%%%%%%%%%%%%%%%%%%%%%%%%%%%%%%%%%%%%%%%%%%%%%%%%%%%%%%%%%%%%%%%%%%%%%%%%
%
%

\subsection{Rounding}
\label{secRounding}
This section proves the propositions about rounding to nearest stated
previously, and states and proves more propositions about rounding.

\subsubsection{Propositions}
\label{secRoundingProps}
In this section we state more propositions regarding rounding to nearest.\\

%
%
%%%%%%%%%%%%%%%%%%%%%%%%%%%%%%%%%%%%%%%%%%%%%%%%%%%%%%%%%%%%%%%%%%%%%%%%%%%%%%%%%%%
% propRoundSign
%%%%%%%%%%%%%%%%%%%%%%%%%%%%%%%%%%%%%%%%%%%%%%%%%%%%%%%%%%%%%%%%%%%%%%%%%%%%%%%%%%%

\pbPropBT{propRoundSign}{Propagation of the sign}
If $\wfl{z} \neq 0$ then $\wsign{\wfl{z}} = \wsign{z}$. For
a general $z \in \wrone{}$, $\wfl{z} = \wsign{z} \wabs{\wfl{z}}$.
\peFullProp{propRoundSign}

%
%
%%%%%%%%%%%%%%%%%%%%%%%%%%%%%%%%%%%%%%%%%%%%%%%%%%%%%%%%%%%%%%%%%%%%%%%%%%%%%%%%%%%
% propScale
%%%%%%%%%%%%%%%%%%%%%%%%%%%%%%%%%%%%%%%%%%%%%%%%%%%%%%%%%%%%%%%%%%%%%%%%%%%%%%%%%%%
%
\pbPropBT{propRoundScale}{Rounding after scaling}
Let $m$ be an integer. If $\wfpf$ is perfect then the function
$\wflr{s}{z} := \beta^{-m} \wfl{\beta^m z}$ rounds to nearest in $\wfpf$.
\peProp{propRoundScale}

%
%
%%%%%%%%%%%%%%%%%%%%%%%%%%%%%%%%%%%%%%%%%%%%%%%%%%%%%%%%%%%%%%%%%%%%%%%%%%%%%%%%%%%
% prop propRoundInterval
%%%%%%%%%%%%%%%%%%%%%%%%%%%%%%%%%%%%%%%%%%%%%%%%%%%%%%%%%%%%%%%%%%%%%%%%%%%%%%%%%%%%

\pbPropBT{propRoundOnInterval}{Rounding in an interval}
If $a,b \in \wfpf$ and $a \leq z \leq b$ then $\wfl{z} \in [a,b]$
and $\wabs{\wfl{z}- z} \leq (b - a)/2$. Moreover, if $z < m := (a + b)/2$ then
$\wfl{z} < b$ and if $z > m$ then $\wfl{z} > a$.
\peFullProp{propRoundOnInterval}

%
%
%%%%%%%%%%%%%%%%%%%%%%%%%%%%%%%%%%%%%%%%%%%%%%%%%%%%%%%%%%%%%%%%%%%%%%%%%%%%%%%%%%%
% propRoundCombination
%%%%%%%%%%%%%%%%%%%%%%%%%%%%%%%%%%%%%%%%%%%%%%%%%%%%%%%%%%%%%%%%%%%%%%%%%%%%%%%%%%%

\pbPropBT{propRoundCombination}{Combination}
For $\wcal{A}_1, \wcal{A}_2 \subset \wrone{}$ with
$\wcal{A}_1 \cup \wcal{A}_2 = \wrone{}$,
let $f_i: \wcal{A}_i \rightarrow \wrone{}$
be such that, for $z_i \in \wcal{A}_i$ and $x \in \wfpf{}$,
$\wfc{f_i}{z_i} \in \wfpf$ and $\wabs{z_i - \wfc{f_i}{z_i}} \leq \wabs{z_i - x}$. The function
   $\wflm: \wrone{} \rightarrow \wrone{}$ given by
   $\wfl{z} = \wfc{f_1}{z}$  for $z \in \wcal{A}_1$ and
   $\wfl{z} = \wfc{f_2}{z}$  for $z \in \wcal{A}_2 \setminus \wcal{A}_1$ rounds
   to nearest in $\wfpf$.
\peProp{propRoundCombination}

%
%
%%%%%%%%%%%%%%%%%%%%%%%%%%%%%%%%%%%%%%%%%%%%%%%%%%%%%%%%%%%%%%%%%%%%%%%%%%%%%%%%%%%
% propRoundExtension
%%%%%%%%%%%%%%%%%%%%%%%%%%%%%%%%%%%%%%%%%%%%%%%%%%%%%%%%%%%%%%%%%%%%%%%%%%%%%%%%%%%

\pbPropBT{propRoundExtension}{Extension}
If $\wcal{A} \subset \wrone{}$ and $f: \wcal{A} \rightarrow \wrone{}$ is such that,
for $z \in \wcal{A}$ and $x \in \wfpf$, $\wfc{f}{z} \in \wfpf$ and
$\wabs{z - \wfc{f}{z}} \leq \wabs{z - x}$
then there exists a function $\wflm$ which rounds to
nearest in $\wfpf$ and is such that $\wfl{z} = \wfc{f}{z}$ for $z \in \wcal{A}$.
\peProp{propRoundExtension}

\subsubsection{Proofs}
\label{secRoundingProofs}
In this section we prove the propositions regarding rounding to nearest.\\
%
%
%%%%%%%%%%%%%%%%%%%%%%%%%%%%%%%%%%%%%%%%%%%%%%%%%%%%%%%%%%%%%%%%%%%%%%%%%%%%%%%%%%%
% Proof of Prop. propRoundIdentity
%%%%%%%%%%%%%%%%%%%%%%%%%%%%%%%%%%%%%%%%%%%%%%%%%%%%%%%%%%%%%%%%%%%%%%%%%%%%%%%%%%%

\pbProofB{Proposition}{propRoundIdentity}
By definition of rounding to nearest,
$0 = \wabs{x - x} \geq \wabs{\wfl{x} - x}$. Therefore, $\wfl{x} = x$.
\peProof{Proposition}{propRoundIdentity}\\

%%%%%%%%%%%%%%%%%%%%%%%%%%%%%%%%%%%%%%%%%%%%%%%%%%%%%%%%%%%%%%%%%%%%%%%%%%%%%%%%%%%
% Proof of Prop. propRoundIdentity
%%%%%%%%%%%%%%%%%%%%%%%%%%%%%%%%%%%%%%%%%%%%%%%%%%%%%%%%%%%%%%%%%%%%%%%%%%%%%%%%%%%
%
%

%
%
%%%%%%%%%%%%%%%%%%%%%%%%%%%%%%%%%%%%%%%%%%%%%%%%%%%%%%%%%%%%%%%%%%%%%%%%%%%%%%%%%%%
% Proof of prop propRoundMonotone
%%%%%%%%%%%%%%%%%%%%%%%%%%%%%%%%%%%%%%%%%%%%%%%%%%%%%%%%%%%%%%%%%%%%%%%%%%%%%%%%%%%

\pbProofB{Proposition}{propRoundMonotone}
Let us show that if $\wfl{z} > \wfl{w}$ then $z > w$.
Indeed, in this case we have that
\[
\wabs{\wfl{w} - z} \geq \wabs{\wfl{z} - z} \geq \wfl{z} - z > \wfl{w} - z.
\]
Therefore, $\wabs{\wfl{w} - z} > \wfl{w} - z$ and this implies that $z > \wfl{w}$.
It follows that
\[
z - \wfl{w} = \wabs{\wfl{w} - z} \geq \wabs{\wfl{z} - z} \geq \wfl{z} - z \ \  \Rightarrow \ \
z \geq \frac{\wfl{z} + \wfl{w}}{2}.
\]
Similarly,
\[
\wabs{w - \wfl{z}} \geq \wabs{w - \wfl{w}} \geq w - \wfl{w} > w - \wfl{z}
\ \ \Rightarrow \ \ w   \leq \wfl{z},
\]
and
\[
\wfl{z} - w = \wabs{\wfl{z} - w} \geq \wabs{\wfl{w} - w} \geq w - \wfl{w} \ \ \Rightarrow \ \
w \leq \frac{\wfl{z} + \wfl{w}}{2}.
\]
As a result, $w \leq \wlr{\wfl{z} + \wfl{w}}/2 \leq z$. Moreover, $w \neq z$ because
$\wfl{z} \neq \wfl{w}$. Therefore, $z > w$ as we have claimed. Logically, we have
proved that $z \leq w \Rightarrow \wfl{z} \leq \wfl{w}$.

When $x \in \wfpf$ we have that $\wfl{x} = x$ (Prop. \ref{propRoundIdentity})
and the argument above shows that $\wfl{z} > x \Rightarrow z > x$ and
$x > \wfl{z} \Rightarrow x > z$. Moreover,
\[
\wabs{x} > \wabs{\wfl{z}}
\Rightarrow
\wabs{x} > \wfl{z}
\Rightarrow
\wabs{x} > z
\]
and using the function $\wrm{m}$ in Prop. \ref{propRoundMinus} we obtain
\[
\wabs{x} > \wabs{\wfl{z}}
\Rightarrow
\wabs{x} > -\wfl{z}
\Rightarrow
\wabs{x} > \wflr{m}{-z}
\Rightarrow
\wabs{x} > -z.
\]
Therefore, $\wabs{x} > \wabs{\wfl{z}} \Rightarrow \wabs{x} > \max \wset{z,-z} = \wabs{z}$.

Finally, if $\wabs{x} < \wabs{\wfl{z}}$ then either
(i) $\wfl{z} < 0$ or (ii) $\wfl{z} > 0$. In both cases
Prop. \ref{propRoundSign} shows that $\wsign{z} = \wsign{\wfl{z}}$.
In case (i) $z$ is positive and
\[
\wabs{x} < \wabs{\wfl{z}}
\Rightarrow
\wabs{x} < \wfl{z}
\Rightarrow
\wabs{x} < z = \wabs{z}.
\]
and in case (ii) $z$ is negative and
\[
\wabs{x} < \wabs{\wfl{z}}
\Rightarrow
\wabs{x} < -\wfl{z}
\Rightarrow
\wabs{x} < \wflr{m}{-z}
\Rightarrow
\wabs{x} < -z = \wabs{z}.
\]
Therefore, $\wabs{x} < \wabs{\wfl{z}} \Rightarrow \wabs{x} < \wabs{z}$ in both
cases and we are done.
\peProof{Proposition}{propRoundMonotone}\\

%%%%%%%%%%%%%%%%%%%%%%%%%%%%%%%%%%%%%%%%%%%%%%%%%%%%%%%%%%%%%%%%%%%%%%%%%%%%%%%%%%%
% Proof of prop propRoundMonotone
%%%%%%%%%%%%%%%%%%%%%%%%%%%%%%%%%%%%%%%%%%%%%%%%%%%%%%%%%%%%%%%%%%%%%%%%%%%%%%%%%%%
%
%

%
%
%%%%%%%%%%%%%%%%%%%%%%%%%%%%%%%%%%%%%%%%%%%%%%%%%%%%%%%%%%%%%%%%%%%%%%%%%%%%%%%%%%%
% Proof of Prop. propRoundMinus
%%%%%%%%%%%%%%%%%%%%%%%%%%%%%%%%%%%%%%%%%%%%%%%%%%%%%%%%%%%%%%%%%%%%%%%%%%%%%%%%%%%

\pbProofB{Proposition}{propRoundMinus}
If $x \in \wfpf$ and $z \in \wrone{}$ then $-x \in \wfpf$ by symmetry and
\[
\wabs{\wflr{m}{z} - z} = \wabs{ \wlr{-\wfl{-z}} - z} =
\wabs{\wfl{-z} - \wlr{-z}} \leq \wabs{\wlr{-x} - \wlr{-z}} = \wabs{x - z}.
\]
Therefore, $\wabs{\wflr{m}{z} - z} \leq \wabs{x -z}$ and $\wrm{m}$ rounds to nearest.
\peProof{Proposition}{propRoundMinus}\\

%%%%%%%%%%%%%%%%%%%%%%%%%%%%%%%%%%%%%%%%%%%%%%%%%%%%%%%%%%%%%%%%%%%%%%%%%%%%%%%%%%%
% End of Proof of Prop. propRoundMinus
%%%%%%%%%%%%%%%%%%%%%%%%%%%%%%%%%%%%%%%%%%%%%%%%%%%%%%%%%%%%%%%%%%%%%%%%%%%%%%%%%%%
%
%

%
%
%%%%%%%%%%%%%%%%%%%%%%%%%%%%%%%%%%%%%%%%%%%%%%%%%%%%%%%%%%%%%%%%%%%%%%%%%%%%%%%%%%%
% Proof of prop propRoundNormal
%%%%%%%%%%%%%%%%%%%%%%%%%%%%%%%%%%%%%%%%%%%%%%%%%%%%%%%%%%%%%%%%%%%%%%%%%%%%%%%%%%%%

\pbProofB{Proposition}{propRoundNormal}
Let us start with $z > 0$.
$w \leq \wlr{\beta - 1} \beta^{\mu}$ and if $\wfloor{w} = \wlr{\beta - 1} \beta^{\mu}$ then
$a = b = \beta^{e + 1 + \wfpmu} \in \wfpbin_{e + 1}$, and this
implies that $a,b \in \wfpf$ because $e + 1$ is also an exponent for $\wfpf$.
Similarly, if $\wceil{w} = \wlr{\beta - 1} \beta^{\mu}$ then $b \in \wfpbin_{e+1} \subset \wfpf$.
If $\wceil{w} < \wlr{\beta - 1} \beta^{\mu}$ then $0 \leq \wfloor{w} \leq \wceil{w} < \wlr{\beta - 1}\beta^{\mu}$ and
$a,b \in \wfpbin_{e} \subset \wfpf$. Therefore, in all cases, $a,b \in \wfpf$.
If $w \in \wz{}$ then $\wfloor{w} = \wceil{w}$ and $z = a = b \in \wfpf$ and
$\wfl{z} = a = b = m$ because $\wfl{x} = x$ when $x \in \wfpf$ by Prop. \ref{propRoundIdentity}.
If $w \not \in \wz{}$ then $\wceil{w} = \wfloor{w} + 1$,
Prop. \ref{propEmptyNormalRange} shows that $(a,b) \cap \wfpf = \emptyset$,
Equation \pRef{thRoundNormal} follows from Prop. \ref{propRoundOnInterval},
and we also have that $(b - a) / 2 \leq \beta^{e}/2$.

For the last paragraph in Prop. \ref{propRoundNormal},
we either have (i) $r \leq w$ or (ii) $r > w$.
In case (i)
\[
r \leq w \leq r + \wabs{w - r} < r + 1/2 \Rightarrow \wfloor{w} = r,  \ \ \wceil{w} = r + 1
\]
and
\[
\frac{1}{2} \wlr{ \wfloor{w} + \wceil{w} } = r + 1/2 > w.
\]
This implies that $a = \beta^{e} \wlr{\beta^{\mu} + r}$,
$b = \beta^{e} \wlr{\beta^{\mu} + r + 1}$ and $z < \wlr{a + b}/2$,
and the results in the previous paragraph show that
$\beta^{e} \wlr{\beta^{\mu} + r} = a = \wfl{z}$.

In case (ii), $r > w \geq 0 \Rightarrow r \geq 1$ and
\[
r - 1/2  \leq w < r  \Rightarrow \wfloor{w} = r - 1,  \ \ \wceil{w} = r
\hspace{0.2cm} \wrm{and} \hspace{0.2cm}
\frac{1}{2} \wlr{ \wfloor{w} + \wceil{w} } = r - 1/2 < w.
\]
This implies that $a = \beta^{e} \wlr{\beta^{\mu} + r - 1}$,
$b = \beta^{e} \wlr{\beta^u + r}$ and $z > \wlr{a + b}/2$,
and the results in the first paragraph of this proof show that
$\beta^{e} \wlr{\beta^u + r} = b = \wfl{z}$.

Finally, for $z < 0$ the arguments above for $\tilde{z} = -z$ and $\wflmx$ equal to the function
$\wrm{m}$ in Prop. \ref{propRoundMinus} and symmetry (Prop. \ref{propSymmetry}) prove
Prop. \ref{propRoundNormal} for $z$.
\peProof{Proposition}{propRoundNormal}\\

%%%%%%%%%%%%%%%%%%%%%%%%%%%%%%%%%%%%%%%%%%%%%%%%%%%%%%%%%%%%%%%%%%%%%%%%%%%%%%%%%%%
% End of Proof of prop propRoundNormal
%%%%%%%%%%%%%%%%%%%%%%%%%%%%%%%%%%%%%%%%%%%%%%%%%%%%%%%%%%%%%%%%%%%%%%%%%%%%%%%%%%%%
%
%

%
%
%%%%%%%%%%%%%%%%%%%%%%%%%%%%%%%%%%%%%%%%%%%%%%%%%%%%%%%%%%%%%%%%%%%%%%%%%%%%%%%%%%%
% Proof of Prop. propRoundSubnormal
%%%%%%%%%%%%%%%%%%%%%%%%%%%%%%%%%%%%%%%%%%%%%%%%%%%%%%%%%%%%%%%%%%%%%%%%%%%%%%%%%%%%

\pbProofB{Proposition}{propRoundSubnormal}
Recall that $\nu = \beta^{\wfpemin+ \mu} \in \wfpbin_{\wfpemin} \subset \wfpsi$,
and by symmetry $-\nu \in \wfpsi$. Let us write $w := \beta^{-\wfpemin} z$ and
$r :=\wfloor{w}$. We have that $a = \beta^{\wfpemin} r$
and if $r = w$ then $a = b = z$ and $\wfl{z} = z$ by Prop. \ref{propRoundIdentity}
and Prop. \ref{propRoundSubnormal} is valid. Let us then assume that $r \neq w$.
This implies that $w \not \in \wz{}$, $r < \beta^{\mu}$, $r + 1 = \wceil{w}$
and $b = \beta^{\wfpemin} \wlr{r + 1}$.
We have that $a \in \wfpf$ because
\begin{eqnarray}
\nonumber
w < 1  - \beta^{\mu} & \Rightarrow & r = -\beta^{\mu} \Rightarrow a = -\nu \in - \wfpbin_{\wfpemin} \subset \wfpf, \\
\nonumber
1 - \beta^{\mu} < w < 0 &  \Rightarrow  & 1 - \beta^{\mu} \leq r < 0  \Rightarrow  a \in -\wfps_{\wfpemin} \subset \wfpf, \\
\nonumber
0 < w < 1 &  \Rightarrow &  r = 0 \Rightarrow a = 0 \in \wfpf, \\
\nonumber
1 < w < \beta^{\mu}  & \Rightarrow &  1 \leq r < \beta^{\mu} \Rightarrow  a \in \wfps_{\wfpemin} \subset \wfpf,
\end{eqnarray}
and $b \in \wfpf$ because
\begin{eqnarray}
\nonumber
-\beta^{\mu} < w < -1 &  \Rightarrow &  1 - \beta^{\mu} < r + 1 \leq -1 \Rightarrow b \in -\wfps_{\wfpemin} \subset \wfpf, \\
\nonumber
-1 < w  < 0 &  \Rightarrow & r + 1 = 0 \Rightarrow b = 0 \in \wfpf, \\
\nonumber
0 < w < \beta^{\mu } - 1  &  \Rightarrow &  1 \leq r + 1 < \beta^{\mu} \Rightarrow   b \in \wfps_{\wfpemin} \subset \wfpf, \\
\nonumber
\beta^{\mu} - 1 < w < \beta^{\mu} &  \Rightarrow &  r + 1 = \beta^{\mu} \Rightarrow   b = \nu \in \wfpbin_{\wfpemin} \in \wfpf.
\end{eqnarray}
Therefore, by monotonicity $\wfl{z} \in [a,b] \cap \wfpf$ and Prop. \ref{propEmptySubnormalRange}
implies that $\wfl{z} \in \wset{a,b}$. It follows that
\[
\wabs{\wfl{z} - z} = \min \wset{z - a, b - z} \leq
\frac{b - a}{2} = \frac{\beta^{\wfpemin} \wlr{r + 1} - \beta^{\wfpemin} \wlr{r}}{2} = \wfpa /2.
\]
Finally, if $z < m$ then $\wabs{b - z} > \wabs{a - z} \Rightarrow \wfl{z} = a$ and
if $z > m$ then $\wabs{a - z} > \wabs{b - z} \Rightarrow \wfl{z} = b$.

\peProof{Proposition}{propRoundSubnormal}\\

%%%%%%%%%%%%%%%%%%%%%%%%%%%%%%%%%%%%%%%%%%%%%%%%%%%%%%%%%%%%%%%%%%%%%%%%%%%%%%%%%%%
% End of proof of orop propRoundSubnormal
%%%%%%%%%%%%%%%%%%%%%%%%%%%%%%%%%%%%%%%%%%%%%%%%%%%%%%%%%%%%%%%%%%%%%%%%%%%%%%%%%%%
%
%

%
%
%%%%%%%%%%%%%%%%%%%%%%%%%%%%%%%%%%%%%%%%%%%%%%%%%%%%%%%%%%%%%%%%%%%%%%%%%%%%%%%%%%%
% Proof of Prop. propRoundSubnormal
%%%%%%%%%%%%%%%%%%%%%%%%%%%%%%%%%%%%%%%%%%%%%%%%%%%%%%%%%%%%%%%%%%%%%%%%%%%%%%%%%%%%

\pbProofB{Proposition}{propRoundBelowAlpha}
Note that, by Prop. \ref{propAlpha}, if
$x \in \wfpf \setminus \wset{0,\pm \wfpa}$ then $\wabs{x} > \wfpa$.
When $\wabs{z} < \alpha /2$, if $x \in \wfpf \setminus \wset{0}$ then
Prop. \ref{propAlpha} implies that $\wabs{x} \geq \wfpa$ and
\[
\wabs{x - z} \geq \wabs{x} - \wabs{z} \geq \wfpa - \wabs{z} > \wfpa / 2 > \wabs{z - 0},
\]
and $\wfl{z} = 0$ because $0 \in \wfpf$.

When $\wabs{z} = \wfpa/2$,
$\wabs{z - 0} = \wabs{z - \wsign{z} \wfpa} = \wfpa / 2$
and $\wabs{z - \wlr{-\wsign{z}}} = 3 \alpha/2$. As a result,
if $x  \in \wfpf \setminus \wset{0,\pm \alpha}$ then
\[
\wabs{x - z} \geq \wabs{x} - \wabs{z} > \wfpa - \wfpa / 2 = \wfpa/2 = \wabs{z - 0},
\]
and the bounds above imply that $\wfl{z} \in \wset{0, \wsign{z} \wfpa}$.

When $\wfpa/2 < \wabs{z} < \wfpa$, $\wabs{z - \wsign{z} \wfpa} = \wfpa - \wabs{z} < \wfpa/2$,
$\wabs{z - 0} = \wabs{z} > \wfpa/2$ and
\[
\wabs{z - \wlr{-\wsign{z}} \wfpa} = \wabs{z} + \wfpa > \wfpa / 2.
\]
Moreover, if $x \in \wfpf \setminus  \wset{0,\pm \wfpa}$ has the same sign as $z$ then
$x > \wfpa$ and
\[
\wabs{x - z} = x - z = \wlr{x - \wfpa} + \wlr{\wfpa - z} > \wabs{\wsign{z} \wfpa - z}.
\]
and if $x$ has the opposite sign of $z$ then
$\wabs{x - z} \geq \wabs{x} > \wfpa > \wabs{\wsign{z} \wfpa - z}$,
and the bounds in this paragraph imply that $\wfl{z} = \wsign{z} \wfpa$.
Finally, if $\wabs{z} = \wfpa$ then $\wfl{z} = z = \wsign{z} \wfpa$ by Prop.
\ref{propRoundIdentity}.
\peProof{Proposition}{propRoundBelowAlpha}\\

%%%%%%%%%%%%%%%%%%%%%%%%%%%%%%%%%%%%%%%%%%%%%%%%%%%%%%%%%%%%%%%%%%%%%%%%%%%%%%%%%%%
% End of proof of orop propRoundSubnormal
%%%%%%%%%%%%%%%%%%%%%%%%%%%%%%%%%%%%%%%%%%%%%%%%%%%%%%%%%%%%%%%%%%%%%%%%%%%%%%%%%%%
%
%

%
%
%%%%%%%%%%%%%%%%%%%%%%%%%%%%%%%%%%%%%%%%%%%%%%%%%%%%%%%%%%%%%%%
% Proof of proposition propRoundAdapt
%%%%%%%%%%%%%%%%%%%%%%%%%%%%%%%%%%%%%%%%%%%%%%%%%%%%%%%%%%%%%%%

\pbProofB{Proposition}{propRoundAdapt}
Let $\wcal{A}$ be the set
$\wset{z \in \wrone{} \ \wrm{with} \ \wabs{z} \geq \nu_{\wfpf}}$
and $f: \wcal{A} \rightarrow \wrone{}$ the
function $\wfc{f}{z} = \wfl{z}$.
We claim that if $x \in \wfpsc$ and $z \in \wcal{A}$ then
$\wabs{z - \wfc{f}{z}} \leq \wabs{x - \wfc{f}{z}}$.
In fact, if $x \in \wfpf$ then $\wabs{x - z} \geq \wabs{\wfl{z} - z} = \wabs{\wfc{f}{z} - z}$,
because $\wflm$ rounds to nearest in $\wfpf$.
If $x \not \in \wfpf$ then
\[
x \in \wfpsc \setminus \wfpf \subset
\wlr{\bigcup_{e = -\infty}^{+\infty} \wlr{ \wfpbin_{e} \bigcup -\wfpbin_e} \setminus
\bigcup_{e = \wfpemin}^{+\infty} \wlr{\wfpbin_{e} \bigcup - \wfpbin_e}}
= \bigcup_{e = -\infty}^{\wfpemin- 1} \wlr{\wfpbin_{e} \bigcup - \wfpbin_e}
\]
and
\[
\wabs{x} < \beta^{\wfpemin - 1} \wlr{\beta^{\mu} + \wlr{\beta - 1} \beta^{\mu}}
= \beta^{\wfpemin + \mu - 1} = \nu_{\wfpf}.
\]
since $\nu_{\wfpf} \in \wfpf$, if $z \geq \nu_{\wfpf}$ then
$z \geq \wabs{x} \geq x$ and
\[
\wabs{x - z} = z - x = \wabs{\nu_{{\wfpf}} - z} + \wabs{\nu_{\wfpf} - x}
\geq \wabs{\wfl{z} - z} +  \wabs{\nu_{\wfpf} - x} > \wabs{\wfc{f}{z} - z}.
\]
Similarly, $\nu_{\wfpf} \in \wfpf$ and if $z \leq - \nu_{\wfpf}$ then
$z \leq - \wabs{x} \leq x$ and
\[
\wabs{x - z} = x - z = \wabs{-\nu_{{\wfpf}} - z} + \wabs{-\nu_{\wfpf} - x}
\geq \wabs{\wfl{z} - z} +  \wabs{-\nu_{\wfpf} - x} > \wabs{\wfc{f}{z} - z},
\]
Therefore, $\wabs{x - z} \geq \wabs{\wfc{f}{z} - z}$ in all cases.
To complete the proof it suffices to take the extension of $f$
to $\wrone{}$ given by Prop. \ref{propRoundExtension}.
\peProof{Proposition}{propRoundCompletion}\\

%%%%%%%%%%%%%%%%%%%%%%%%%%%%%%%%%%%%%%%%%%%%%%%%%%%%%%%%%%%%%%%%%%%%%%%%%%%%%%%%%%%
% End of proof of propRoundCompletion
%%%%%%%%%%%%%%%%%%%%%%%%%%%%%%%%%%%%%%%%%%%%%%%%%%%%%%%%%%%%%%%%%%%%%%%%%%%%%%%%%%%
%
%

%
%
%%%%%%%%%%%%%%%%%%%%%%%%%%%%%%%%%%%%%%%%%%%%%%%%%%%%%%%%%%%%%%%%%%%%%%%%%%%%%%%%%%%
% Proof of propRoundIEEEX
%%%%%%%%%%%%%%%%%%%%%%%%%%%%%%%%%%%%%%%%%%%%%%%%%%%%%%%%%%%%%%%%%%%%%%%%%%%%%%%%%%%

\pbProofB{Proposition}{propRoundIEEEX}
For $k = 1,\dots,n$ let $\wflmx_{k}$ be the adapter
of $\wflm_k$ in Prop. \ref{propRoundAdapt}.
On the one hand, by the definition of $\wflmx_k$, we have that
if $x,y \in \wfpf$ and $\wabs{x + y} \geq \nu_{\wfpsi}$ then
\pbDef{priee}
\wflk{k}{x + y } = \wflxkf{k}{x + y}.
\peDef{priee}
On the other hand, Lemma \ref{lemIEEESum} shows that
Equation \pRef{priee} holds when $\wabs{x + y} \leq \nu_{\wfpsi}$.
Therefore, Equation \pRef{priee} holds for all $x,y \in \wfpsi$.
For $\wvec{x} \in \wfpsi^{n+1}$ define $\wvec{z} \in \wrn{n}$
as $z_1 := x_0 + x_1$ and $z_k := x_k$ for $2 \leq k < n$.
We now prove by induction
that $\wfpsumkf{k}{\wvec{z},\wflmt} = \wfpsumkf{k}{\wvec{z},\wflmtx}$.
By definition,  $\wfpsumkf{0}{\wvec{z},\wflmt} = 0 = \wfpsumkf{0}{\wvec{z},\wflmtx} $.
Let us then analyze $k > 0$ assuming that
$\wfpsumkf{k - 1}{\wvec{z},\wflmtx} = \wfpsumkf{k - 1}{\wvec{z},\wflmt} \in \wfpsi$.
Using Equation \pRef{priee} we deduce that
\[
\wfpsumkf{k}{\wvec{z},\wflmtx} =
  \wflxkf{k}{\wfpsumkf{k-1}{\wvec{z},\wflmt} + z_{k}}
= \wflk{k}{\wfpsumkf{k-1}{\wvec{x},\wflmt} + z_{k}} =
\wfpsumkf{k}{\wvec{z},\wflmt} \in \wfpsi
\]
and we are done.
\peProof{Proposition}{propRoundIEEEX}\\

%%%%%%%%%%%%%%%%%%%%%%%%%%%%%%%%%%%%%%%%%%%%%%%%%%%%%%%%%%%%%%%%%%%%%%%%%%%%%%%%%%%
% Proof of prop. RoundIEEEX
%%%%%%%%%%%%%%%%%%%%%%%%%%%%%%%%%%%%%%%%%%%%%%%%%%%%%%%%%%%%%%%%%%%%%%%%%%%%%%%%%%%
%
%

%
%
%%%%%%%%%%%%%%%%%%%%%%%%%%%%%%%%%%%%%%%%%%%%%%%%%%%%%%%%%%%%%%%%%%%%%%%%%%%%%%%%%%%
% Proof of propRoundMPFRX
%%%%%%%%%%%%%%%%%%%%%%%%%%%%%%%%%%%%%%%%%%%%%%%%%%%%%%%%%%%%%%%%%%%%%%%%%%%%%%%%%%%

\pbProofB{Proposition}{propRoundMPFRX}
For $k = 1,\dots,n$, let $\wflmxk{k}$ be the adapter
of $\wflmk{k}$ in Prop. \ref{propRoundAdapt}.
By the definition of $\wflmxk{k}$ we have that
if $x,y \in \wfpf$ and $\wabs{x + y} \geq \wfpa_{\wfpsi} = \nu_{\wfpsi}$ then
\pbDef{prmp}
\wflk{k}{x + y } = \wflxkf{k}{x + y},
\peDef{prmp}
and, of course, this equation is also satisfied when $x + y = 0$.
For $\wvec{x} \in \wfpsi^{n+1}$ define $\wvec{z} \in \wrn{n}$
as $z_1 := x_0 + x_1$ and $z_k := x_k$ for $2 \leq k < n$.
We now prove by induction that if $y_k := \wfpsumkf{k - 1}{\wvec{z},\wflmt} + z_k \geq 0$
for $k = 0,\dots,n$ then
$\wfpsumkf{k}{\wvec{z},\wflmt} = \wfpsumkf{k}{\wvec{z},\wflmtx}$.
By definition,  $\wfpsumkf{0}{\wvec{z},\wflmt} = 0 = \wfpsumkf{0}{\wvec{z},\wflmtx} $.
Let us then analyze $k > 0$ assuming that
$\wfpsumkf{k - 1}{\wvec{z},\wflmtx} = \wfpsumkf{k - 1}{\wvec{z},\wflmt} \in \wfpsm$.
The assumption that $y_k \geq 0$ and Prop. \ref{propAlpha} implies
that either $y_k = 0$ or $y_k \geq \wfpa_{\wfpsm} = \wfpnu_{\wfpsm}$, and in both
cases Equation \pRef{prmp} holds for $x + y = y_k$.
It follows that
\[
\wfpsumkf{k}{\wvec{z},\wflmtx} =
\wflxkf{k}{\wfpsumkf{k-1}{\wvec{z},\wflmt} + z_{k}}
= \wflxkf{k}{y_k}
= \wflk{k}{y_k} =
\wfpsumkf{k}{\wvec{z},\wflmt} \in \wfpsm,
\]
and we are done.
\peProof{Proposition}{propRoundMPFRX}\\

%%%%%%%%%%%%%%%%%%%%%%%%%%%%%%%%%%%%%%%%%%%%%%%%%%%%%%%%%%%%%%%%%%%%%%%%%%%%%%%%%%%
% Proof of prop. RoundMPFRX
%%%%%%%%%%%%%%%%%%%%%%%%%%%%%%%%%%%%%%%%%%%%%%%%%%%%%%%%%%%%%%%%%%%%%%%%%%%%%%%%%%%
%
%

%
%
%%%%%%%%%%%%%%%%%%%%%%%%%%%%%%%%%%%%%%%%%%%%%%%%%%%%%%%%%%%%%%%
% Proof of proposition propFlat
%%%%%%%%%%%%%%%%%%%%%%%%%%%%%%%%%%%%%%%%%%%%%%%%%%%%%%%%%%%%%%%

\pbProofB{Proposition}{propFlat}
If $w = 0$ then we can take $\delta = \beta^{e-1}/2$,
because in this case $z \in \wfpbin_{e} \subset \wfpf$ and $\wflk{1}{z} = z$
by Prop. \ref{propRoundIdentity} and,
according to Prop. \ref{propNormalFormDis},
if $\wabs{y - z} < \delta$ then either
\begin{itemize}
\item[(i)] $y = \wsign{z} \beta^{e} \wlr{\beta^{\mu} + v}$
with
\[
0 \leq v = \beta^{-e} \wabs{y - z} <  \beta^{-e} \delta < 1/2 \Rightarrow \wfloor{v} = 0
\]
and $\wflk{2}{y} = \wflk{1}{z} = z$ by Prop. \ref{propRoundNormal}, or
\item[(ii)] $y = \wsign{z} \beta^{e - 1} \wlr{\beta^{\mu} + v}$
for
\[
\wlr{\beta - 1} \beta^{\mu} - \beta^{1 - e} \wabs{y - z} = v <
\wlr{\beta - 1} \beta^{\mu} \Rightarrow
\]
\[
\wlr{\beta - 1} \beta^{\mu} - 1/2 < v <  \wlr{\beta - 1} \beta^{\mu}
\Rightarrow \wceil{v} =  \wlr{\beta - 1} \beta^{\mu}
\]
and, by Prop. \ref{propRoundNormal},
\[
\wflk{2}{y} = \wsign{z} \beta^{e - 1} \wlr{\beta^{\mu} +  \wlr{\beta - 1} \beta^{\mu}} =
\wsign{z} \beta^{e + \mu} = z = \wflk{1}{z}
\]
\end{itemize}
Let us then assume that $w > 0$ and write $m := \wfloor{w} + 1/2$
and show that
\[
\delta = \beta^{e} \min \wset{w, \, \wlr{\beta - 1} \beta^{\mu} - w, \, 1/2 - \wabs{m - w}, \, \wabs{m - w}}
\]
is a valid choice. Note that $\delta > 0$, because
$\wabs{m - w} \leq 1/2$ for a general $w$ and $w \neq 1/2$ for
the particular $w$ we discuss here.
If $\wabs{y - z} < \delta$ then Prop. \ref{propNormalFormCont} implies that
$y = \wsign{z} \beta^{e} \wlr{\beta^{\mu} + v}$ with
\[
\wabs{v - w} = \beta^{-e} \wabs{y - z} <  \beta^{-e} \delta \leq \min \wset{ 1/2 - \wabs{m - w}, \wabs{m - w}}.
\]
On the one hand, if $w < m$ then $\wabs{m - w} = m - w$,
\[
\wfloor{w} = m - 1/2 < m - \wlr{\wabs{w - m} + \wabs{v - w}}  \leq v
\leq w + \wabs{w - v} < w + \wabs{m - w} =  m,
\]
$\wfloor{v} = \wfloor{w}$
and Prop. \ref{propRoundNormal} implies that
$\wflk{2}{y} = \wflk{1}{z} = \wsign{z} \beta^{e} \wlr{\beta^{\mu} + \wfloor{w}}$.
On the other hand, if $w > m$ then $\wabs{m - w} = w - m$,
\[
m = w - \wabs{w - m} < w - \wabs{w - v}
\leq
v
\leq
m + \wlr{\wabs{w - m} + \wabs{v - w}} < m + 1/2 =  \wceil{w},
\]
$\wceil{v} = \wceil{w}$
and Prop. \ref{propRoundNormal} implies that
$\wflk{2}{y} = \wflk{1}{z} = \wsign{z} \beta^{e} \wlr{\beta^{\mu} + \wceil{w}}$.
\peProof{Proposition}{propFlat}\\

%%%%%%%%%%%%%%%%%%%%%%%%%%%%%%%%%%%%%%%%%%%%%%%%%%%%%%%%%%%%%%%%%%%%%%
% End of proof of prop propFlat
%%%%%%%%%%%%%%%%%%%%%%%%%%%%%%%%%%%%%%%%%%%%%%%%%%%%%%%%%%%%%%%%%%%%%%
%
%

%
%
%%%%%%%%%%%%%%%%%%%%%%%%%%%%%%%%%%%%%%%%%%%%%%%%%%%%%%%%%%%%%%%%%%%%%%%%%%%%%%%%%%%
% propSumScale
%%%%%%%%%%%%%%%%%%%%%%%%%%%%%%%%%%%%%%%%%%%%%%%%%%%%%%%%%%%%%%%%%%%%%%%%%%%%%%%%%%%

\pbProofB{Proposition}{propSumScale}
For $k = 1,\dots,n$ Props. \ref{propRoundMinus} and \ref{propRoundScale}
show that the function
$\wflxkf{k}{z} := \sigma \beta^{-m} \wflk{k}{\sigma \beta^m z}$
rounds to nearest in $\wfpsc$, and we define $\wflmtx := \wset{\wflmxk{1},\dots,\wflmxk{n}}$.
We now prove by induction in $k = 0, \dots, n$ that
\pbDef{sumSPFoo}
\wfpsumkf{k}{\sigma \beta^m \wvec{z},\wflmt} = \sigma \beta^m \wfpsumkf{k}{\wvec{z},\wflmtx},
\peDef{sumSPFoo}
For $k = 0$,
$\wfpsumkf{0}{\sigma \beta^m \wvec{z},\wflmt} = 0 = \sigma \beta^m \wfpsumkf{k}{\wvec{z},\wflmtx}$
by definition. Assuming that \pRef{sumSPFoo} holds for $k \geq 0$ we have that
\[
\wfpsumkf{k+1}{\sigma \beta^m \wvec{z},\wflmt} =
\wflk{k+1}{\wfpsumkf{k}{\sigma \beta^m \wvec{z},\wflmt} + \sigma \beta^m  z_k}
\]
\[
= \wflk{k+1}{\sigma \beta^m \wlr{\wfpsumkf{k}{\wvec{z},\wflmtx} + z_k}}
= \sigma \beta^m \wflxkf{k+1}{\wfpsumkf{k}{\wvec{z},\wflmtx} + z_k}
=
\wfpsumkf{k+1}{\wvec{z},\wflmtx},
\]
and we are done.
\peProof{Proposition}{propSumScale}\\

%%%%%%%%%%%%%%%%%%%%%%%%%%%%%%%%%%%%%%%%%%%%%%%%%%%%%%%%%%%%%%%%%%%%%%
% End of proof of prop propSumScale
%%%%%%%%%%%%%%%%%%%%%%%%%%%%%%%%%%%%%%%%%%%%%%%%%%%%%%%%%%%%%%%%%%%%%%
%
%

%
%
%%%%%%%%%%%%%%%%%%%%%%%%%%%%%%%%%%%%%%%%%%%%%%%%%%%%%%%%%%%%%%%%%%%%%%%%%%%%%%%%%%%
% propRoundSign
%%%%%%%%%%%%%%%%%%%%%%%%%%%%%%%%%%%%%%%%%%%%%%%%%%%%%%%%%%%%%%%%%%%%%%%%%%%%%%%%%%%

\pbProofB{Proposition}{propRoundSign}
Prop. \ref{propRoundIdentity} shows that $\wfl{0} = 0$. Therefore,
if $\wfl{z} \neq 0$ then either (i) $z > 0$ or (ii) $z < 0$. In case (i)
\[
\wabs{\wfl{z} - z} \leq \wabs{0 - z} \Rightarrow
z - \wfl{z} \leq z \Rightarrow \wfl{z} \geq 0 \Rightarrow \wsign{\wfl{z}} = 1 = \wsign{z}.
\]
In case (ii)
\[
\wabs{\wfl{z} - z} \leq \wabs{0 - z} \Rightarrow
\wfl{z} - z \leq -z \Rightarrow \wfl{z} \leq 0.
\]
Since $\wfl{z} \neq 0$ this implies that $\wsign{\wfl{z}} = -1 = \wsign{z}$.
It follows that if
$\wfl{z} \neq 0$ then $\wfl{z} = \wsign{\wfl{z}} \wabs{\wfl{z}} = \wsign{z} \wabs{\wfl{z}}$
and it is clear that this equality also holds when $\wfl{z} = 0$.
\peProof{Proposition}{propRoundSign}\\

%%%%%%%%%%%%%%%%%%%%%%%%%%%%%%%%%%%%%%%%%%%%%%%%%%%%%%%%%%%%%%%%%%%%%%%%%%%%%%%%%%%
% End of Proof of Prop. propRoundSign
%%%%%%%%%%%%%%%%%%%%%%%%%%%%%%%%%%%%%%%%%%%%%%%%%%%%%%%%%%%%%%%%%%%%%%%%%%%%%%%%%%%
%
%

%
%
%%%%%%%%%%%%%%%%%%%%%%%%%%%%%%%%%%%%%%%%%%%%%%%%%%%%%%%%%%%%%%%
% Proof of prop propRoundScale
%%%%%%%%%%%%%%%%%%%%%%%%%%%%%%%%%%%%%%%%%%%%%%%%%%%%%%%%%%%%%%%

\pbProofB{Proposition}{propRoundScale}
Suppose $x \in \wfpf{}$ and $z \in \wrone{}$.
When $\wfpf$ is perfect we have that
$\beta^{m} x \in \wfpf$ by Prop. \ref{propSysScale}
and since $\wflm$ rounds to nearest we have
\[
\wabs{\wflr{s}{z} - z} \, = \,  \wabs{ \wlr{\beta^{-m} \wfl{\beta^m z}} - z}
 \,  =  \,    \beta^{-m} \wabs{\wfl{\beta^m z} - \wlr{\beta^{m} z}}
 \]
\[
\leq \,  \beta^{-m} \wabs{\wfl{\beta^m z} - \wlr{\beta^{m} x}}
   \,  = \,         \wabs{ \wlr{\beta^{-m} \wfl{\beta^m z}} - x}
   \,  = \,  \wabs{\wflr{s}{z} - x}.
\]
Therefore, $\wrm{s}$ rounds to nearest in $\wfpf$.
\peProof{Proposition}{propRoundScale}\\

%%%%%%%%%%%%%%%%%%%%%%%%%%%%%%%%%%%%%%%%%%%%%%%%%%%%%%%%%%%%%%%
% End of Proof of proposition propRoundScale
%%%%%%%%%%%%%%%%%%%%%%%%%%%%%%%%%%%%%%%%%%%%%%%%%%%%%%%%%%%%%%%
%
%

%
%
%%%%%%%%%%%%%%%%%%%%%%%%%%%%%%%%%%%%%%%%%%%%%%%%%%%%%%%%%%%%%%%%%%%%%%%%%%%%%%%%%%%
% Proof of prop propRoundOnInterval
%%%%%%%%%%%%%%%%%%%%%%%%%%%%%%%%%%%%%%%%%%%%%%%%%%%%%%%%%%%%%%%%%%%%%%%%%%%%%%%%%%%%

\pbProofB{Proposition}{propRoundOnInterval}
Since $x = a,b \in \wfpf$, the definition of rounding to nearest yields
$\wabs{z - a} \geq \wabs{z - \wfl{z}}$ and $\wabs{z - b} \geq \wabs{z - \wfl{z}}$.
If $y < a$ then $y < z$ and
\[
\wabs{z - y} = z - y > z - a = \wabs{z - a} \geq \wabs{z - \wfl{z}}
\ \ \Rightarrow  \ \ \wabs{z - y} > \wabs{z - \wfl{z}}
\]
Therefore, $\wfl{z} \neq y$. Similarly,
if $y > b$ then $y > z$ and
\[
\wabs{z - y} = y - z > b - z = \wabs{z - b} \geq \wabs{z - \wfl{z}}
\ \ \Rightarrow  \ \ \wabs{z - y} > \wabs{z - \wfl{z}}
\]
As a result, $\wfl{z} \neq y$ and
$\wfl{z} \in \wrone{} \setminus \wlr{\wset{y < a} \cup \wset{y > b}} =  [a,b]$.
If $z \leq m$ then
\[
\wabs{\wfl{z} - z} \leq \wabs{a - z} = z - a \leq m - a = \delta := \wlr{b - a}/2.
\]
and if $z \geq m$ then
\[
\wabs{\wfl{z} - z} \leq \wabs{b - z} = b - z \leq b - m = \delta.
\]
Therefore, $\wabs{\wfl{z} - z} \leq \delta$.
If $ z < m$ then
\[
\wfl{z} \leq \wabs{\wfl{z} - z} + \wlr{z - a} + a \leq \delta + z < \delta + m = b,
\]
and $\wfl{z} < b$. If $ z > m$ then
\[
\wfl{z} \geq b - \wlr{b - z} - \wabs{z - \wfl{z}} \geq z - \delta > m - \delta = a,
\]
and $\wfl{z} > a$.
\peProof{Proposition}{propRoundOnInterval}\\

%%%%%%%%%%%%%%%%%%%%%%%%%%%%%%%%%%%%%%%%%%%%%%%%%%%%%%%%%%%%%%%%%%%%%%%%%%%%%%%%%%%
% End of Proof of prop propRoundOnInterval
%%%%%%%%%%%%%%%%%%%%%%%%%%%%%%%%%%%%%%%%%%%%%%%%%%%%%%%%%%%%%%%%%%%%%%%%%%%%%%%%%%%%
%
%

%
%
%%%%%%%%%%%%%%%%%%%%%%%%%%%%%%%%%%%%%%%%%%%%%%%%%%%%%%%%%%%%%%%
% Proof of proposition propRoundCombination
%%%%%%%%%%%%%%%%%%%%%%%%%%%%%%%%%%%%%%%%%%%%%%%%%%%%%%%%%%%%%%%

\pbProofB{Proposition}{propRoundCombination}
If $z \in \wrone{}$ then either (i) $z \in \wcal{A}_1$ or
(ii) $z \in \wcal{A}_2 \setminus \wcal{A}_1$. In case (i), for $x \in \wfpf$ we
have that $\wabs{x - z} \geq \wabs{\wflk{1}{z} - z}$ by hypothesis.
Therefore, $\wabs{x - z} \geq \wabs{\wflk{1}{z} - z} = \wabs{\wfl{z} - z}$ in case (i).
In case (ii), for $x \in \wfpf$ we
have that $\wabs{x - z} \geq \wabs{\wflk{2}{z} - z} = \wabs{\wfl{z} - z}$.
As a result, $\wabs{x - z} \geq \wabs{\wfl{z} - z}$ in both cases and $\wflm$ rounds to nearest in $\wfpf$.
\peProof{Proposition}{propRoundCombination}\\

%%%%%%%%%%%%%%%%%%%%%%%%%%%%%%%%%%%%%%%%%%%%%%%%%%%%%%%%%%%%%%%%%%%%%%%%%%%%%%%%%%%
% propRoundCombination
%%%%%%%%%%%%%%%%%%%%%%%%%%%%%%%%%%%%%%%%%%%%%%%%%%%%%%%%%%%%%%%%%%%%%%%%%%%%%%%%%%%
%
%

%
%
%%%%%%%%%%%%%%%%%%%%%%%%%%%%%%%%%%%%%%%%%%%%%%%%%%%%%%%%%%%%%%%%%%%%%%%%%%%%%%%%%%%
% proof of propRoundExtension
%%%%%%%%%%%%%%%%%%%%%%%%%%%%%%%%%%%%%%%%%%%%%%%%%%%%%%%%%%%%%%%%%%%%%%%%%%%%%%%%%%%

\pbProofB{Proposition}{propRoundExtension}
We assume that there exists
 $f_2: \wrone{} \rightarrow \wrone{}$ which
rounds to nearest in $\wfpf$.
Take $\wcal{A}_1 = \wcal{A}$ and $\wcal{A}_2 = \wrone \setminus \wcal{A}$.
Prop. \ref{propRoundCombination} with $f_1 = f$
implies that there exists $\wflm$ which rounds to nearest in $\wfpf$ and is such that
$\wfl{z} = \wfc{f}{z}$ for $z \in \wcal{A}$.
\peProof{Proposition}{propRoundExtension} \\

%%%%%%%%%%%%%%%%%%%%%%%%%%%%%%%%%%%%%%%%%%%%%%%%%%%%%%%%%%%%%%%%%%%%%%%%%%%%%%%%%%%
% End of Proof of prop propRoundExtension
%%%%%%%%%%%%%%%%%%%%%%%%%%%%%%%%%%%%%%%%%%%%%%%%%%%%%%%%%%%%%%%%%%%%%%%%%%%%%%%%%%%%
%
%

\subsection{Tightness}
\label{secTight}
In this section we prove the propositions regarding tightness,
and present and prove additional propositions about this subject.

\subsubsection{Propositions}
\label{secTightProps}
In this section we present additional propositions regarding tightness.

%
%
%%%%%%%%%%%%%%%%%%%%%%%%%%%%%%%%%%%%%%%%%%%%%%%%%%%%%%%%%%%%%%%%%%%%%%%%%%%%%%%%%%%
% Proof TightContinuous
%%%%%%%%%%%%%%%%%%%%%%%%%%%%%%%%%%%%%%%%%%%%%%%%%%%%%%%%%%%%%%%%%%%%%%%%%%%%%%%%%%%

\pbPropBT{propTightContinuous}{Tightness and continuity}
Let $\wcal{A}$, $\wcal{B}$ and $\wcal{C}$
be topological spaces and $\wcal{R}$ a set.
If $g: \wcal{A} \times \wcal{B}  \rightarrow \wcal{C}$ is continuous
and $h: \wcal{A} \times \wcal{R} \rightarrow \wcal{B}$ is tight
then $f: \wcal{A} \times \wcal{R} \rightarrow \wcal{C}$
given by $\wfc{f}{a, r} = \wfc{g}{a,\wfc{h}{a,r}}$ is tight.
In particular, if $\wcal{R}$ is a tight set of functions
from $\wcal{A}$ to $\wcal{B}$ then the function
$f: \wcal{A} \times \wcal{R} \rightarrow \wcal{B}$ given by
$\wfc{f}{a,r} = \wfc{g}{a,\wfc{r}{a}}$ is tight.
\peProp{propTightContinuous}

%
%
%%%%%%%%%%%%%%%%%%%%%%%%%%%%%%%%%%%%%%%%%%%%%%%%%%%%%%%%%%%%%%%%%%%%%%%%%%%%%%%%%%%
% Proof Tight chain
%%%%%%%%%%%%%%%%%%%%%%%%%%%%%%%%%%%%%%%%%%%%%%%%%%%%%%%%%%%%%%%%%%%%%%%%%%%%%%%%%%%

\pbPropBT{propTightChain}{Tight chain rule}
Let $\wcal{A}$, $\wcal{B}$ and $\wcal{C}$ be topological spaces and
let $\wcal{T}$ and $\wcal{U}$ be sets. If the functions
$h: \wcal{A} \times \wcal{T} \rightarrow \wcal{B}$ and
$g: \wcal{B} \times \wcal{U} \rightarrow \wcal{C}$
are tight then the function
$f: \wcal{A} \times \wlr{\wcal{T} \times \wcal{U}} \rightarrow \wcal{C}$ given by
$\wfc{f}{a, \wlr{t,u}} := \wfc{g}{\wfc{h}{a,t}, u}$
is tight.
\peProp{propTightChain}

\subsubsection{Proofs}
\label{secTightProofs}
This section contains the proofs of the propositions regarding tightness.\\

%
%
%%%%%%%%%%%%%%%%%%%%%%%%%%%%%%%%%%%%%%%%%%%%%%%%%%%%%%%%%%%%%%%%%%%%%%
% Proof of prop propWholeTight
%%%%%%%%%%%%%%%%%%%%%%%%%%%%%%%%%%%%%%%%%%%%%%%%%%%%%%%%%%%%%%%%%%%%%%

\pbProofB{Proposition}{propWholeIsTight}
Let $\wcal{R}$ be the set of all functions
which round to nearest in $\wfpf$ and let
$\wcal{S} = \wset{\wlr{z_k, \wflm_k}, k \in \wn{}} \subset \wrone{} \times \wcal{R}$
be a sequence with $\lim_{k \rightarrow \infty} z_k = z$.
Props. \ref{propRoundSubnormal}, \ref{propRoundBelowAlpha} and \ref{propFlat}
imply that there exist
$a,b \in \wfpf$ and
$\delta > 0$ such that if $\wabs{y - z} < \delta$ then
$\wfl{y} \in \wset{a,b}$ for $\wflm \in \wcal{R}$.
Let  $m \in \wn{}$ be such that $k > m \Rightarrow \wabs{z_k - z} < \delta$
and define
$\wcal{A} := \wset{k \geq m \ \wrm{with} \ \wflk{k}{z_k} = a}$
and
$\wcal{B} := \wset{k \geq m \ \wrm{with} \ \wflk{k}{z_k} = b}$.
Since $\wcal{A} \bigcup \wcal{B} = \wset{k \geq m, k \in \wn{}}$
is infinite, $\wcal{A}$ or $\wcal{B}$ is infinite.
By exchanging $a$ and $b$ if necessary, we may assume that $\wcal{A}$
is infinite, and $\wset{\wlr{z_{n_k}, \wflm_{n_k}}, n_k \in \wcal{A}}$
is a subsequence of $\wcal{S}$.
We claim that the function $\wflm: \wrone{} \rightarrow \wrone{}$ given
by $\wfl{w} = \wflk{m}{w}$ for $w \neq z$ and $\wfl{z} = a$ rounds to
nearest in $\wfpf$. Indeed, if $z' \in \wfpf \setminus \wset{z}$ and $w \in \wrone{}$ then
\[
\wabs{w - \wfl{z'}} = \wabs{w - \wflk{m}{z'}} \geq \wabs{z' - \wflk{m}{z'}} = \wabs{z' - \wfl{z'}}
\]
because $\wflmk{m}$ rounds to nearest in $\wfpf$, and
\[
\wabs{w - \wfl{z}} = \wabs{w - a} =
\wabs{w - \wflk{n_k}{z_k}}
\geq \wabs{z_k - \wflk{n_k}{z_k}}
=  \wabs{z_k - a}
=  \wabs{z_k - \wfl{z}}.
\]
because the $\wflmk{n_k}$ round to nearest in $\wfpf$. Taking the limit $k \rightarrow \infty$
in the equation above we obtain $\wabs{w - \wfl{z}} \geq \wabs{z - \wfl{z}}$,
and $\wflm$ rounds to nearest in $\wfpf$. Finally,
\[
\lim_{k \rightarrow \infty} \wfc{\varphi}{z_{n_k},\wflmk{n_k}} =
 \lim_{k \rightarrow \infty} \wflk{n_k}{z_{n_k}} =
 a = \wfl{z}
 \]
 and $\wcal{R}$ is tight.
\peProof{Proposition}{propWholeIsTight}\\

%%%%%%%%%%%%%%%%%%%%%%%%%%%%%%%%%%%%%%%%%%%%%%%%%%%%%%%%%%%%%%%%%%%%%%
% Proof of prop propWholeIsTight
%%%%%%%%%%%%%%%%%%%%%%%%%%%%%%%%%%%%%%%%%%%%%%%%%%%%%%%%%%%%%%%%%%%%%%
%
%

%
%
%%%%%%%%%%%%%%%%%%%%%%%%%%%%%%%%%%%%%%%%%%%%%%%%%%%%%%%%%%%%%%%%%%%%%%
% Proof of prop propSumsAreTight
%%%%%%%%%%%%%%%%%%%%%%%%%%%%%%%%%%%%%%%%%%%%%%%%%%%%%%%%%%%%%%%%%%%%%%

\pbProofB{Proposition}{propSumsAreTight}
For $n = 0$,
$\wfc{T_0}{\wvec{z},\wflmt} = 0$ and Prop. \ref{propSumsAreTight} follows
from Prop. \ref{propTightContinuous},
because constant functions are continuous.
Assuming that Prop. \ref{propSumsAreTight} holds for $n \geq 0$,
let us show that it holds for $n + 1$. By induction
and Prop. \ref{propTightContinuous} the function
$h: \wrn{n+1} \times \wcal{R}^n \rightarrow \wrn{n+1} \times \wrone{}$
given  by $\wfc{h}{\wvec{w},\wflmt} := \wlr{\wfc{T_n}{\wrm{P}_n \wvec{w}, \wflmt},w_{n+1}}$
is tight. The function $g: \wlr{\wrn{n+1} \times \wrone} \times \wcal{R} \rightarrow \wrn{n+2}$
given by $\wfc{g}{\wlr{\wvec{w},z}, \wflm} := \wlr{\wvec{w}, \wfl{w_{n+1} + z}}$
is also tight by Prop. \ref{propTightContinuous} because
$\wcal{R}$ is tight.
Finally, Prop. \ref{propSumsAreTight} follows
from Prop. \ref{propTightChain} for $f = T_n$, $g$ and $h$ because
$\wfc{T_{n+1}}{\wvec{z},\wflmt} = \wfc{g}{\wfc{h}{\wvec{w},\wrm{P}_n \wflmt},\wflmk{n+1}}$.
\peProof{Proposition}{propSumsAreTight}\\

%%%%%%%%%%%%%%%%%%%%%%%%%%%%%%%%%%%%%%%%%%%%%%%%%%%%%%%%%%%%%%%%%%%%%%
% End of proof of prop propSumsAreTight
%%%%%%%%%%%%%%%%%%%%%%%%%%%%%%%%%%%%%%%%%%%%%%%%%%%%%%%%%%%%%%%%%%%%%%
%
%

%
%
%%%%%%%%%%%%%%%%%%%%%%%%%%%%%%%%%%%%%%%%%%%%%%%%%%%%%%%%%%%%%%%%%%%%%%%%%%%%%%%%%%%
% Proof TightContinuous
%%%%%%%%%%%%%%%%%%%%%%%%%%%%%%%%%%%%%%%%%%%%%%%%%%%%%%%%%%%%%%%%%%%%%%%%%%%%%%%%%%%

\pbProofB{Proposition}{propTightContinuous}
Let $\wset{\wlr{a_k, r_k},  k \in \wn{}} \subset \wcal{A} \times \wcal{R}$ be a sequence with
$\lim_{k \rightarrow \infty}a_k = a$. Since $h$ is tight,
there exists $r \in \wcal{R}$ and a subsequence
$\wset{\wlr{a_{n_k}, r_{n_k}}, k \in \wn{}}$ such that
$\lim_{k \rightarrow \infty} \wfc{h}{a_{n_k}, r_{n_k}} = \wfc{h}{a,r}$.
By continuity of $g$,
\[
\lim_{k \rightarrow \infty} \wfc{f}{a_{n_k},r_{n_k}} =
\lim_{k \rightarrow \infty} \wfc{g}{a_{n_k}, \wfc{h}{a_{n_k},r_{n_k}}} = \wfc{g}{a,\wfc{h}{a,r}} = \wfc{f}{a,r},
\]
and $f$ is tight. To handle the particular case, note that when
$\wcal{R}$ is set of tight functions as in the hypothesis the function
$h: \wcal{A} \times \wcal{R} \rightarrow \wcal{B}$ given by
$\wfc{h}{a,r} = \wfc{r}{a}$ is tight.
\peProof{Proposition}{propTightContinuous}\\

%
%
%%%%%%%%%%%%%%%%%%%%%%%%%%%%%%%%%%%%%%%%%%%%%%%%%%%%%%%%%%%%%%%%%%%%%%
% Proof of prop propTightChain
%%%%%%%%%%%%%%%%%%%%%%%%%%%%%%%%%%%%%%%%%%%%%%%%%%%%%%%%%%%%%%%%%%%%%%

\pbProofB{Proposition}{propTightChain}
Let $\wset{\wlr{a_k,\wlr{t_k,u_k}}, k \in \wn{}}
\subset \wcal{A} \times \wlr{\wcal{T} \times \wcal{U}}$ be a sequence such that
$\lim_{k \rightarrow \infty} a_k = a$.
Since $h$ is tight, there exists $t \in \wcal{T}$
and a subsequence $\wset{\wlr{a_{n_k}, t_{n_k}}, k \in \wn{}}$ of
$\wset{\wlr{a_k,t_k}, k \in \wn{}}$ such that
$b_{n_k} := \wfc{h}{a_{n_k}, t_{n_k}}$ satisfies
$\lim_{k \rightarrow \infty} b_{n_k} = \wfc{h}{a,t} =: b$.
Since $b_{n_k}$ converges to $b$ and $g$ is tight, there exists
$u \in \wcal{U}$ and
a subsequence $\wset{\wlr{b_{m_k},u_{m_k}}, k \in \wn{}}$ of
$\wset{\wlr{b_{n_k},u_{n_k}}, k \in \wn{}}$
such that
\[
\wfc{g}{b,u} =
\lim_{k \rightarrow \infty} \wfc{g}{b_{m_k},u_{m_k}}
= \lim_{k \rightarrow \infty} \wfc{g}{\wfc{h}{a_{m_k}, t_{m_k}},u_{m_k}}.
\]
This leads to
\[
\wfc{f}{a,\, \wlr{t,u} \, } =
\wfc{g}{\wfc{h}{a,t}, \, u} =
\wfc{g}{b, \, u} =
\]
\[
\lim_{k \rightarrow \infty} \wfc{g}{\wfc{h}{a_{m_k}, \, t_{m_k}}, \, u_{m_k}}
=
\lim_{k \rightarrow \infty} \wfc{f}{a_{m_k}, \, \wlr{t_{m_k}, u_{m_k}}},
\]
and $f$ is tight.
\peProof{Proposition}{propTightChain}\\

%%%%%%%%%%%%%%%%%%%%%%%%%%%%%%%%%%%%%%%%%%%%%%%%%%%%%%%%%%%%%%%%%%%%%%
% End of proof of prop propTightChain
%%%%%%%%%%%%%%%%%%%%%%%%%%%%%%%%%%%%%%%%%%%%%%%%%%%%%%%%%%%%%%%%%%%%%%
%
%

\subsection{Examples}
\label{secExample}
In this section we verify the examples 2 to 5.
Example 1 needs no verification.\\

%
%
%%%%%%%%%%%%%%%%%%%%%%%%%%%%%%%%%%%%%%%%%%%%%%%%%%%%%%%%%%%%%%%%%%%%%%%%%
% Verification of example exQuadGrowth
%%%%%%%%%%%%%%%%%%%%%%%%%%%%%%%%%%%%%%%%%%%%%%%%%%%%%%%%%%%%%%%%%%%%%%%%%

\pbVerifyB{exQuadGrowth} Our parcels are
$y_0 := 1$ and $y_k := 1 + 2^{\wfloor{\wfc{\log_2}{k+1}}} u$ for $k = 1, \dots, n := 2^m - 1$
and we break ties downward. If $1 \leq 2^{\ell} - 1 \leq k < 2^{\ell + 1} - 1$
then $y_k = 1 + 2^{\ell} u$ and we now show by induction that, for $k \geq 1$,
\pbDef{exqgI}
\sum_{i = 0}^k y_i = k + 1 + \frac{4^{\ell} + 2}{3} u  + \wlr{k + 1 - 2^{\ell}} 2^{\ell} u
\hspace{1cm} \wrm{and} \hspace{1cm}
\wfl{\sum_{i = 0}^k y_i} = k + 1.
\peDef{exqgI}
Indeed, for $k = 1$ we have $\ell = 1$ and $y_0 + y_1 = 2 + 2 u$, the first
equality in Equation \pRef{exqgI} is clearly correct and the second
holds because we break ties downward.

If $\wfl{\sum_{i = 0}^k y_i} = k + 1$ and  $2^{\ell} - 2 \leq k < 2^{\ell + 1} - 2$
then $y_{k+1} = 1 + 2^{\ell} u$ and
\[
\wfl{\sum_{i = 0}^{k+1} y_i} = \wfl{k + 1 + 1 + 2^{\ell} u}
= k + 2 + 2^{\ell} u = k + 2 = \wlr{k+1} + 1,
\]
because $k + 2 \geq 2^{\ell}$ and we break ties downward.
Therefore, $\wfl{\sum_{i = 0}^k y_i} = k + 1$.

Let us now assume that the first Equation in \pRef{exqgI} holds for
$k$ and show that it holds for $k+1$.
When $2^{\ell} - 1 \leq k < 2^{\ell + 1} - 2$ we have that
$2^{\ell} - 1 \leq k + 1 < 2^{\ell + 1} - 1$ and
\pbDef{exqgIA}
\sum_{i = 0}^{k+1} y_i =
\wlr{\sum_{i = 0}^{k} y_i} + y_{k+1} =
k + 1 + \frac{4^{\ell} + 2}{3} u  + \wlr{k + 1 - 2^{\ell}} 2^{\ell} u + 1 + 2^{\ell} u
\peDef{exqgIA}
\[
= \wlr{k + 1} + 1 + \frac{4^{\ell + 1} + 2}{3} u  + \wlr{\wlr{k + 1} + 1 - 2^{\ell+1}} 2^{\ell+1} u
\]
and the first equality in Equation \pRef{exqgI} holds for $k + 1$. For $k = 2^{\ell + 1} - 2$, we have that
\[
2^{\ell + 1} - 1 = k + 1 < 2^{\wlr{\ell + 1} +1} - 1
\]
and Equation \pRef{exqgIA} leads to
\[
\sum_{i = 0}^{k+1} y_i = k + 2 + \frac{4^{\ell} + 2}{3} u  + 4^{\ell} u
= \wlr{k + 1} + 1 + \frac{4^{\ell + 1} + 2}{3} u + \wlr{\wlr{k + 1} + 1 - 2^{\ell}} 2^{\ell + 1} u
\]
because $k + 2 - 2^{\ell + 1} = 0$,
and the first equality in Equation \pRef{exqgI} is satisfied for $k + 1$.

Finally, for $n = 2^m - 1$ we have that $\ell = m$ and
\[
\frac{1}{u} \wlr{\sum_{k = 0}^n y_k - \wfl{\sum_{k = 0}^n y_k}} =
\frac{4^{m} + 2}{3}  = \frac{\wlr{n + 1}^2 + 2}{3} = \frac{n^2 + 2n + 3}{3}.
\]
The last equation in Example \ref{exQuadGrowth} follows from the equation above and
the fact that $y_k < 2$ when $2^m u < 1$.
\peVerify{exQuadGrowth}\\

%%%%%%%%%%%%%%%%%%%%%%%%%%%%%%%%%%%%%%%%%%%%%%%%%%%%%%%%%%%%%%%%%%%%%%%%%
% End of the verification of example exQuadGrowth
%%%%%%%%%%%%%%%%%%%%%%%%%%%%%%%%%%%%%%%%%%%%%%%%%%%%%%%%%%%%%%%%%%%%%%%%%
%
%
%
%
%%%%%%%%%%%%%%%%%%%%%%%%%%%%%%%%%%%%%%%%%%%%%%%%%%%%%%%%%%%%%%%%%%%%%%%%%
% Verification of example exSharpSumNearB
%%%%%%%%%%%%%%%%%%%%%%%%%%%%%%%%%%%%%%%%%%%%%%%%%%%%%%%%%%%%%%%%%%%%%%%%%

\pbVerifyB{exSharpSumNearB} Let us define $\rho := u^{-k}$.
Since $x_k = \rho^k$ and we break ties downward, we have
$\wfl{\sum_{i = 0}^k x_i} = \rho^k$ and
\[
\sum_{i = 0}^k x_i = \frac{\rho^{k+1} - 1}{\rho -1}
\hspace{0.2cm} \wrm{and} \hspace{0.2cm}
\sum_{k=1}^n \sum_{i = 0}^k x_i = \frac{\rho^{n+2} - \rho^2 - n \wlr{\rho - 1}}{\wlr{\rho -1}^2}
= \frac{1}{u^n} \frac{1 - u^n - n u^{n+1} \wlr{1 - u}}{\wlr{1 - u}^2}.
\]
It follows that
\[
\sum_{k = 0}^n x_k - \wfl{\sum_{k = 0}^n x_k} = \frac{\rho^{n+1} - 1}{\rho - 1} - \rho^n =
\frac{\rho^n - 1}{\rho - 1} = \frac{1}{u^{n-1}} \frac{1 - u^n}{1 - u} = \kappa_n u \sum_{k=1}^n \sum_{i = 0}^k x_i
\]
for
\[
\kappa_n := \frac{\wlr{1 - u}\wlr{{1 - u^{n}}}}{1 - u^{n} - n u^{n+1} \wlr{1 - u}}.
\]
If $2 n u < 1$ then $1 - u \leq \kappa_n \leq \wlr{1 - u} \wlr{1 + u^n}$ because
\[
 0 < \frac{\kappa_n}{1 - u} - 1 =
 \frac{n u^{n+1} \wlr{1 - u}}{1 - u^{n} - n u^{n+1} \wlr{1 - u}}
 = u^n \frac{n u \wlr{1 - u}}{1 - u^n - n u^{n+1} \wlr{1 - u}} < u^n.
\]
\peVerify{exSharpSumNearB} \\
%%%%%%%%%%%%%%%%%%%%%%%%%%%%%%%%%%%%%%%%%%%%%%%%%%%%%%%%%%%%%%%%%%%%%%%%%
% End of the verification of example exSharpSumNearB
%%%%%%%%%%%%%%%%%%%%%%%%%%%%%%%%%%%%%%%%%%%%%%%%%%%%%%%%%%%%%%%%%%%%%%%%%
%
%

%
%
%%%%%%%%%%%%%%%%%%%%%%%%%%%%%%%%%%%%%%%%%%%%%%%%%%%%%%%%%%%%%%%%%%%%%%%%%
% Verification of example exSharpSumNearC
%%%%%%%%%%%%%%%%%%%%%%%%%%%%%%%%%%%%%%%%%%%%%%%%%%%%%%%%%%%%%%%%%%%%%%%%%

\pbVerifyB{exSharpSumNearC}
Recall that $x_0 := u$, $x_1 := 1$ and
\[
x_k := \beta^{e_k} \wlr{1 + u} - \beta^{e_{k-1}} \wlr{1 + 2 u}
\]
for $k \geq 2$, with $0 = e_1 < e_2 \dots < e_n \in \wz{}$.
Induction using the basic properties of rounding to nearest in Prop. \ref{propRoundNormal} shows that
\[
s_k := \sum_{i = 0}^k x_i = \beta^{e_k} \wlr{1 + u} - u \sum_{i=1}^{k-1} \beta^{e_i},
\hspace{1.5cm}
\hat{s}_{k} := \wfl{\sum_{i = 0}^k x_i} = \beta^{e_k} \wlr{1 + 2u}
\]
for $k \geq 1$  and
\[
\sum_{k = 1}^n s_k = \wlr{1 + u} \sum_{k = 1}^n \beta^{e_k} - u \sum_{k =1}^n \sum_{i = 1}^{k - 1} \beta^{e_i}
= \wlr{1 + u} \sigma_{n} - u \sum_{k = 1}^n \sigma_{k - 1},
\]
for $\sigma_k := \sum_{i = 1}^{k} \beta^{e_i}$ (we assume that $\sum_{i}^k a_k = 0$ when $k < i$.)
Therefore
\[
\hat{s}_{n} - s_{n} = \wlr{\beta^{e_n} + \sum_{k = 1}^{n-1} \beta^{e_k}} u
= u \sum_{k = 1}^{n} \beta^{e_k} = u \sigma_n,
\]
and
\pbDef{exncA}
\frac{\hat{s}_{n} - s_{n}}{\sum_{k=1}^n s_k} =
\frac{\sigma_n u}{\sigma_n + u \wlr{\sigma_n - \sum_{k=0}^{n-1} \sigma_k}}
= \frac{u}{1 + u \wlr{1 - \sum_{k = 1}^{n - 1} v_k }}
\hspace{0.5cm} \wrm{for} \hspace{0.5cm}
v_k := \frac{\sigma_k}{\sigma_n}.
\peDef{exncA}

Note that
\[
\sigma_{\wlr{k+1}} - 1 = \sum_{i = 1}^{k + 1} \beta^{e_i} - 1 =
 \sum_{i = 2}^{k + 1} \beta^{e_i}
\geq  \beta \sum_{i = 2}^{k + 1} \beta^{e_{\wlr{i-1}}}
=  \beta \sum_{i = 1}^{k} \beta^{e_{i}}
= \beta \sigma_k.
\]
Since $\sigma_0 = 0$ and $1/\sigma_n = v_1 = \sigma_1 / \sigma_n$,
dividing the last equation by $\sigma_n$ we obtain
\pbDef{excnLPA}
v_1 + \beta v_k -  v_{k+1} \leq 0 \hspace{1cm} \wrm{for} \ \ k = 1,\dots, n - 2,
\hspace{0.5cm} \wrm{and} \hspace{0.5cm}
v_1 + \beta v_{n-1} \leq 1.
\peDef{excnLPA}
We end the verification of Example \ref{exSharpSumNearC}
using a duality argument to prove that
\pbDef{lpineq}
\sum_{k  = 1}^{n - 1} v_k \leq \frac{1}{\beta - 1} - \frac{n}{\beta^n - 1}.
\peDef{lpineq}
This equation combined with Equation \pRef{exncA} shows that the value of $\tau_n$ mentioned in
Example \ref{exSharpSumNearC} is appropriate.
We use basic facts about duality in linear programming \cite{Chvatal}
    applied to the problem with variables $v_k$,  objective function $\sum_{k = 1}^{n-1} v_k$ and
constraints given by $v_k \geq 0$ and  Equation \pRef{excnLPA}.
This problem can be written as
\pbDef{excnLP}
\left\{
\begin{array}{cccccccc}
\wrm{maximize}     & \wones{}^{\wtr} \wvec{v}  & =      & \sum_{k = 1}^{n-1} v_k  &   &     &        & \\
\wrm{subject \ to} & \wvec{A} \wvec{v}         & \leq   & \wvec{e},               &   & v_k & \geq   & 0.
\end{array}
\right.
\peDef{excnLP}
where the matrix $A$ has $a_{11} := \beta + 1$, $a_{i1} = 1$ for $1 < i < n$,
$a_{ii} = \beta$ for $2 \leq i < n$, $a_{i,i+1} = -1$ for $1 \leq i < n-1$ and
the remaining $a_{ij}$ are $0$. The vector $\wones$  has all its entries equal to $1$
and $e_i = 0$ for $1 \leq i < n - 1$ and $e_{n-1} = 1$.

This problem has a feasible solution
\[
v_k = \frac{\beta^{k} - 1}{\beta^{n}- 1},
\hspace{1cm} \wrm{for} \hspace{1cm} k = 1,\dots,n-1
\]
and
\pbDef{excnVal}
\sum_{k = 0}^{n-1} v_k = \frac{1}{\beta^n - 1} \wlr{\frac{\beta^{n} - 1}{\beta - 1} - n}
= \frac{1}{\beta - 1} - \frac{n}{\beta^n - 1}.
\peDef{excnVal}
Its dual has $n - 1$ variables, which we call $y_1, \dots, y_{n-1}$, and is
\pbDef{excnDual}
\left\{
\begin{array}{cccccccc}
\wrm{minimize}     & \wvec{e}^{\wtr} \wvec{y} & =      & y_{n-1}  &      &     &        & \\
\wrm{subject \ to} & \wvec{A}^{\wtr} \wvec{y} & \geq   & \wones{},& \ \  & y_k & \geq   & 0.
\end{array}
\right.
\peDef{excnDual}
We claim that the vector $\wvec{y} \in \wrn{n-1}$ with entries
\[
y_{n-1} = \frac{1}{\beta - 1} - \frac{n}{\beta^n - 1}
\hspace{1cm} \wrm{and} \hspace{1cm}
y_{k} = \beta^{n - k - 1} y_{n-1} + \frac{1}{\beta - 1} \ \ \ \wrm{for} \ k = 1 \dots n - 2,
\]
is a feasible solution of the dual problem.
Indeed, $y_{n-1} \geq 0$ because
\[
\frac{\beta^n - 1}{\beta - 1} = \sum_{k = 0}^{n-1} \beta^k \geq n,
\]
and the other entries of $\wvec{y}$ are clearly non negative because $y_{n-1} \geq 0$.
The first inequality in the system $\wvec{A}^{\wtr} \wvec{y} \geq \wones{}$
is satisfied because
\[
\wlr{\beta + 1} y_1 + \sum_{k = 2}^{n-1} y_k =
 \wlr{\wlr{\beta + 1} \beta^{n-2} + \sum_{k=2}^{n-2} \beta^{n-k-1} + 1} y_{n-1} + \frac{\beta + n-2}{\beta - 1}
\]
\[
= \wlr{\sum_{k=0}^{n-1} \beta^k} y_{n-1} + \frac{\beta + n-2}{\beta - 1}
\]
\[
= \frac{\beta^{n} - 1}{\beta - 1} \wlr{\frac{1}{\beta - 1} - \frac{n}{\beta^n - 1}} + \frac{\beta + n-2}{\beta - 1}
= \frac{\beta^{n} - \beta}{\wlr{\beta - 1}^2} + 1 \geq 1,
\]
and the remaining inequalities are satisfied as equalities, because
\[
- y_{k-1} + \beta y_k =
-\beta^{n - k} y_{n-1} - \frac{1}{\beta - 1} +
 \beta \beta^{n - k - 1} y_{n-1} + \frac{\beta}{\beta - 1} = 1.
\]

The value of the objective function of the dual problem for $\wvec{y}$, $y_{n-1}$,
is equal to the value of the objective function of the primal problem in \pRef{excnVal}.
Therefore, this is the optimal value of both problems
and Equation \pRef{lpineq} holds.
The linear programming problem
above also shows that the worst case in Equation \pRef{ssBnC} is achieved
for $e_k = k - 1$, because these exponents lead to the $v_k$ in the solution
of the primal problem.
\peVerify{exSharpSumNearC}\\

%%%%%%%%%%%%%%%%%%%%%%%%%%%%%%%%%%%%%%%%%%%%%%%%%%%%%%%%%%%%%%%%%%%%%%%%%
% End of the verification of example exSharpSumNearC
%%%%%%%%%%%%%%%%%%%%%%%%%%%%%%%%%%%%%%%%%%%%%%%%%%%%%%%%%%%%%%%%%%%%%%%%%
%
%

%
%
%%%%%%%%%%%%%%%%%%%%%%%%%%%%%%%%%%%%%%%%%%%%%%%%%%%%%%%%%%%%%%%%%%%%%%%%%%
% Verification of example exSignedSumB
%%%%%%%%%%%%%%%%%%%%%%%%%%%%%%%%%%%%%%%%%%%%%%%%%%%%%%%%%%%%%%%%%%%%%%%%%
\pbVerifyB{exSignedSumB}
Recall that $x_0 := u$, $x_1 := 1$ and $x_k := - 2^{1 - k} \wlr{1 + 3 u}$
for $k > 1$. It follows by induction that
\[
\sum_{i = 0}^k x_i = 2^{1 - k} \wlr{1 + 3 u} - 2 u
\hspace{1.0cm} \wrm{and} \hspace{1.0cm}
\wfl{\sum_{i=0}^k x_i} = 2^{1 - k} \wlr{1 + 2 u}.
\]
Since $2^{n} u \leq 1$, we have
\[
\sum_{k = 1}^n \wabs{\sum_{i = 0}^k x_i} = 2 \wlr{1 -  2^{-n}} \wlr{1 + 3u} - 2 n u > 0,
\]
and Equations \pRef{thXSSA} follows from the expressions above.
Finally, since $2^{-n} \geq u$, we have that $n u < 1$ and
\[
\kappa_n - \wlr{1 - u} =
u \frac{\wlr{2^{-n} - u} n + 3 \wlr{1 - 2^{-n}}u }{\wlr{1 -  2^{-n}} \wlr{1 + 3u} - n u} > 0
\]
and
\[
1 - \kappa_n =
u \frac{1 - 2^{-n} \wlr{n + 1}}{\wlr{1 -  2^{-n}} \wlr{1 + 3u} - n u} \geq 0.
\]
\peVerify{exSignedSumB}\\

%%%%%%%%%%%%%%%%%%%%%%%%%%%%%%%%%%%%%%%%%%%%%%%%%%%%%%%%%%%%%%%%%%%%%%%%%%
% End of the verification of example exSignedSumB
%%%%%%%%%%%%%%%%%%%%%%%%%%%%%%%%%%%%%%%%%%%%%%%%%%%%%%%%%%%%%%%%%%%%%%%%%
%
%

%%%%%%%%%%%%%%%%%%%%%%%%%%%%%%%%%%%%%%%%%%%%%%%%%%%%%%%%%%%%%%%
% The examples are verified only in the full version
%%%%%%%%%%%%%%%%%%%%%%%%%%%%%%%%%%%%%%%%%%%%%%%%%%%%%%%%%%%%%%%
} %% latexFull
%%%%%%%%%%%%%%%%%%%%%%%%%%%%%%%%%%%%%%%%%%%%%%%%%%%%%%%%%%%%%%%
% The examples are verified only in the full version
%%%%%%%%%%%%%%%%%%%%%%%%%%%%%%%%%%%%%%%%%%%%%%%%%%%%%%%%%%%%%%%

\end{document}